\documentclass[11pt]{article}
\usepackage{}
\usepackage{amssymb}
\usepackage[hang,flushmargin]{footmisc}
\usepackage{multirow}
\usepackage{amsfonts}
\usepackage{amsmath}
\usepackage{graphicx}
\usepackage{amsmath,amssymb,latexsym,color}
\usepackage[mathscr]{eucal}
\usepackage{cite}
 \makeatletter
    
    \newcommand{\Rmnum}[1]
    {\expandafter\@slowromancap\romannumeral #1@}
    \makeatother
\def\wz{\tilde}

\newtheorem{thm}{Theorem}[section]

\newtheorem{prop}[thm]{Proposition}
\newtheorem{lemma}[thm]{Lemma}
\newcounter{foo}[subsection]
\newcounter{fooo}[section]
\newtheorem{step}[foo]{Step}
\newtheorem{stepp}[fooo]{Step}

\newtheorem{cor}[thm]{Corollary}

\newtheorem{example}{Example}[section]
\newtheorem{defin}[thm]{Definition}

\newtheorem{remark}[thm]{Remark}

\newcommand{\qed}{\hfill\Box\medskip}

\usepackage{CJK}
\begin{document}
\begin{CJK*}{GBK}{song}

\renewcommand{\baselinestretch}{1.3}
\title{Thick weakly distance-regular digraphs}

\author{
Yuefeng Yang\textsuperscript{a}\quad
Kaishun Wang\textsuperscript{b}\\
\\{\footnotesize  \textsuperscript{a} \em  School of Science, China University of Geosciences, Beijing, 100083, China}\\{\footnotesize  \textsuperscript{b} \em Sch. Math. Sci. {\rm \&} Lab. Math. Com. Sys., Beijing Normal University, Beijing, 100875, China }  }
\date{}
\maketitle
\footnote{\scriptsize
{\em E-mail address:} yangyf@cugb.edu.cn(Y. Yang), wangks@bnu.edu.cn(K. Wang).}

\begin{abstract}

A weakly distance-regular digraph is thick if its attached scheme is regular. In this paper, we show that each commutative thick weakly distance-regular digraph has a thick weakly distance-regular subdigraph such that the corresponding quotient digraph falls into six families of thick weakly distance-regular digraphs up to isomorphism.

\medskip
\noindent {\em AMS classification:} 05E30

\noindent {\em Key words:} Weakly distance-regular digraph; thick; association scheme

\end{abstract}

\section{Introduction}

All the digraphs considered in this paper are finite, simple and strongly connected. Let $\Gamma$ be a digraph and $V\Gamma$ be its vertex set. For any $x,y\in V\Gamma$, let $\partial_\Gamma(x,y)$ denote the \emph{distance} from $x$ to $y$ in $\Gamma$. The pair $\wz{\partial}_{\Gamma}(x,y)=(\partial_{\Gamma}(x,y),\partial_{\Gamma}(y,x))$ is called the \emph{two-way distance} from $x$ to $y$. If no confusion occurs, we write $\partial(x,y)$ (resp. $\wz{\partial}(x,y)$) instead of $\partial_\Gamma(x,y)$ (resp. $\wz{\partial}_\Gamma(x,y)$). Let $\wz{\partial}(\Gamma)$ be the set of all pairs $\wz{\partial}(x,y)$.

A digraph $\Gamma$ is said to be \emph{weakly distance-regular} if, for any $\wz{h}$, $\wz{i}$, $\wz{j}\in\wz{\partial}(\Gamma)$, the number of $z\in V\Gamma$ such that $\wz{\partial}(x,z)=\wz{i}$ and $\wz{\partial}(z,y)=\wz{j}$ is constant whenever $\wz{\partial}(x,y)=\wz{h}$. This constant is denoted by $p_{\wz{i},\wz{j}}^{\wz{h}}$. The integers
$p_{\wz{i},\wz{j}}^{\wz{h}}$ are called the \emph{intersection numbers}. We say
that $\Gamma$ is \emph{commutative} if $p_{\tilde{i},\tilde{j}}^{\tilde{h}}=p_{\tilde{j},\tilde{i}}^{\tilde{h}}$ for all $\tilde{i}$, $\tilde{j}$,
$\tilde{h}\in\tilde{\partial}(\Gamma)$. $\Gamma$ is \emph{thin} (resp. \emph{quasi-thin}) if the maximum value of its intersection numbers is $1$ (resp. $2$).

As a natural generalization of distance-regular graphs (see \cite{AEB98, DKT16} for the theory of distance-regular graphs), Wang and Suzuki \cite{KSW03} introduced the concept of weakly distance-regular digraphs. Since then some special families of weakly distance-regular digraphs were classified. See \cite{KSW03,HS04} for valency $2$, \cite{KSW04,YYF16,YYF18} for valency $3$, \cite{HS04} for thin case, \cite{YYF16+} for quasi-thin case.

To state our main theorem, we need additional notations and terminologies.

Let $\Gamma$ be a weakly distance-regular digraph and $R=\{\Gamma_{\wz{i}}\mid\wz{i}\in\wz{\partial}(\Gamma)\}$, where $\Gamma_{\wz{i}}=\{(x,y)\in V\Gamma\times V\Gamma\mid\wz{\partial}(x,y)=\wz{i}\}$. Then $(V\Gamma,R)$ is an association scheme (see \cite{EB84,PHZ96,PHZ05} for the theory of association schemes), which is called
the \emph{attached scheme} of $\Gamma$. The size of $\Gamma_{\tilde{i}}(x):=\{y\in V\Gamma\mid\tilde{\partial}(x,y)=\tilde{i}\}$
depends only on $\tilde{i}$, denoted by $k_{\tilde{i}}$. For each $\tilde{i}:=(a,b)\in\tilde{\partial}(\Gamma)$, we define $\tilde{i}^{*}=(b,a)$, and write $k_{a,b}$ (resp. $\Gamma_{a,b}$) instead of $k_{(a,b)}$ (resp. $\Gamma_{(a,b)}$).

For two nonempty subsets $E$ and $F$ of $R$, define
\begin{eqnarray*}
EF&:=&\{\Gamma_{\tilde{h}}\mid\sum_{\Gamma_{\tilde{i}}\in E}\sum_{\Gamma_{\tilde{j}}\in F}p_{\tilde{i},\tilde{j}}^{\tilde{h}}\neq0\},
\end{eqnarray*}
and write $\Gamma_{\tilde{i}}\Gamma_{\tilde{j}}$ instead of $\{\Gamma_{\tilde{i}}\}\{\Gamma_{\tilde{j}}\}$. We say $F$ {\em closed} if $\Gamma_{\wz{i}^{*}}\Gamma_{\wz{j}}\subseteq F$ for any $\Gamma_{\wz{i}}$ and $\Gamma_{\wz{j}}$ in $F$. Let $\langle F\rangle$ be the minimum closed subset containing $F$. Denote $F(x)=\{y\in V\Gamma\mid(x,y)\in \cup_{f\in F}f\}$ and $T=\{q\mid(1,q-1)\in\wz{\partial}(\Gamma)\}$. For a subset $I$ of $T$ and a vertex $x\in V\Gamma$, let $\Delta_{I}(x)$ be the digraph $(F_{I}(x),\cup_{q\in I}\Gamma_{1,q-1})$, where $F_{I}=\langle\{\Gamma_{1,q-1}\}_{q\in I}\rangle$. Since the digraph $\Delta_{I}(x)$ does not depend on the choice of vertex $x$ up to isomorphism, if
no confusion occurs, we write $\Delta_{I}$ instead of $\Delta_{I}(x)$.

For any nonempty closed subset $F$ of $R$, let $V\Gamma/F:=\{F(x)\mid x\in V\Gamma\}$.  The \emph{quotient digraph} of $\Gamma$ over $F$, denoted by $\Gamma/F$, is defined as the digraph with vertex set $V\Gamma/F$ in which $(F(x),F(y))$ is an arc whenever there is an arc in $\Gamma$ from $F(x)$ to $F(y)$.

A weakly distance-regular digraph $\Gamma$ is \emph{thick} if its attached scheme is regular, that means $p_{\wz{i},\wz{i}}^{\wz{h}}, p_{\wz{i},\wz{i}^{*}}^{\wz{h}}\in\{0,k_{\wz{i}}\}$ for any $\wz{i},\wz{h}\in\wz{\partial}(\Gamma)$. The regular association schemes have been investigated
by Yoshikawa in \cite{MY15,MY16,MY18,MY20}. In this paper, we study commutative thick weakly distance-regular digraphs, and obtain the following result.

\begin{thm}\label{Main4}
Let $\Gamma$ be a commutative thick weakly distance-regular digraph, and $q$ the maximum integer in $T$. Then $\Delta_{T\setminus I}$ is a thick weakly distance-regular digraph for some nonempty subset $I$ of $\{q-1,q\}$, and $\Gamma/F_{T\setminus I}$ is isomorphic to one of the following digraphs:\vspace{-0.3cm}
\begin{itemize}
\item [{\rm(i)}] ${\rm Cay}(\mathbb{Z}_{p},\{1\})$;\vspace{-0.3cm}

\item [{\rm(ii)}] ${\rm Cay}(\mathbb{Z}_{p},\{1\})[\overline{K}_{2}]$;\vspace{-0.3cm}

\item [{\rm(iii)}] ${\rm Cay}(\mathbb{Z}_{2q-2},\{1,2\})$;\vspace{-0.3cm}

\item [{\rm(iv)}] ${\rm Cay}(\mathbb{Z}_{2q-2},\{1,2\})[\overline{K}_{2}]$;\vspace{-0.3cm}

\item [{\rm(v)}] ${\rm Cay}(\mathbb{Z}_{2^{\alpha+1}(q-1)}\times\mathbb{Z}_{2^{1-\alpha}},\{(2^{\alpha}+\beta,1),(2^{\alpha}-\beta,\alpha),(2^{\alpha+1},\alpha+1)\})$;\vspace{-0.3cm}

\item [{\rm(vi)}] ${\rm Cay}(\mathbb{Z}_{2^{\alpha+1}(q-1)}\times\mathbb{Z}_{2^{1-\alpha}},\{(2^{\alpha}+\beta,1),(2^{\alpha}-\beta,\alpha),(2^{\alpha+1},\alpha+1)\})[\overline{K}_{2}]$.\vspace{-0.3cm}
\end{itemize}
Here, $p\in\{q-1,q\}$, $\alpha=2(\frac{2}{{\rm gcd}(q-1,4)}-\lfloor\frac{2}{{\rm gcd}(q-1,4)}\rfloor)-1$, $\beta=\frac{2^{\alpha}\alpha(1-q)}{{\rm gcd}(q-1,4)}$ and $\overline{K}_{2}$ is the coclique of $2$ vertices.
\end{thm}

Routinely, all digraphs in Theorem \ref{Main4} (i)--(vi) are thick weakly distance-regular digraphs. Since all the digraphs in (i)--(iii) and (v) have valency at most $3$, one can check from \cite[Theorem 1.1]{KSW03}, \cite[Theorem 1.3]{KSW04} and \cite[Theorem 1 and Proposition 9]{YYF16}. Since the digraphs in (iv) and (vi) are lexicographic products from the digraphs in (iii) and (v) to $\overline{K}_{2}$, respectively, they are thick weakly distance-regular digraphs from \cite[Section 2]{YYF16} and \cite[Propositions 2.3 and 2.4]{KSW03}.

This paper is organized as follows. In Section 2, we introduce some basic results about commutative thick weakly distance-regular digraphs. In Section 3, applying the results in Section 2, we give a characterization of mixed arcs. In Section 4, applying the results in Section 3, we give some results about the two-way distance which are used frequently in this paper. In Section 5, applying the results in Sections 3 and 4, we determines the relationship between different types of arcs. In Section 6, applying the results in Sections 3--5, we give a proof of Theorem \ref{Main4}.

\section{Basic properties}

In the remaining of this paper, $\Gamma$ always denotes a commutative thick weakly distance-regular digraph, let $P_{\wz{i},\wz{j}}(x,y)=\Gamma_{\wz{i}}(x)\cap\Gamma_{\wz{j}^{*}}(y)$ for all $\wz{i},\wz{j}\in\wz{\partial}(\Gamma)$ and $x,y\in V\Gamma$.

\begin{lemma}\label{jiben2}
If $\Gamma_{\wz{i}}^{l}\Gamma_{\wz{j}}\cap\Gamma_{\wz{i}}^{l-1}\Gamma_{\wz{j}}^{2}\neq\emptyset$ with $l>1$, then $\wz{i}=\wz{j}$.
\end{lemma}
\noindent\textbf{Proof.}~By the assumption, we may assume $(x_{l},x_{l+1}),(y_{l-1},y_{l}),(y_{l},y_{l+1})\in\Gamma_{\wz{j}}$ and $(x_{h},x_{h+1}),(y_{h'},y_{h'+1})\in\Gamma_{\wz{i}}$ for $0\leq h\leq l-1$ and $0\leq h'\leq l-2$, where $x_{0}=y_{0}$ and $x_{l+1}=y_{l+1}$. Since $P_{\wz{i},\wz{i}}(y_{0},x_{2})=\Gamma_{\wz{i}}(y_{0})$, we have $(y_{1},x_{2})\in\Gamma_{\wz{i}}$. By induction, one gets $(y_{h},x_{h+1})\in\Gamma_{\wz{i}}$ for any $0\leq h\leq l-1$. By $P_{\wz{j},\wz{j}}(y_{l-1},y_{l+1})=\Gamma_{\wz{j}^{*}}(y_{l+1})$, we get $\wz{\partial}(y_{l-1},x_{l})=\wz{i}=\wz{j}$.$\qed$

\begin{lemma}\label{jiben3}
If $\Gamma_{1,p-1}\in\Gamma_{1,q-1}\Gamma_{q-1,1}$, then $\Gamma_{1,p-1}^{2}\subseteq\Gamma_{1,q-1}\Gamma_{q-1,1}$.
\end{lemma}
\textbf{Proof.}~Since $p_{(1,q-1),(q-1,1)}^{(1,p-1)}=k_{1,q-1}$, we have $\Gamma_{1,p-1}\Gamma_{1,q-1}=\{\Gamma_{1,q-1}\}$. This implies $\Gamma_{1,p-1}^{2}\subseteq\Gamma_{1,p-1}\Gamma_{1,q-1}\Gamma_{q-1,1}\subseteq\Gamma_{1,q-1}\Gamma_{q-1,1}$.$\qed$

We recall the definitions of pure arcs and mixed arcs introduced in \cite{YYF16+}. An arc $(u,v)$ of $\Gamma$ is of \emph{type} $(1,q-1)$ if $\partial(v,u)=q-1$. An arc of type $(1,q-1)$ is said to be \emph{pure}, if every circuit of length $q$ containing it consists of arcs of type $(1,q-1)$; otherwise, this arc is said to be \emph{mixed}. We say that $(1,q-1)$ is pure if any arc of type $(1,q-1)$ is pure; otherwise, we say that $(1,q-1)$ is mixed. We say that the configuration C$(q)$ (resp. D$(q)$) exists if $p_{(1,q-1),(1,q-1)}^{(1,q-2)}\neq0$ (resp. $p_{(1,q-2),(q-2,1)}^{(1,q-1)}\neq0$) and $(1,q-2)$ is pure.

\begin{lemma}\label{jiben4}
If C$(q)$ exists, then any shortest path between distinct vertices contains at most two arcs of type $(1,q-1)$.
\end{lemma}
\textbf{Proof.}~Suppose to, the contrary that $(x_{0},x_{1},\ldots,x_{l})$ is a shortest path such that $(x_{0},x_{1}),(x_{1},x_{2}),(x_{2},x_{3})\in\Gamma_{1,q-1}$. Since C$(q)$ exists, there exists a vertex $x_{2}'\in P_{(1,q-2),(q-1,1)}(x_{0},x_{1})$. By $P_{(1,q-1),(1,q-1)}(x_{1},x_{3})=\Gamma_{1,q-1}(x_{1})$, we have $(x_{2}',x_{3})\in\Gamma_{1,q-1}$, contrary to the fact that $(x_{0},x_{1},\ldots,x_{l})$ is a shortest path. $\qed$

\begin{lemma}\label{(1,q-1)}
If $p_{(1,q-1),(1,q-1)}^{(2,q-2)}=m>0$, then $\Delta_{\{q\}}$ is isomorphic to the digraph ${\rm Cay}(\mathbb{Z}_{q}\times\mathbb{Z}_{m},\{(1,0),(1,1),\ldots,(1,m-1)\})$ and $\Gamma_{1,q-1}^{l}=\{\Gamma_{l,q-l}\}$ for $1\leq l\leq q-1$.
\end{lemma}
\textbf{Proof.}~Choose a circuit $(x_{0,0},x_{1,0},\ldots,x_{q-1,0})$ such that $(x_{0,0},x_{1,0}),(x_{1,0},x_{2,0})\in\Gamma_{1,q-1}$ and $(x_{0,0},x_{2,0})\in\Gamma_{2,q-2}$, where the first subscription of $x$ are taken modulo $q$. Since $m=k_{1,q-1}$, one has $\Gamma_{1,q-1}(x_{0,0})=\Gamma_{q-1,1}(x_{2,0})$. Note that $(x_{1,0},x_{3,0})\in\Gamma_{2,q-2}$. It follows that $\Gamma_{1,q-1}(x_{1,0})=\Gamma_{q-1,1}(x_{3,0})$ and $(x_{2,0},x_{3,0})\in\Gamma_{1,q-1}$. By induction, we get $(x_{i,0},x_{i+2,0})\in\Gamma_{2,q-2}$, $\Gamma_{1,q-1}(x_{i,0})=\Gamma_{q-1,1}(x_{i+2,0})$ and $(x_{i+1},x_{i+2})\in\Gamma_{1,q-1}$ for each $i$. Write $\Gamma_{1,q-1}(x_{i,0})=\{x_{i+1,j}\mid0\leq j\leq m-1\}$. Since $\wz{\partial}(x_{i-1,0},x_{i+1,j})=(2,q-2)$ for any $j$, we get $\wz{\partial}(x_{i,j'},x_{i+1,j})=(1,q-1)$ for any $j'$. The first statement is valid. By $\wz{\partial}(x_{i,j},x_{i+l,j'})=(l,q-l)$ for $1\leq l\leq q-1$, we obtain $\Gamma_{1,q-1}^{l}=\{\Gamma_{l,q-l}\}$.$\qed$

\begin{lemma}\label{mix-(2,q-2) 1}
If {\rm C}$(q)$ exists, then the following hold:\vspace{-0.3cm}
\begin{itemize}
\item [{\rm(i)}] $\Gamma_{1,q-1}\Gamma_{1,q-2}^{i-1}=\{\Gamma_{i,q-i}\}$ for $1\leq i\leq q-1$.\vspace{-0.3cm}

\item [{\rm(ii)}] $\Gamma_{1,q-1}^{2}=\{\Gamma_{1,q-2}\}$ or $\Gamma_{1,q-1}^{2}=\{\Gamma_{1,q-2},\Gamma_{2,q-1}\}$.\vspace{-0.3cm}

\item [{\rm(iii)}] Any circuit of length $q$ containing an arc of type $(1,q-1)$ consists of arcs of types $(1,q-1)$ and $(1,q-2)$.
\end{itemize}
\end{lemma}
\noindent\textbf{Proof.}~(i) Pick a path $(x,x_{0},x_{1},\ldots,x_{i-1})$ such that $\wz{\partial}(x,x_{0})=(1,q-1)$ and $\wz{\partial}(x_{j},x_{j+1})=(1,q-2)$ for $0\leq j\leq i-2$. Since $(1,q-2)$ is pure, from Lemma \ref{(1,q-1)}, we have $(x_{0},x_{i-1})\in\Gamma_{i-1,q-i}$. If $i=q-1$, from $p_{(1,q-1),(1,q-1)}^{(1,q-2)}=k_{1,q-1}$, then $x\in P_{(1,q-1),(1,q-1)}(x_{q-2},x_{0})$, as desired. Suppose that $i<q-1$. Pick a vertex $x_{q-2}\in P_{(q-1,1),(1,q-2)}(x,x_{0})$. In view of Lemma \ref{(1,q-1)}, we obtain $\wz{\partial}(x_{q-2},x_{i-1})=(i,q-i-1)$, which implies $\wz{\partial}(x,x_{i-1})=(i,q-i)$. Thus, (i) is valid.

(ii) Suppose $|\Gamma_{1,q-1}^{2}|>1$. Let $(z_{0},z_{1},z_{2})$ be a path consisting of arcs of type $(1,q-1)$ with $(z_{0},z_{2})\notin\Gamma_{1,q-2}$. Pick a vertex $z_{3}\in P_{(1,q-2),(q-1,1)}(z_{1},z_{2})$. By (i), one has $\Gamma_{1,q-1}\Gamma_{1,q-2}=\{\Gamma_{2,q-2}\}$, which implies $\wz{\partial}(z_{0},z_{3})=(2,q-2)$ and $p_{(1,q-1),(1,q-1)}^{(2,q-2)}=0$. Since $(z_{0},z_{2})\notin\Gamma_{1,q-2}$, we get $\partial(z_{2},z_{0})=q-1$. By $P_{(1,q-1),(1,q-1)}(z_{0},z_{2})=\Gamma_{1,q-1}(z_{0})$, one has $\wz{\partial}(z_{0},z_{2})=(2,q-1)$. This proves (ii).

(iii) Let $(y_{0},y_{1},\ldots,y_{q-1})$ be a circuit with $\wz{\partial}(y_{q-1},y_{0})=(1,q-1)$. If $(y_{i},y_{i+1})\in\Gamma_{1,q-1}$ for $0\leq i\leq q-2$, by $q>2$, then $(y_{0},y_{2})\in\Gamma_{2,q-2}$, and so $\Gamma_{2,q-2}\in\Gamma_{1,q-1}^2$, contrary to (ii). Hence, $(y_{i},y_{i+1})\in\Gamma_{1,p-1}$ with $p\neq q$ for some $i$. By the commutativity of $\Gamma$, we may assume $i=0$. It suffices to show that $p=q-1$. Since $p_{(1,q-1),(1,q-1)}^{(1,q-2)}=k_{1,q-1}$, one has $\wz{\partial}(y_{q-1},y_{1})=(2,q-2)$.  By (i), there exists a vertex $y\in P_{(1,q-2),(1,q-1)}(y_{q-1},y_{1})$, which implies $y_{0}\in P_{(1,q-1),(1,q-1)}(y_{q-1},y)$. Since $y\in P_{(1,q-1),(1,q-1)}(y_{0},y_{1})$, from (ii), we get $p=q-1$.$\qed$

The commutativity of $\Gamma$ will be used frequently in the sequel, so we no longer
refer to it for the sake of simplicity.

\begin{lemma}\label{Gamma{1,q-1}=2}
Suppose that {\rm C}$(q)$ exists and $\Gamma_{1,q-1}^{2}=\{\Gamma_{1,q-2},\Gamma_{2,q-1}\}$. Then the following hold:\vspace{-0.3cm}
\begin{itemize}
\item [{\rm(i)}] If $2\leq i\leq q-1$, then  $\Gamma_{1,q-1}^{2i-1}=\Gamma_{1,q-2}^{i-1}\Gamma_{1,q-1}=\{\Gamma_{i,q-i}\}$.\vspace{-0.3cm}

\item [{\rm(ii)}] If $2\leq i\leq q-2$, then $\Gamma_{1,q-1}^{2i}=\Gamma_{1,q-2}^{i-1}\Gamma_{1,q-1}^{2}=\{\Gamma_{i,q-i-1},\Gamma_{i+1,q-i}\}$.
\end{itemize}
\end{lemma}
\textbf{Proof.}~Let $(x_{0},x_{1},\ldots,x_{j})$ be a path consisting of arcs of type $(1,q-1)$ for $3\leq j\leq 2q-3$. By induction, there exists a vertex $x_{2h}'\in P_{(1,q-1),(1,q-1)}(x_{2h-1},x_{2h+1})=\Gamma_{1,q-1}(x_{2h-1})$ such that $(x_{2h-2}',x_{2h}')\in\Gamma_{1,q-2}$ for $1\leq h\leq(j-1)/2$, where $x_{0}'=x_{0}$. If $j=2i-1$ for some $i\in\{2,3,\ldots,q-1\}$, from Lemma \ref{mix-(2,q-2) 1} (i), then (i) holds.

Assume that $j=2i$ for some $i\in\{2,3,\ldots,q-2\}$. If $\wz{\partial}(x_{2i-2}',x_{2i})=(1,q-2)$, from Lemma \ref{(1,q-1)}, then $\wz{\partial}(x_{0},x_{2i})=(i,q-i-1)$. Now suppose $\wz{\partial}(x_{2i-2}',x_{2i})=(2,q-1)$. Pick a path $(x_{2i-2}',x_{2i}',\ldots,x_{2q-2}'=x_{0})$ consisting of arcs of type $(1,q-2)$ and a vertex $x_{2i+1}\in\Gamma_{1,q-1}(x_{2i})$. Since $P_{(1,q-1),(1,q-1)}(x_{2i-1},x_{2i+1})=\Gamma_{1,q-1}(x_{2i-1})$ and $p_{(1,q-1),(1,q-1)}^{(1,q-2)}=k_{1,q-1}$, one obtains $(x_{2i-1},x_{2i}'),(x_{2i}',x_{2i+1}),(x_{2i+1},x_{2i+2}')\in\Gamma_{1,q-1}$. The fact that $P_{(1,q-2),(1,q-2)}(x_{2i-2}',x_{2i+2}')=\Gamma_{q-2,1}(x_{2i+2}')$ implies that $(x_{2i},x_{2i+2}')\in\Gamma_{2,q-1}$. Since $\partial(x_{0},x_{2i})\leq i+1$ and $\partial(x_{2i+2}',x_{2i})\leq\partial(x_{2i+2}',x_{0})+\partial(x_{0},x_{2i})\leq q-i-2+\partial(x_{0},x_{2i})$, we have $\partial(x_{0},x_{2i})=i+1$. By (i), we get $\wz{\partial}(x_{0},x_{2i+1})=(i+1,q-i-1)$. Since $q-1=\partial(x_{2i},x_{2i-2}')\leq\partial(x_{2i},x_{0})+i-1$, one has $q-i\leq\partial(x_{2i},x_{0})\leq1+\partial(x_{2i+1},x_{0})=q-i$. Thus, (ii) holds. $\qed$

\begin{lemma}\label{mix-(2,q-2) 2}
If {\rm D}$(q)$ exists, then $\Gamma_{1,q-1}^{2}=\{\Gamma_{2,q-2}\}$ or $\Gamma_{1,q-1}^{2}=\{\Gamma_{2,q-1}\}$.
\end{lemma}
\noindent\textbf{Proof.}~Pick a path $(x,y,z)$ consisting of arcs of type $(1,q-1)$ and a vertex $w\in\Gamma_{1,q-2}(x)$. Since $p_{(1,q-2),(q-2,1)}^{(1,q-1)}=k_{1,q-2}$, we get $y,z\in\Gamma_{q-2,1}(w)$. Note that $q-2\leq\partial(z,x)\leq1+\partial(w,x)=q-1$. Since $p_{(1,q-2),(1,q-2)}^{(1,q-2)}=p_{(1,q-1),(1,q-1)}^{(1,q-1)}=0$, one has $\partial(x,z)=2$. If $\Gamma_{2,q-2}\in\Gamma_{1,q-1}^{2}$, then $\Gamma_{1,q-1}^{2}=\{\Gamma_{2,q-2}\}$ from Lemma \ref{(1,q-1)}; if $\Gamma_{2,q-2}\notin\Gamma_{1,q-1}^{2}$, then $\wz{\partial}(x,z)=(2,q-1)$ and $\Gamma_{1,q-1}^{2}=\{\Gamma_{2,q-1}\}$.$\qed$

\section{Characterization of mixed arcs}

The main result of this section is the following important result which characterizes
mixed arcs.

\begin{thm}\label{Main1}
The configuration {\rm C}$(q)$ or {\rm D}$(q)$ exists if and only if $(1,q-1)$ is mixed.
\end{thm}

In order to prove Theorem \ref{Main1}, we need a lemma.

\begin{lemma}\label{lemma Main1}
Suppose that $(1,q-1)$ is pure or one of the configurations C($q$) and D($q$) exists. If $p_{(1,s-1),(1,t-1)}^{(1,q-1)}\neq0$ with $s\neq t$, then $q\in\{s,t\}$.
\end{lemma}
\textbf{Proof.}~Let $x_{0},x,x',x_{1}$ be distinct vertices such that $\wz{\partial}(x_{0},x_{1})=(1,q-1)$, $x\in P_{(1,s-1),(1,t-1)}(x_{0},x_{1})$ and $x'\in P_{(1,t-1),(1,s-1)}(x_{0},x_{1})$. Assume the contrary, namely, $q\notin\{s,t\}$. Note that $q>2$.

We claim that $\wz{\partial}(x,y)\neq\wz{\partial}(x',y)$ for all $y\in\Gamma_{1,q-1}(x_{1})$. Let $y\in\Gamma_{1,q-1}(x_{1})$. Since $p_{(1,s-1),(1,t-1)}^{(1,q-1)}\neq0$, we have $p_{(1,s-1),(1,s-1)}^{(1,q-1)}=0$, which implies $(x_{0},x_{1}')\notin\Gamma_{1,q-1}$ for all $x_{1}'\in\Gamma_{1,s-1}(x)$. The fact that $P_{(1,q-1),(1,q-1)}(x_{0},y)=\Gamma_{q-1,1}(y)$ implies $(x_{1}',y)\notin\Gamma_{1,q-1}$ for all $x_{1}'\in\Gamma_{1,s-1}(x)$. Hence, $P_{(1,s-1),(1,q-1)}(x,y)=\emptyset$. Since $x_{1}\in P_{(1,s-1),(1,q-1)}(x',y)$, the claim is valid.

Suppose that $(1,q-1)$ is pure. Pick a vertex $y\in P_{(2,q-2),(q-1,1)}(x_{0},x_{1})$. Since $q\notin\{s,t\}$, we have $\partial(y,x)=\partial(y,x')=q-1$, which implies $\partial(x,y)=\partial(x',y)=2$, contrary to the claim.

Suppose that C($q$) exists. Note that there exists $y\in P_{(1,q-2),(q-1,1)}(x_{0},x_{1})$. Since $x\in P_{(1,s-1),(1,t-1)}(x_{0},x_{1})$ and $p_{(1,q-1),(1,q-1)}^{(1,q-2)}=k_{1,q-1}$, we have $q-1\notin\{s,t\}$. By Lemma \ref{mix-(2,q-2) 1} (iii), one gets $\partial(y,x)=\partial(y,x')=q-1$. Since $P_{(1,q-1),(1,q-1)}(x_{0},y)=\Gamma_{q-1,1}(y)$, we obtain $x,x'\in\Gamma_{q-1,2}(y)$, contrary to the claim.

Suppose that D$(q)$ exists. Since $s\neq t$ and $p_{(1,q-2),(q-2,1)}^{(1,q-1)}=k_{1,q-2}$, we have $q-1\notin\{s,t\}$. Pick vertices $x_{2}\in P_{(1,q-2),(q-2,1)}(x_{0},x_{1})$ and $y\in P_{(1,q-1),(1,q-2)}(x_{1},x_{2})$.

If $\wz{\partial}(x,x_{2})=\wz{\partial}(x_{0},y)$, from $x_{1}\in P_{(1,t-1),(1,q-2)}(x,x_{2})$, then there exists $x_{1}''\in P_{(1,t-1),(1,q-2)}(x_{0},y)$, which implies $x_{0}\in P_{(1,q-2),(1,q-2)}(x_{1}'',x_{2})=\Gamma_{q-2,1}(x_{2})$, contrary to $t\neq q-1$. Then $\wz{\partial}(x,x_{2})\neq\wz{\partial}(x_{0},y)$. Similarly, $\wz{\partial}(x',x_{2})\neq\wz{\partial}(x_{0},y)$.

Since $(1,q-2)$ is pure and $\partial(x_{2},x)\leq1+\partial(x_{2},x_{0})$, one has $\partial(x_{2},x)=q-2$ or $q-1$. The fact $q-1\neq s$ and $p_{(q-2,1),(1,q-2)}^{(1,q-1)}=k_{1,q-2}$ imply $(x,x_{2})\in\Gamma_{1,q-2}\cup\Gamma_{2,q-2}\cup\Gamma_{2,q-1}$. By Lemma \ref{mix-(2,q-2) 2}, we obtain $(x_{0},y)\in\Gamma_{2,q-2}\cup\Gamma_{2,q-1}$. Since $\wz{\partial}(x,x_{2})\neq\wz{\partial}(x_{0},y)$, one has $q-2\leq\partial(y,x)\leq\min\{\partial(x_{2},x),\partial(y,x_{0})\}+1=q-1$.  By $P_{(1,q-1),(1,q-1)}(x_{0},y)=\Gamma_{q-1,1}(y)$, we get $\wz{\partial}(x,y)\neq\wz{\partial}(x_{0},y)$ and $(x,y)\notin\Gamma_{1,q-1}$. Since $t\neq q-1$ and $p_{(1,q-2),(q-2,1)}^{(1,q-1)}=k_{1,q-2}$, we obtain $(x,y)\notin\Gamma_{1,q-2}$, and so $(x,y)\in\Gamma_{2,q-2}\cup\Gamma_{2,q-1}$. Similarly, $\wz{\partial}(x',y)\neq\wz{\partial}(x_{0},y)$ and $(x',y)\in\Gamma_{2,q-2}\cup\Gamma_{2,q-1}$. It follows that $\wz{\partial}(x,y)=\wz{\partial}(x',y)$, contrary to the claim.$\qed$

Now we are ready to give a proof of Theorem \ref{Main1}.

\noindent\textbf{Proof of Theorem \ref{Main1}.}~If C($q$) or D($q$) exists, it is obvious that $(1,q-1)$ is mixed. We prove the converse. By way of contradiction, we may assume that $q$ is the minimum integer such that $(1,q-1)$ is mixed, and neither C($q$) nor D($q$) exists. Since $(1,1)$ is pure, $q>2$. Pick a circuit $(z_{0},z_{1},\ldots,z_{q-1})$ such that $(z_{q-1},z_{0})\in\Gamma_{1,q-1}$.

\textbf{Case 1}. $\Gamma_{1,q-1}\in\Gamma_{l-1,1}^{q-1}$ for some $l\in\{2,3,\ldots,q-1\}$.

Without loss of generality, we assume $(z_{i},z_{i+1})\in\Gamma_{1,l-1}$ for $0\leq i\leq q-2$. By the minimality of $q$, $(1,l-1)$ is pure or one of the configurations C$(l)$ and D$(l)$ exists.

Suppose that $(1,l-1)$ is pure. If $l=2$, from $P_{(1,1),(1,1)}(z_{0},z_{2})=\Gamma_{1,1}(z_{2})$, then $q=3$ and $p_{(1,1),(1,1)}^{(1,2)}\neq0$, contrary to the fact that D($q$) does not exist. Then $q>l>2$. By Lemma \ref{(1,q-1)}, we have $\wz{\partial}(z_{0},z_{l-1})=(l-1,1)$. Since D($q$) does not exist, $q>l+1$. Hence, $z_{0}\in P_{(1,l-1),(1,l-1)}(z_{l-1},z_{l+1})=\Gamma_{1,l-1}(z_{l-1})$, a contradiction.

Suppose that C$(l)$ or D$(l)$ exists. Since $q>l>2$, from Lemma \ref{jiben4}, D$(l)$ exists. Pick a vertex $z_{0}'\in\Gamma_{1,l-2}(z_{0})$. Since $p_{(1,l-2),(l-2,1)}^{(1,l-1)}=k_{1,l-2}$, we have $\wz{\partial}(z_{i},z_{0}')=(1,l-2)$ for $0\leq i\leq q-1$. Then $q-1=\partial(z_{0},z_{q-1})\leq1+\partial(z_{0}',z_{q-1})=l-1$, contrary to $l<q$.

\textbf{Case 2}. $\Gamma_{1,q-1}\notin\Gamma_{l-1,1}^{q-1}$ for any $l\in\{2,3,\ldots,q-1\}$.

Without loss of generality, we may assume $\wz{\partial}(z_{0},z_{1})=(1,p-1)$ with $p<q$.

\textbf{Case 2.1}. $\partial(z_{i+1},z_{i})\neq q-1$ for $1\leq i\leq q-2$.

Without loss of generality, we may assume $\partial(z_{q-1},z_{q-2})=s-1$ with $s\notin\{p,q\}$.

Suppose $\partial(z_{q-1},z_{1})=1$. By the minimality of $q$ and Lemma \ref{lemma Main1}, we have $\partial(z_{1},z_{q-1})=q-2=p-1$, and so $p=q-1$. Since $s\neq q-1$, $(1,q-2)$ is mixed. If D$(q-1)$ exists, from $p_{(1,q-3),(q-3,1)}^{(1,q-2)}=k_{1,q-3}$, then there exists a vertex $z_{q-2}'$ such that $z_{q-1},z_{0},z_{1}\in\Gamma_{q-3,1}(z_{q-2}')$, which implies $q-1=\partial(z_{0},z_{q-1})\leq1+\partial(z_{q-2}',z_{q-1})=q-2$, a contradiction. By the minimality of $l$, C$(q-1)$ exists. Since $\partial(z_{1},z_{q-2})=q-3$, from Lemma \ref{mix-(2,q-2) 1} (iii), we have $s=q-2$. Since $(1,q-3)$ is pure, by Lemma \ref{lemma Main1}, one gets $\wz{\partial}(z_{q-2},z_{0})=(2,q-2)$, and so $p_{(1,q-1),(1,q-1)}^{(2,q-2)}=0$. Pick a vertex $z_{2}'\in P_{(1,q-3),(q-2,1)}(z_{q-1},z_{1})$. Since $\partial(z_{0},z_{q-1})\leq\partial(z_{0},z_{2}')+q-3$, from Lemma \ref{mix-(2,q-2) 1} (ii), one obtains $\wz{\partial}(z_{0},z_{2}')=(2,q-2)$. By $z_{q-1}\in P_{(1,q-3),(1,q-1)}(z_{q-2},z_{0})$, there exists $z_{1}'\in P_{(1,q-1),(1,q-3)}(z_{0},z_{2}')$. Since $p_{(1,q-1),(1,q-1)}^{(2,q-2)}=0$, one has $\wz{\partial}(z_{q-1},z_{1}')=(1,q-2)$, which implies $p_{(1,q-1),(1,q-1)}^{(1,q-2)}=k_{1,q-1}$. Then $z_{0}\in P_{(1,q-1),(1,q-1)}(z_{q-1},z_{1})$, a contradiction. Then $\wz{\partial}(z_{q-1},z_{1})=(2,q-2)$. Similarly, $\wz{\partial}(z_{q-2},z_{0})=(2,q-2)$.

Since $z_{0}\in P_{(1,q-1),(1,p-1)}(z_{q-1},z_{1})$, there exists $z_{q-1}'\in P_{(1,q-1),(1,p-1)}(z_{q-2},z_{0})$. Note that the number of arcs of type $(1,p-1)$ in the circuit $(z_{q-1}',z_{0},z_{1},\ldots,z_{q-2})$ is more than the number of arcs of type $(1,p-1)$ in the circuit $(z_{0},z_{1},\ldots,z_{q-1})$. Repeat this process, there exists a circuit of length $q$ consisting of an arc of type $(1,q-1)$ and $q-1$ arcs of type $(1,p-1)$, contrary to $\Gamma_{1,q-1}\notin\Gamma_{l-1,1}^{q-1}$ for $1<l<q$.

\textbf{Case 2.2}.~$\partial(z_{i+1},z_{i})=q-1$ for some $i\in\{1,2,\ldots,q-2\}$.

Without loss of generality, we may assume $\partial(z_{q-1},z_{q-2})=q-1$. Note that $\wz{\partial}(z_{q-1},z_{1})=(2,q-2)$. By $z_{0}\in P_{(1,q-1),(1,p-1)}(z_{q-1},z_{1})$, we have $p_{(1,q-1),(1,q-1)}^{(2,q-2)}=0$, which implies $\wz{\partial}(z_{q-2},z_{0})=(1,q-2)$ and $p_{(1,q-1),(1,q-1)}^{(1,q-2)}=k_{1,q-1}$. Since C($q$) does not exist, $(1,q-2)$ is mixed. By the minimality of $q$, C$(q-1)$ or D$(q-1)$ exists.

Suppose that C$(q-1)$ exists. Pick a vertex $z_{1}''\in P_{(1,q-3),(q-2,1)}(z_{q-2},z_{0})$ and a circuit $(z_{q-2},z_{1}'',z_{2}''\ldots,z_{q-3}'')$ consisting of arcs of type $(1,q-3)$. Note that $(z_{q-3}'',z_{q-1}),(z_{q-1},z_{1}'')\in\Gamma_{2,q-2}$. Since $z_{q-2}\in P_{(1,q-3),(1,q-1)}(z_{q-3}'',z_{q-1})$, there exists $z_{0}''\in P_{(1,q-1),(1,q-3)}(z_{q-1},z_{1}'')$. By $p_{(1,q-2),(1,q-2)}^{(1,q-3)}=k_{1,q-2}$, we get $(z_{0}'',z_{0})\in\Gamma_{1,q-2}$. Since $p_{(1,q-1),(1,q-1)}^{(1,q-2)}=k_{1,q-1}$, one has $z_{q-1}\in P_{(1,q-1),(1,q-1)}(z_{0}'',z_{0})$, a contradiction.

Suppose that D$(q-1)$ exists. Pick vertices $w\in P_{(1,q-3),(q-3,1)}(z_{q-2},z_{0})$ and $z\in P_{(1,q-2),(1,q-3)}(z_{0},w)$. Since $q-1=\partial(z_{q-1},z_{q-2})\leq\partial(z_{q-1},w)+q-3$, we have $\wz{\partial}(z_{q-1},w)=(2,q-2)$. If $\wz{\partial}(z_{q-2},z)=(2,q-2)$, from $z_{0}\in P_{(1,q-2),(1,q-2)}(z_{q-2},z)$, then there exists $w'\in P_{(1,q-2),(1,q-2)}(z_{q-1},w)$, which implies $z_{q-1},w'\in\Gamma_{1,q-3}(z_{0})$ since $p_{(1,q-3),(q-3,1)}^{(1,q-2)}=k_{1,q-3}$, contrary to $q\neq2$. By Lemma \ref{mix-(2,q-2) 2}, $\wz{\partial}(z_{q-2},z)=(2,q-3)$. Since $z_{0}\in P_{(1,q-2),(1,q-2)}(z_{q-2},z)=\Gamma_{q-2,1}(z)$, we get $\wz{\partial}(z_{q-1},z)=(2,q-2)$. By $z_{0}\in P_{(1,q-1),(1,q-2)}(z_{q-1},z)$, there exists $w''\in P_{(1,q-1),(1,q-2)}(z_{q-1},w)$. The fact $p_{(1,q-3),(q-3,1)}^{(1,q-2)}=k_{1,q-3}$ implies $z_{q-2}\in\Gamma_{q-3,1}(w'')$, contrary to $\partial(w'',z_{q-2})\geq q-2$.$\qed$

In Lemmas \ref{(1,q-1)}, \ref{mix-(2,q-2) 2} and Lemma \ref{mix-(2,q-2) 1} (ii), we determine all the possible of the set $\Gamma_{1,q-1}^{2}$ for $q>2$ from Theorem \ref{Main1}.

Combining Theorem \ref{Main1} and Lemma \ref{lemma Main1}, we obtain the following corollary.

\begin{cor}\label{bkb}
If $p_{(1,s-1),(1,t-1)}^{(1,q-1)}\neq0$ with $s\neq t$, then $q\in\{s,t\}$.
\end{cor}

\section{The two-way distance}

In this section, we give some results about the two-way distance which are used frequently in the remaining of this paper.

\begin{lemma}\label{(2,q-2)}
If $p_{(1,p-1),(1,q-1)}^{(2,q-2)}\neq0$ with $q\neq p$, then {\rm C}$(q)$ exists and $p=q-1$.
\end{lemma}
\noindent\textbf{Proof.}~Suppose to, the contrary that C$(q)$ does not exist or $p\neq q-1$. Note that $(1,q-1)$ is mixed. By Lemma \ref{mix-(2,q-2) 1} (iii) and Theorem \ref{Main1}, D($q$) exists. Since $p_{(1,q-2),(q-2,1)}^{(1,q-1)}=k_{1,q-2}$, one obtains $p<q-1$. By Lemma \ref{mix-(2,q-2) 2}, we have $\Gamma_{1,q-1}^{2}=\{\Gamma_{2,q-1}\}$. Choose vertices $x_{q-1},x_{0},x_{1}$ such that $(x_{q-1},x_{1})\in\Gamma_{2,q-2}$ and $x_{0}\in P_{(1,q-1),(1,p-1)}(x_{q-1},x_{1})$.  Pick a path $(x_{1},x_{2},\ldots,x_{q-1})$. If $(x_{q-2},x_{q-1})\in\Gamma_{1,q-2}$, from $p_{(1,q-2),(q-2,1)}^{(1,q-1)}=k_{1,q-2}$, then $(x_{q-2},x_{0})\in\Gamma_{1,q-2}$, which implies $p=q-1$ since $(1,q-2)$ is pure, a contradiction. Since $\Gamma_{1,q-1}^{2}=\{\Gamma_{2,q-1}\}$, we obtain $\partial(x_{i+1},x_{i})\notin\{q-2,q-1\}$ for $0\leq i\leq q-2$. Since $\partial(x_{0},x_{q-2})=q-2$, from Corollary \ref{bkb}, we get $\partial(x_{q-2},x_{0})=2$. Since $p_{(1,p-1),(1,q-1)}^{(2,q-2)}\neq0$, there exists a vertex $x_{q-1}'\in P_{(1,q-1),(1,p-1)}(x_{q-2},x_{0})$. Repeat this process, there exists a circuit of length $q$ containing an arc of type $(1,q-1)$ and $(q-1)$ arcs of type $(1,p-1)$.

Without loss of generality, we may assume $(x_{i},x_{i+1})\in\Gamma_{1,p-1}$ for $0\leq i\leq q-2$. Since $q-1>p\geq2$, from Lemma \ref{jiben4} and Theorem \ref{Main1}, D$(p)$ exists or  $(1,p-1)$ is pure. If D$(p)$ exists, then there exists a vertex $x\in\Gamma_{1,p-2}(x_{i})$ for $0\leq i\leq q-1$, which implies $q-1=\partial(x_{0},x_{q-1})\leq1+\partial(x,x_{q-1})=p-1$, a contradiction; if $(1,p-1)$ is pure, from Lemma \ref{(1,q-1)}, then $(x_{0},x_{p-1})\in\Gamma_{p-1,1}$, which implies $x_{0}\in P_{(1,p-1),(1,p-1)}(x_{p-1},x_{p+1})$, a contradiction.$\qed$

\begin{lemma}\label{(1,q-1),(1,p-1)}
Let $\Gamma_{1,p-1}^{2}\cap\Gamma_{1,q-1}\Gamma_{q-1,1}\neq\emptyset$ and $\Gamma_{1,p-1}\notin\Gamma_{1,q-1}\Gamma_{q-1,1}$ with $q\neq p$.  \vspace{-0.3cm}
\begin{itemize}
\item [{\rm(i)}] If $(1,q-1)$ is pure, then $\Gamma_{1,p-1}\Gamma_{1,q-1}=\{\Gamma_{2,q}\}$.\vspace{-0.3cm}

\item [{\rm(ii)}] If {\rm C}$(q)$ exists with $(p,q)\neq(2,3)$, then $\Gamma_{1,p-1}\Gamma_{1,q-1}=\{\Gamma_{2,q}\}$.\vspace{-0.3cm}

\item [{\rm(iii)}] If {\rm D}$(q)$ exists with $(p,q)\neq(2,3)$, then $\Gamma_{1,p-1}\Gamma_{1,q-1}=\{\Gamma_{2,q-l}\}$, where $l=p_{(1,q-2),(q-2,1)}^{(1,p-1)}/k_{1,q-2}$.
\end{itemize}
\end{lemma}
\noindent\textbf{Proof.}~Pick vertices $x,y,z$ such that $\wz{\partial}(x,y)=(1,p-1)$ and $\wz{\partial}(y,z)=(1,q-1)$. Since $\Gamma_{1,p-1}^{2}\cap\Gamma_{1,q-1}\Gamma_{q-1,1}\neq\emptyset$, there exists a vertex $w\in P_{(p-1,1),(1,q-1)}(x,z)$.

Suppose that $\Gamma_{1,q-1}\in\Gamma_{1,p-1}^{2}$. Note that $x\in P_{(1,p-1),(1,p-1)}(w,z)$ and $x\in P_{(1,p-1),(1,p-1)}(y,z)$. It follows that $p=2$. Since $(x,y,z)$ is a path, one gets $q=3$, contrary to the fact that {\rm D}$(q)$ exists with $(p,q)\neq(2,3)$. Hence, $\Gamma_{1,q-1}\notin\Gamma_{1,p-1}^2$.

By $\Gamma_{1,q-1}\notin\Gamma_{1,p-1}^2$, we have $(x,z)\notin\Gamma_{1,p-1}$, and so $\Gamma_{1,q-1}\notin\Gamma_{1,p-1}\Gamma_{p-1,1}$. Since $\Gamma_{1,p-1}\notin\Gamma_{1,q-1}\Gamma_{q-1,1}$, from Corollary \ref{bkb}, one gets $\partial(x,z)=2$. Note that $\partial(z,x)\leq1+\partial(z,w)=q$. If $\partial(z,x)=q-2$, from Lemma \ref{(2,q-2)}, then C$(q)$ exists and $p=q-1$, which imply $z\in P_{(1,q-1),(1,q-1)}(w,x)$ and $q=3$, contrary to $(p,q)\neq(2,3)$. Then $\partial(z,x)=q-1$ or $q$, and so $\Gamma_{1,p-1}\Gamma_{1,q-1}\subseteq\{\Gamma_{2,q-1},\Gamma_{2,q}\}$. If $p_{(1,q-2),(q-2,1)}^{(1,q-1)}=p_{(1,q-2),(q-2,1)}^{(1,p-1)}=k_{1,q-2}$, then there exists a vertex $u$ such that $x,y,z\in\Gamma_{q-2,1}(u)$, which implies $\partial(z,x)\leq1+\partial(u,x)=q-1$, and so $\Gamma_{1,p-1}\Gamma_{1,q-1}=\{\Gamma_{2,q-1}\}$. It suffices to show that $p_{(1,q-2),(q-2,1)}^{(1,q-1)}=p_{(1,q-2),(q-2,1)}^{(1,p-1)}=k_{1,q-2}$ when $\Gamma_{2,q-1}\in\Gamma_{1,p-1}\Gamma_{q-1}$. Without loss of generality, we may assume $\wz{\partial}(x,z)=(2,q-1)$. By Theorem \ref{Main1}, we consider three cases.

\textbf{Case 1}. $(1,q-1)$ is pure.

Pick a path $(z=x_{0},x_{1},\ldots,x_{q-1}=x)$.

Suppose $|\{i\mid\partial(x_{i+1},x_{i})=q-1~\textrm{and}~0\leq i\leq q-2\}|\geq q-2$. Without loss of generality, we may assume $\partial(x_{i+1},x_{i})=q-1$ for $0\leq i\leq q-3$. Lemma \ref{(1,q-1)} implies $\wz{\partial}(x_{q-2},y)=(1,q-1)$. Since $\Gamma_{1,p-1}\notin\Gamma_{1,q-1}\Gamma_{q-1,1}$, from Corollary \ref{bkb}, we get $\wz{\partial}(x_{q-2},x)=(1,p-1)$ and $x\in P_{(1,p-1),(1,p-1)}(x_{q-2},y)$, contrary to $\Gamma_{1,q-1}\notin\Gamma_{1,p-1}^2$.

Suppose $|\{i\mid\partial(x_{i+1,x_{i}})=q-1~\textrm{and}~0\leq i\leq q-2\}|<q-2$. Without loss of generality, we may assume $\partial(x_{1},z)\neq q-1$ and $\partial(x_{2},x_{1})\neq q-1$. Since $(1,q-1)$ is pure, we have $\wz{\partial}(y,x_{1})=(2,q-1)$. By $y\in P_{(1,p-1),(1,q-1)}(x,z)$, there exists a vertex $z'\in P_{(1,p-1),(1,q-1)}(y,x_{1})$. Similarly, $(z',x_{2})\in\Gamma_{2,q-1}$ and there exists $z_{1}'\in P_{(1,p-1),(1,q-1)}(z',x_{2})$. Since $\Gamma_{1,p-1}^{2}\cap\Gamma_{1,q-1}\Gamma_{q-1,1}\neq\emptyset$, there exists $y'\in P_{(1,p-1),(1,p-1)}(x,z')$ such that $(y',x_{2})\in\Gamma_{1,q-1}$, contrary to the fact that $(x_{2},x_{3},\ldots,x_{q-1},y')$ is a path of length $q-2$.

\textbf{Case 2}. C($q$) exists.

Since $y\in P_{(1,p-1),(1,q-1)}(x,z)$ and $q\neq p$, we have $p_{(1,q-1),(1,q-1)}^{(2,q-1)}=0$. Pick a vertex $v\in P_{(1,q-2),(q-1,1)}(y,z)$. By Lemma \ref{mix-(2,q-2) 1} (ii), one gets $\Gamma_{1,q-1}^{2}=\{\Gamma_{1,q-2}\}$, which implies that $\wz{\partial}(w,v)=(1,q-2)$. If $p=q-1$, then $w\in P_{(1,p-1),(1,p-1)}(x,v)=\Gamma_{p-1,1}(v)$, and so $p=2$, contrary to the fact that C$(q)$ exists with $(p,q)\neq(2,3)$. Hence, $p\neq q-1$. Since $(1,q-2)$ is pure, one has $q-2\leq\partial(v,x)\leq1+\partial(v,w)=q-1$.

Suppose $\partial(x,v)=1$. Since $q\notin\{p,q-1\}$, from Corollary \ref{bkb}, one has $\wz{\partial}(x,v)=(1,q-2)$. By $p_{(1,q-1),(1,q-1)}^{(1,q-2)}=k_{1,q-1}$, we get $z\in P_{(1,q-1),(1,q-1)}(x,v)$, a contradiction.

Suppose $\partial(x,v)=2$. By Lemma \ref{mix-(2,q-2) 1} (i) or $y\in P_{(1,p-1),(1,q-1)}(x,z)$, there exists a vertex $y'\in\Gamma_{q-1,1}(v)$ such that $(x,y')\in\Gamma_{1,q-2}\cup\Gamma_{1,p-1}$. Since $p_{(1,q-1),(1,q-1)}^{(1,q-2)}=k_{1,q-1}$, one has $y'\in P_{(1,q-1),(1,q-1)}(y,v)$. In view of $p\neq q-1$ and Corollary \ref{bkb}, we get $(x,y')\in\Gamma_{1,p-1}$. By $x\in P_{(p-1,1),(1,p-1)}(y,y')$, one obtains $p_{(p-1,1),(1,p-1)}^{(1,q-1)}=k_{1,p-1}$, which implies $x\in P_{(p-1,1),(1,p-1)}(y,z)$, a contradiction.

\textbf{Case 3}. D($q$) exists.

Since $p_{(1,q-2),(q-2,1)}^{(1,q-1)}=k_{1,q-2}$ and $\partial(x,z)=2$, we have $p\neq q-1$. In view of $y\in P_{(1,p-1),(1,q-1)}(x,z)$, one gets $p_{(1,q-1),(1,q-1)}^{(2,q-1)}=0$, which implies $\Gamma_{1,q-1}^{2}=\{\Gamma_{2,q-2}\}$ from Lemma \ref{mix-(2,q-2) 2}. Choose a vertex $u$ such that $w,y,z\in\Gamma_{q-2,1}(u)$. Since $(1,q-2)$ is pure, we obtain $q-2\leq\partial(u,x)\leq1+\partial(u,w)=q-1$. By Corollary \ref{bkb}, one has $(x,u)\in\Gamma_{2,q-2}\cup\Gamma_{2,q-1}\cup\Gamma_{1,q-2}$.

Suppose $(x,u)\in\Gamma_{2,q-1}\cup\Gamma_{2,q-2}$. By $y\in P_{(1,p-1),(1,q-1)}(x,z)$ or $p_{(1,q-1),(1,q-1)}^{(2,q-2)}=k_{1,q-1}$, there exists a vertex $y'\in\Gamma_{q-1,1}(u)$ such that $(x,y')\in\Gamma_{1,p-1}\cup\Gamma_{1,q-1}$. Since $p_{(1,q-2),(q-2,1)}^{(1,q-1)}=k_{1,q-2}$, one has $y\in P_{(q-2,1),(1,q-2)}(y',u)$. In view of $p\neq q-1$ and Corollary \ref{bkb}, we get $(x,y')\in\Gamma_{1,p-1}$. Since $x\in P_{(p-1,1),(1,p-1)}(y,y')$, we obtain $p_{(p-1,1),(1,p-1)}^{(1,q-2)}=k_{1,p-1}$. Then $x\in P_{(p-1,1),(1,p-1)}(y,u)$, contrary to $\partial(x,u)=2$.

Thus, $(x,u)\in\Gamma_{1,q-2}$ and $p_{(1,q-2),(q-2,1)}^{(1,p-1)}=k_{1,q-2}$. The desired result holds.$\qed$

\begin{lemma}\label{Gamma 2,q-1}
Let $\Gamma_{1,q-1}^{2}=\{\Gamma_{2,q-2}\}$. If $\Gamma_{1,p-1}^{2}\cap\Gamma_{1,q-1}\Gamma_{q-1,1}\neq\emptyset$ and $\Gamma_{1,p-1}\Gamma_{1,q-1}=\{\Gamma_{2,q}\}$, then $\Gamma_{1,p-1}\Gamma_{1,q-1}^{2}=\{\Gamma_{3,q-1}\}$.
\end{lemma}
\noindent\textbf{Proof.}~Pick a path $(x_{q-1},x,x',x_{0})$ such that $(x_{q-1},x)\in\Gamma_{1,p-1}$ and $(x,x'),(x',x_{0})\in\Gamma_{1,q-1}$. Since $\Gamma_{1,p-1}^{2}\cap\Gamma_{1,q-1}\Gamma_{q-1,1}\neq\emptyset$, there exists $x_{q-2}\in P_{(p-1,1),(1,q-1)}(x_{q-1},x')$. Lemma \ref{(1,q-1)} implies that there exists a path $(x_{0},x_{1},\ldots,x_{q-2})$ consisting of arcs of type $(1,q-1)$. By $\Gamma_{1,p-1}\Gamma_{1,q-1}=\{\Gamma_{2,q}\}$, one has $(x_{q-1},x'),(x_{q-3},x_{q-1})\in\Gamma_{2,q}$. Since $(x_{q-1},x,x',x_{0},x_{1},\ldots,x_{q-3})$ is a shortest path, we get $\wz{\partial}(x_{q-1},x_{0})=(3,q-1)$. Then $\Gamma_{1,p-1}\Gamma_{1,q-1}^{2}=\{\Gamma_{3,q-1}\}$.$\qed$

\begin{lemma}\label{(1,q-1),(1,p-1)uniform}
Let $\Gamma_{1,p_{i}-1}^{2}\cap\Gamma_{1,q-1}\Gamma_{q-1,1}\neq\emptyset$ and $\Gamma_{1,p_{i}-1}\notin\Gamma_{1,q-1}\Gamma_{q-1,1}$ with $p_{i}\neq q$ for $i\in\{1,2\}$. If $ (2,3)\notin\{(p_{1},q),(p_{2},q)\}$ or $(1,2)$ is pure, then $\Gamma_{1,p_{1}-1}\Gamma_{1,q-1}=\Gamma_{1,p_{2}-1}\Gamma_{1,q-1}$.
\end{lemma}
\noindent\textbf{Proof.}~Suppose not. By Lemma \ref{(1,q-1),(1,p-1)}, we may assume $\Gamma_{1,p_{1}-1}\Gamma_{1,q-1}=\{\Gamma_{2,q-1}\}$ and $\Gamma_{1,p_{2}-1}\Gamma_{1,q-1}=\{\Gamma_{2,q}\}$. It follows that D($q$) exists, $p_{(1,q-2),(q-2,1)}^{(1,p_{1}-1)}=k_{1,q-2}$ and $p_{(1,q-2),(q-2,1)}^{(1,p_{2}-1)}=0$. Then $q-1\notin\{p_{1},p_{2}\}$. Pick vertices $x,y,z$ such that $\wz{\partial}(x,y)=(1,p_{2}-1)$ and $\wz{\partial}(y,z)=(1,q-1)$. By $\Gamma_{1,p_{2}-1}^{2}\cap\Gamma_{1,q-1}\Gamma_{q-1,1}\neq\emptyset$, there exists $w\in P_{(p_{2}-1,1),(1,q-1)}(x,z)$. Since $p_{(1,q-2),(q-2,1)}^{(1,q-1)}=k_{1,q-2}$, there exists a vertex $u$ such that $w,y,z\in\Gamma_{q-2,1}(u)$. The fact that $x\in P_{(1,p_{2}-1),(1,p_{2}-1)}(w,y)$ implies $\Gamma_{1,p_{2}-1}^{2}\cap\Gamma_{1,q-2}\Gamma_{q-2,1}\neq\emptyset$. Since $(1,q-2)$ is pure and $p_{(1,q-2),(q-2,1)}^{(1,p_{2}-1)}=0$, from Lemma \ref{(1,q-1),(1,p-1)} (i), we obtain $\wz{\partial}(x,u)=(2,q-1)$. By $\Gamma_{1,p_{1}-1}\Gamma_{1,q-1}=\{\Gamma_{2,q-1}\}$, there exists a vertex $y'\in P_{(1,p_{1}-1),(1,q-1)}(x,u)$. Since $p_{(1,q-2),(q-2,1)}^{(1,q-1)}=k_{1,q-2}$, we get $y\in P_{(q-2,1),(1,q-2)}(y',u)$. Corollary \ref{bkb} implies $p_{1}\in\{p_{2},q-1\}$, a contradiction.$\qed$

\begin{lemma}\label{(1,q-1),(1,p-1) C(q)}
Let $\Gamma_{1,p-1}^{2}\cap\Gamma_{1,q-1}\Gamma_{q-1,1}\neq\emptyset$. If C$(p)$ exists, then $\Gamma_{1,p-2}\in\Gamma_{1,q-1}\Gamma_{q-1,1}$.
\end{lemma}
\textbf{Proof.}~Suppose not. Since $\Gamma_{1,p-2}\in\Gamma_{1,p-1}^{2}$, from Lemma \ref{jiben3}, we get $p_{(1,q-1),(q-1,1)}^{(1,p-1)}=0$. By Lemma \ref{mix-(2,q-2) 1} (ii), one has $\Gamma_{1,p-1}^{2}=\{\Gamma_{1,p-2},\Gamma_{2,p-1}\}$ and $p_{(1,q-1),(q-1,1)}^{(2,p-1)}=k_{1,q-1}$.  Let $(x_{0},x_{1},x_{2},x_{3},x_{4})$ be a path consisting of arcs of type $(1,p-1)$ such that $(x_{0},x_{2}),(x_{2},x_{4})\in\Gamma_{2,p-1}$. It follows that there exists a vertices $x$ such that $x_{0},x_{2},x_{4}\in\Gamma_{q-1,1}(x)$.  By Lemma \ref{Gamma{1,q-1}=2} (ii), there exist vertices $x_{2}',x_{3}'$ such that $(x_{0},x_{2}')\in\Gamma_{1,p-2}$ and $(x_{2}',x_{3}'),(x_{3}',x_{4})\in\Gamma_{1,p-1}$. The fact $P_{(2,p-1),(2,p-1)}(x_{0},x_{4})=\Gamma_{p-1,2}(x_{4})$ implies $(x_{2}',x_{4})\in\Gamma_{1,p-2}$. Observe that $\Gamma_{1,p-2}^{2}\cap\Gamma_{1,q-1}\Gamma_{q-1,1}\neq\emptyset$. If $q=p-1$, from $p_{(1,p-1),(1,p-1)}^{(1,p-2)}=k_{1,p-1}$, then $x_{1}\in P_{(1,p-1),(1,p-1)}(x_{0},x)$ and $x_{1}\in P_{(1,p-1),(1,p-1)}(x_{2},x)$, contrary to $p\neq2$. Since $\Gamma_{1,p-2}\notin\Gamma_{1,q-1}\Gamma_{q-1,1}$ and $q\neq p-1$, from Lemma \ref{(1,q-1),(1,p-1)uniform}, there exists $x_{4}'\in P_{(1,p-1),(1,q-1)}(x_{2}',x)$. By $p_{(1,p-1),(1,p-1)}^{(1,p-2)}=k_{1,p-1}$, we have $\wz{\partial}(x_{4}',x_{4})=(1,p-1)$, contrary to $p_{(1,q-1),(q-1,1)}^{(1,p-1)}=0$.$\qed$

\begin{lemma}\label{(1,q-1),(1,p-1) C(q) 2}
Let $\Gamma_{1,p-1}^{2}\cap\Gamma_{1,q-1}\Gamma_{q-1,1}\neq\emptyset$ with $p\neq q-1$. If C$(p+1)$ exists and $\Gamma_{1,p-1}\Gamma_{1,q-1}=\{\Gamma_{2,q-1}\}$, then $\Gamma_{1,p}\in\Gamma_{1,q-2}\Gamma_{q-2,1}$.
\end{lemma}
\textbf{Proof.}~Since $\Gamma_{1,p-1}\Gamma_{1,q-1}=\{\Gamma_{2,q-1}\}$, one has $\Gamma_{1,p-1}\notin\Gamma_{1,q-1}\Gamma_{q-1,1}$. By Lemma \ref{(1,q-1),(1,p-1)}, $(1,q-2)$ is pure and $p_{(1,q-2),(q-2,1)}^{(1,q-1)}=p_{(1,q-2),(q-2,1)}^{(1,p-1)}=k_{1,q-2}$. Since $\Gamma_{1,p-1}\in\Gamma_{1,p}^{2}\cap\Gamma_{1,q-2}\Gamma_{q-2,1}$, from Lemma \ref{(1,q-1),(1,p-1)} (i), we get $\Gamma_{1,p}\Gamma_{1,q-2}=\{\Gamma_{1,q-2}\}$ or $\{\Gamma_{2,q-1}\}$.

Suppose $\Gamma_{1,p}\Gamma_{1,q-2}=\{\Gamma_{2,q-1}\}$. Let $x,y,y',z$ be vertices  such that $(x,z)\in\Gamma_{2,q-1}$, $y\in P_{(1,p-1),(1,q-1)}(x,z)$ and $y'\in P_{(1,p),(1,q-2)}(x,z)$. By $p_{(1,q-2),(q-2,1)}^{(1,q-1)}=k_{1,q-2}$, one has $\wz{\partial}(y',y)=(1,q-2)$. Since $y'\in P_{(1,p),(1,q-2)}(x,y)$, from Corollary \ref{bkb}, we get $p=q-1$, a contradiction$.\qed$

\begin{lemma}\label{(1,q-1),(1,p-1) D(p)}
Let $\Gamma_{1,p-1}^{2}\cap\Gamma_{1,q-1}\Gamma_{q-1,1}\neq\emptyset$. If D$(q)$ exists and $\Gamma_{1,p-1}\Gamma_{1,q-1}=\{\Gamma_{2,q}\}$, then $\Gamma_{1,q-1}^{2}=\{\Gamma_{2,q-2}\}$.
\end{lemma}
\textbf{Proof.}~Let $x,y,z,w,u$ be vertices such that $(x,y),(y,z)\in\Gamma_{1,p-1}$, $(z,w),(x,w)\in\Gamma_{1,q-1}$ and $(w,u)\in\Gamma_{1,q-2}$. By $p_{(1,q-2),(q-2,1)}^{(1,q-1)}=k_{1,q-2}$, one has $x,z\in\Gamma_{q-2,1}(u)$ and $p\neq q-1$. Since D$(q)$ exists and $\Gamma_{1,p-1}\Gamma_{1,q-1}=\{\Gamma_{2,q}\}$, from Lemma \ref{(1,q-1),(1,p-1)} (iii), we get $\Gamma_{1,p-1}\notin\Gamma_{1,q-2}\Gamma_{q-2,1}$. The fact that $y\in P_{(1,p-1),(1,p-1)}(x,z)$ and $u\in P_{(1,q-2),(q-2,1)}(x,z)$ imply $\Gamma_{1,p-1}^{2}\cap\Gamma_{1,q-2}\Gamma_{q-2,1}\neq\emptyset$. By  Lemma \ref{(1,q-1),(1,p-1)} (i), one obtains $(y,u)\in\Gamma_{2,q-1}$. If $p_{(1,q-1),(1,q-1)}^{(2,q-1)}=k_{1,q-1}$, then there exists a vertex $z'\in P_{(1,q-1),(1,q-1)}(y,u)$, which implies $z\in P_{(q-2,1),(1,q-2)}(z',u)$ since $p_{(1,q-2),(q-2,1)}^{(1,q-1)}=k_{1,q-2}$, contrary to Corollary \ref{bkb}. The desired result follows from Lemma \ref{mix-(2,q-2) 2}.$\qed$

\begin{lemma}\label{tongyi}
If $\Gamma_{1,p-1}^{2}\cap\Gamma_{1,h-1}\Gamma_{h-1,1}\neq\emptyset$ and $\Gamma_{1,h-1}^{2}\cap\Gamma_{1,q-1}\Gamma_{q-1,1}\neq\emptyset$ with $q\neq p$ and $h>2$, then $\Gamma_{1,p-1}^{2}\cap\Gamma_{1,q-1}\Gamma_{q-1,1}\neq\emptyset$.
\end{lemma}
\textbf{Proof.}~Choose vertices $x,y,z,w,v$ such that $(x,y),(y,z)\in\Gamma_{1,p-1}$, $(x,w),(z,w)\in\Gamma_{1,h-1}$ and $(z,v)\in\Gamma_{1,q-1}$. Then there exists a vertex $u\in P_{(1,h-1),(1,q-1)}(w,v)$. It suffices to show that $(x,v)\in\Gamma_{1,q-1}$. Suppose not. Since $v\in P_{(1,q-1),(q-1,1)}(z,u)$, we have $\wz{\partial}(x,u)\neq\wz{\partial}(z,u)$, which implies $|\Gamma_{1,h-1}^{2}|>1$. By Lemmas \ref{(1,q-1)}, \ref{mix-(2,q-2) 2} and Lemma \ref{mix-(2,q-2) 1} (ii), $(1,h-2)$ is pure and $\Gamma_{1,h-1}^{2}=\{\Gamma_{1,h-2},\Gamma_{2,h-1}\}$. In view of Lemma \ref{(1,q-1),(1,p-1) C(q)}, we have $\wz{\partial}(z,u)=(1,h-2)$ and $\wz{\partial}(x,u)=(2,h-1)$, which imply $p_{(1,q-1),(q-1,1)}^{(2,h-1)}=0$.

Suppose $p_{(1,h-1),(h-1,1)}^{(1,p-1)}=k_{1,h-1}$. This implies $w\in P_{(1,h-1),(h-1,1)}(y,z)$, and so $(y,u)\in\Gamma_{1,h-2}\cup\Gamma_{2,h-1}$. By $u\in P_{(1,h-2),(h-2,1)}(y,z)\cup P_{(2,h-1),(h-1,2)}(x,y)$, we get $p_{(1,h-2),(h-2,1)}^{(1,p-1)}=k_{1,h-2}$ or $p_{(2,h-1),(h-1,2)}^{(1,p-1)}=k_{2,h-1}$. Then $u\in P_{(1,h-2),(h-2,1)}(x,y)$ or $u\in P_{(2,h-1),(h-1,2)}(y,z)$, a contradiction. Thus, $p_{(1,h-1),(h-1,1)}^{(1,p-1)}=0$.

Suppose $p=h-1$. By $p_{(1,h-1),(1,h-1)}^{(1,h-2)}=k_{1,h-1}$, we have $w\in P_{(1,h-1),(1,h-1)}(x,y)$. Since $(y,z,w)$ is a circuit, one has $p=2$. Since $P_{(1,1),(1,1)}(x,z)=\Gamma_{1,1}(z)$, we obtain $(x,u)\in\Gamma_{1,1}$, contrary to $\partial(x,u)=2$. Hence, $p\neq h-1$.

Since $h>2$ and $\Gamma_{1,p-1}^{2}\cap\Gamma_{1,h-1}\Gamma_{h-1,1}\neq\emptyset$, one has $p\neq h$, which implies $(y,u)\notin\Gamma_{1,h-1}$ from Corollary \ref{bkb}. Since $p_{(1,h-1),(h-1,1)}^{(1,p-1)}=0$, we get $(y,w)\notin\Gamma_{1,h-1}$, and so $P_{(1,h-1),(1,h-1)}(y,u)=\emptyset$. By $\Gamma_{1,h-1}^{2}=\{\Gamma_{1,h-2},\Gamma_{2,h-1}\}$, one obtains $(y,u)\notin\Gamma_{1,h-2}\cup\Gamma_{2,h-1}$. Since $(1,h-2)$ is pure and $\wz{\partial}(x,u)=(2,h-1)$, we have $h-2\leq\partial(u,y)\leq1+\partial(u,x)=h$. It follows that $(y,u)\in\Gamma_{1,h}\cup\Gamma_{2,h-2}\cup\Gamma_{2,h}$.

If $(y,u)\in\Gamma_{1,h}$, from Corollary \ref{bkb}, then $p=h+1$, and so $(1,h)$ is pure by Theorem \ref{Main1}, contrary to $\wz{\partial}(x,u)=(2,h-1)$ and $w\in P_{(1,h-1),(1,h-1)}(x,w)$. Hence, $(y,u)\in\Gamma_{2,h-2}\cup\Gamma_{2,h}$. Note that $p\neq h-1$, $p_{(1,h-1),(h-1,1)}^{(1,p-1)}=0$ and $\Gamma_{1,p-1}^{2}\cap\Gamma_{1,h-1}\Gamma_{h-1,1}\neq\emptyset$. By Lemma \ref{mix-(2,q-2) 1} (i) or Lemma \ref{(1,q-1),(1,p-1)} (ii), there exists $z'\in\Gamma_{h-1,1}(u)$ such that $(y,z')\in\Gamma_{1,h-2}\cup\Gamma_{1,p-1}$. By $p_{(1,h-1),(1,h-1)}^{(1,h-2)}=k_{1,h-1}$, we get $z'\in P_{(1,h-1),(1,h-1)}(z,u)$. Since $p\neq h$, from Corollary \ref{bkb}, we obtain $(y,z')\in\Gamma_{1,p-1}$ and $p_{(p-1,1),(1,p-1)}^{(1,h-1)}=k_{1,p-1}$. Then $y\in P_{(p-1,1),(1,p-1)}(z',u)$, a contradiction.$\qed$

\section{Relationship between different types of arcs}

The main result of this section is the following important result which determine
the relationship between different types of arcs.

\begin{thm}\label{Main2}
Let $q,p\in T$ with $q\neq p$. Then one of the following holds:\vspace{-0.3cm}
\begin{itemize}
\item [{\rm(i)}] {\rm C}$(q)$ exists and $p=q-1$;\vspace{-0.3cm}

\item [{\rm(ii)}] {\rm C}$(p)$ exists and $q=p-1$;\vspace{-0.3cm}

\item [{\rm(iii)}] $\Gamma_{1,q-1}^{2}\cap\Gamma_{1,p-1}\Gamma_{p-1,1}\neq\emptyset$;\vspace{-0.3cm}

\item [{\rm(iv)}] $\Gamma_{1,p-1}^{2}\cap\Gamma_{1,q-1}\Gamma_{q-1,1}\neq\emptyset$.
\end{itemize}
\end{thm}

Let $\mathscr{B}$ denote the set consisting of $(p,p-1)$ and $(p-1,p)$ where C$(p)$ exists, and $\mathscr{C}=\{(p,q)\mid\Gamma_{1,q-1}^{2}\cap\Gamma_{1,p-1}\Gamma_{p-1,1}\neq\emptyset~\textrm{or}~\Gamma_{1,p-1}^{2}\cap\Gamma_{1,q-1}\Gamma_{q-1,1}\neq\emptyset\}$. The following result is key in the proof of Theorem \ref{Main2}.

\begin{lemma}\label{C}
Let $(q,p)\notin\mathscr{C}$ and $\Gamma_{2,l}\in\Gamma_{1,q-1}\Gamma_{1,p-1}$ with $q\neq p$. If $h\in\{q,p\}$ or $\{(p,h),(q,h)\}\subseteq\mathscr{C}\setminus\mathscr{B}$ for all $h$ with  $\Gamma_{m,l-1}\in\Gamma_{2,l}\Gamma_{1,h-1}$, then $(q,p)\in\mathscr{B}$.
\end{lemma}

Before we give a proof of Lemma \ref{C}, we show how it implies Theorem \ref{Main2}.

\noindent\textbf{Proof of Theorem \ref{Main2}.}~In order to prove Theorem \ref{Main2}, it suffices to show that $(q,p)\in\mathscr{B}\cup\mathscr{C}$. We shall prove it by contradiction. Since $\Gamma_{0,0}\in\Gamma_{1,1}^{2}\cap\Gamma_{1,s-1}\Gamma_{s-1,1}$ for $s\in T$, we have $q,p>2$. Since $(q,p)\notin\mathscr{C}$, from Lemma \ref{jiben3}, we have $p_{(1,q-1),(q-1,1)}^{(1,p-1)}=p_{(1,p-1),(p-1,1)}^{(1,q-1)}=0$. By Corollary \ref{bkb}, we can set $l=\min\{r\mid p_{(1,i-1),(1,j-1)}^{(2,r)}\neq0,~i\neq j,~(i,j)\notin\mathscr{B}\cup\mathscr{C}\}$. Without loss of generality, we may assume $p_{(1,q-1),(1,p-1)}^{(2,l)}\neq0$. Choose vertices $x,y,y'$ and $z$ with $(x,z)\in\Gamma_{2,l}$, $y\in P_{(1,q-1),(1,p-1)}(x,z)$ and $y'\in P_{(1,p-1),(1,q-1)}(x,z)$. By Lemma \ref{C}, there exists $x_{1}\in\Gamma_{1,h-1}(z)$ such that $\partial(x_{1},x)=l-1$, where $h\notin\{q,p\}$ and $\{(p,h),(q,h)\}\nsubseteq\mathscr{C}\setminus\mathscr{B}$.

\textbf{Case 1}. $\{(p,h),(q,h)\}\nsubseteq\mathscr{B}\cup\mathscr{C}$.

Without loss of generality, we may assume $(p,h)\notin\mathscr{B}\cup\mathscr{C}$. In view of Lemma \ref{jiben3}, one has $p_{(1,h-1),(h-1,1)}^{(1,p-1)}=p_{(1,p-1),(p-1,1)}^{(1,h-1)}=0$. By Corollary \ref{bkb}, we get $\partial(y,x_{1})=2$. It follows from the minimality of $l$ that $\partial(x_{1},y)=l$. By $z\in P_{(1,p-1),(1,h-1)}(y,x_{1})$, there exists a vertex $y''\in P_{(1,p-1),(1,h-1)}(x,z)$.

Suppose $(q,h)\notin\mathscr{B}\cup\mathscr{C}$. Similarly, $\wz{\partial}(y',x_{1})=(2,l)$. By $z\in P_{(1,q-1),(1,h-1)}(y',x_{1})$, there exist vertices $z'\in P_{(1,q-1),(1,h-1)}(y,x_{1})$ and $y'''\in P_{(1,q-1),(1,h-1)}(x,z)$. Since $P_{(1,h-1),(1,h-1)}(y''',x_{1})=\Gamma_{h-1,1}(x_{1})$ and $P_{(1,q-1),(1,q-1)}(x,z')=\Gamma_{1,q-1}(x)$, we obtain $\wz{\partial}(y''',z')=(1,q-1)=(1,h-1)$, a contradiction.

Suppose $(q,h)\in\mathscr{B}\cup\mathscr{C}$. Since $p_{(1,p-1),(p-1,1)}^{(1,h-1)}=0$, we have $(y',y'')\notin\Gamma_{1,h-1}$. By $p_{(1,p-1),(p-1,1)}^{(1,q-1)}=0$, one gets $(y'',y')\notin\Gamma_{1,q-1}$. Then $p_{(1,h-1),(1,h-1)}^{(1,q-1)}=p_{(1,q-1),(1,q-1)}^{(1,h-1)}=0$. Thus, $(q,h)\notin\mathscr{B}$ and $(q,h)\in\mathscr{C}$. Observe $\Gamma_{2,l}\in\Gamma_{1,p-1}\Gamma_{1,q-1}\cap\Gamma_{1,p-1}\Gamma_{1,h-1}$ and $(q,p),(h,p)\notin\mathscr{B}\cup\mathscr{C}$. It follows that the proofs for the case $\Gamma_{1,q-1}^{2}\cap\Gamma_{1,h-1}\Gamma_{h-1,1}\neq\emptyset$ and the case $\Gamma_{1,h-1}^{2}\cap\Gamma_{1,q-1}\Gamma_{q-1,1}\neq\emptyset$ are similar. Without loss of generality, we may assume $\Gamma_{1,q-1}^{2}\cap\Gamma_{1,h-1}\Gamma_{h-1,1}\neq\emptyset$. Then there exists $w\in P_{(1,h-1),(1,q-1)}(y'',y')$.

By $P_{(1,h-1),(1,h-1)}(y'',x_{1})=\Gamma_{1,h-1}(y'')$, we obtain $(w,x_{1})\in\Gamma_{1,h-1}$. Since $p_{(1,h-1),(h-1,1)}^{(1,p-1)}=p_{(1,p-1),(p-1,1)}^{(1,h-1)}=0$, from Corollary \ref{bkb}, one has $\partial(x,w)=2$. Note that $\partial(w,x)\leq1+\partial(x_{1},x)=l$. By the minimality of $l$, we get $\partial(w,x)=l$. Since $\Gamma_{2,l}\in\Gamma_{1,p-1}\Gamma_{1,q-1}$, there exists $w'\in P_{(1,p-1),(1,q-1)}(x,w)$. It follows that $w\in P_{(1,q-1),(1,q-1)}(w',y')$ and $x\in P_{(p-1,1),(1,p-1)}(w',y')$, contrary to $(q,p)\notin\mathscr{C}$.

\textbf{Case 2}. $\{(p,h),(q,h)\}\subseteq\mathscr{B}\cup\mathscr{C}$.

Without loss of generality, we may assume $(p,h)\in\mathscr{B}$ and $(q,h)\in\mathscr{C}$.

\textbf{Case 2.1}. $h=p+1$ and C($p+1$) exists.

Since $p_{(1,q-1),(q-1,1)}^{(1,p-1)}=0$ and $(q,p+1)\in\mathscr{C}$, from Lemma \ref{(1,q-1),(1,p-1) C(q)}, one gets $\Gamma_{1,q-1}^{2}\cap\Gamma_{1,p}\Gamma_{p,1}\neq\emptyset$. Choose vertices $w\in P_{(1,p),(1,p)}(y,z)$ and $w'\in P_{(q-1,1),(1,p)}(x,w)$. Since $(q,p)\notin\mathscr{C}$, we have $(w',z)\notin\Gamma_{1,p-1}$. By Lemma \ref{mix-(2,q-2) 1} (ii), one gets $\Gamma_{1,p}^{2}=\{\Gamma_{1,p-1},\Gamma_{2,p}\}$ and $(w',z)\in\Gamma_{2,p}$. In view of Lemma \ref{mix-(2,q-2) 1} (iii), we obtain $p-2<\partial(x_{1},x)=l-1\leq\partial(z,w')$, which implies $l=p$ or $p+1$.

Suppose $l=p$. Since $p_{(1,q-1),(1,p-1)}^{(2,p)}\neq0$, there exists $y''\in P_{(1,q-1),(1,p-1)}(w',z)$. By $P_{(1,q-1),(1,q-1)}(w',y)=\Gamma_{1,q-1}(w')$, we have $\wz{\partial}(y'',y)=(1,q-1)$, which implies $\Gamma_{1,q-1}\in\Gamma_{1,p-1}\Gamma_{p-1,1}$,  a contradiction.

Suppose $l=p+1$. Since $\Gamma_{1,p}^{2}=\{\Gamma_{1,p-1},\Gamma_{2,p}\}$, one has $(x,w)\notin\Gamma_{1,p}$ and $\Gamma_{1,q-1}\notin\Gamma_{1,p}\Gamma_{p,1}$. By Lemma \ref{(1,q-1),(1,p-1)} (ii), we get $\Gamma_{1,q-1}\Gamma_{1,p}=\{\Gamma_{2,p+1}\}$. Then there exists $y'''\in P_{(1,q-1),(1,p)}(x,z)$. In view of $p_{(1,p),(1,p)}^{(1,p-1)}=k_{1,p}$, one obtains $\wz{\partial}(y,y''')=(1,p)$ and $p_{(q-1,1),(1,q-1)}^{(1,p)}=k_{1,q-1}$, which implies $x\in P_{(q-1,1),(1,q-1)}(y''',z)$, a contradiction.

\textbf{Case 2.2}. $h=p-1$ and C($p$) exists.

Since $z\in P_{(1,p-1),(1,p-2)}(y,x_{1})$, there exists a vertex $z'\in P_{(1,p-2),(1,p-1)}(y,x_{1})$. By $p_{(1,p-1),(1,p-1)}^{(1,p-2)}=k_{1,p-1}$, we have $\wz{\partial}(z,z')=(1,p-1)$. If $\Gamma_{1,q-1}^{2}\cap\Gamma_{1,p-2}\Gamma_{p-2,1}\neq\emptyset$, then there exists $x'\in P_{(q-1,1),(1,p-2)}(x,z')$, which implies $z\in P_{(1,p-1),(1,p-1)}(x',z')$, contrary to $(q,p)\notin\mathscr{C}$. Thus, $\Gamma_{1,p-2}^{2}\cap\Gamma_{1,q-1}\Gamma_{q-1,1}\neq\emptyset$.

Since $(q,p)\notin\mathscr{C}$, we obtain $(x,z')\notin\Gamma_{1,q-1}$, and so $\Gamma_{1,p-2}\notin\Gamma_{1,q-1}\Gamma_{q-1,1}$. By Lemma \ref{(1,q-1),(1,p-1)}, we get $\Gamma_{1,p-2}\Gamma_{1,q-1}=\{\Gamma_{2,q-1}\}$ or $\{\Gamma_{2,q}\}$, and $\partial(x,z')=2$. If $\partial(z',x)=l$, from $y\in P_{(1,q-1),(1,p-2)}(x,z')$, then there exists $y''\in P_{(1,q-1),(1,p-2)}(x,z)$, which implies $(y'',y)\in\Gamma_{1,p-1}$ since $p_{(1,p-1),(1,p-1)}^{(1,p-2)}=k_{1,p-1}$, contrary to $p_{(1,q-1),(q-1,1)}^{(1,p-1)}=0$. By $\partial(z',x)\leq1+\partial(x_{1},x)$ and $\partial(z,x)\leq1+\partial(z',x)$, one has $\partial(z',x)=l-1$. If $\Gamma_{1,p-2}\Gamma_{1,q-1}=\{\Gamma_{2,q-1}\}$, from Lemmas \ref{(1,q-1),(1,p-1)} and \ref{(1,q-1),(1,p-1) C(q) 2}, then $p_{(1,q-2),(q-2,1)}^{(1,q-1)}=p_{(1,q-2),(q-2,1)}^{(1,p-1)}=k_{1,q-2}$, which implies that there exists a vertex $u$ with $x,y,z\in\Gamma_{q-2,1}(u)$, contrary to $q=\partial(z,x)\leq1+\partial(u,x)=q-1$. Hence, $\Gamma_{1,p-2}\Gamma_{1,q-1}=\{\Gamma_{2,q}\}$ and $\partial(z',x)=q=l-1$. Since $\Gamma_{1,p-2}^{2}\cap\Gamma_{1,q-1}\Gamma_{q-1,1}\neq\emptyset$, there exists $x_{2}\in P_{(1,q-1),(p-2,1)}(x,z')$. By $p_{(1,p-1),(1,p-1)}^{(1,p-2)}=k_{1,p-1}$, we have $\wz{\partial}(x_{1},x_{2})=(1,p-1)$.

Observe that $\Gamma_{1,p-2}\Gamma_{1,q-1}=\{\Gamma_{2,q}\}$. By Theorem \ref{Main1} and Lemmas \ref{(1,q-1),(1,p-1) D(p)}, \ref{mix-(2,q-2) 2}, $p_{(1,q-1),(1,q-1)}^{(2,q-2)}=k_{1,q-1}$ or C($q$) exists. In view of Lemma \ref{(1,q-1)}, if $p_{(1,q-1),(1,q-1)}^{(2,q-2)}=k_{1,q-1}$, then there exist $x_{3}\in P_{(2,q-2),(q-1,1)}(x,x_{2})$ and $x_{q}\in\Gamma_{q-1,1}(x)$ such that $\partial(x_{3},x_{q})=q-3$; if C($q$) exists, then there exist $x_{3}\in P_{(1,q-2),(q-1,1)}(x,x_{2})$ and $x_{q}\in\Gamma_{q-2,1}(x)$ such that $\partial(x_{3},x_{q})=q-3$. By the minimality of $l$, we have $\wz{\partial}(x_{1},x_{3})=(2,l)$ and $\wz{\partial}(x,x_{1})=\wz{\partial}(x_{q},z)=(3,l-1)$. Since $z\in P_{(2,l),(1,p-2)}(x,x_{1})$, there exist $y_{2}\in P_{(2,l),(1,p-2)}(x_{q},z)$ and $y_{1}\in P_{(1,q-1),(1,p-1)}(x_{q},y_{2})$. By $p_{(1,p-1),(1,p-1)}^{(1,p-2)}=k_{1,p-1}$, we get $\wz{\partial}(y_{2},y)=(1,p-1)$. Since $q\neq h$, from Lemmas \ref{(1,q-1)}, \ref{mix-(2,q-2) 2} and Lemma \ref{mix-(2,q-2) 1} (ii), one obtains $(y_{1},y)\notin\Gamma_{1,q-1}$, and so $P_{(1,q-1),(q-1,1)}(y_{1},x)=\emptyset$. Then C($q$) exists and $x_{p}\in\Gamma_{q-2,1}(x)$. Since $p_{(1,q-1),(1,q-1)}^{(1,q-2)}=k_{1,q-1}$, one gets $\wz{\partial}(y_{1},x)=(1,q-1)$. Note that $x\in P_{(1,q-1),(1,q-1)}(y_{1},y)$ and $y_{2}\in P_{(1,p-1),(1,p-1)}(y_{1},y)$. Since C$(q)$ and C$(p)$ both exist, from Lemma \ref{mix-(2,q-2) 1} (ii), we get $q=p$, a contradiction.$\qed$

In the remaining of this section, we give a proof of Lemma \ref{C}. Since $(q,p)\notin\mathscr{C}$, from Lemma \ref{jiben3}, we have $p_{(1,q-1),(q-1,1)}^{(1,p-1)}=p_{(1,p-1),(p-1,1)}^{(1,q-1)}=0$. Since $\Gamma_{0,0}\in\Gamma_{1,1}^{2}\cap\Gamma_{1,s-1}\Gamma_{s-1,1}$ for $s\in T$, we have $q,p>2$. Without loss of generality, we always assume $l=\min\{j\mid p_{(1,q-1),(1,p-1)}^{(2,j)}\neq0\}$. Choose vertices $x$, $y$ and $z$ with $\wz{\partial}(x,y)=(1,q-1)$, $\wz{\partial}(y,z)=(1,p-1)$ and $\wz{\partial}(x,z)=(2,l)$. Then there exists a vertex $y'\in P_{(1,p-1),(1,q-1)}(x,z)$. Next, we give three auxiliary lemmas.

\begin{lemma}\label{leq l-2}
If there exists a shortest path from $z$ to $x$ containing more than $l-2$ arcs of type $(1,q-1)$ or $(1,p-1)$, then $(q,p)\in\mathscr{B}$.
\end{lemma}
\textbf{Proof.}~Since the proofs are similar for both cases, we may assume that $(x_{0}=z,x_{1},\ldots,x_{l}=x)$ is a path with $(x_{i},x_{i+1})\in\Gamma_{1,p-1}$ for $0\leq i\leq l-2$. Suppose, to the contrary that $(q,p)\notin\mathscr{B}$. Lemma \ref{(2,q-2)} implies $l\geq3$.

\textbf{Case 1}. $p_{(1,p-1),(1,p-1)}^{(2,p-2)}=0$.

If C($p$) exists, from $l\geq3$, then there exists $z'\in P_{(1,p-1),(1,p-1)}(y,x_{1})$ such that $\wz{\partial}(z',x_{2})=(1,p-2)$, contrary to the minimality of $l$.

Since $(1,p-1)$ is mixed, from Theorem \ref{Main1}, $p_{(1,p-2),(p-2,1)}^{(1,p-1)}=k_{1,p-2}$ and $(1,p-2)$ is pure. By Lemma \ref{mix-(2,q-2) 2}, one gets $\Gamma_{1,p-1}^{2}=\{\Gamma_{2,p-1}\}$, which implies $p_{(1,q-1),(1,p-1)}^{(2,p-1)}=0$. Pick a vertex $w$ such that $y,x_{i}\in\Gamma_{p-2,1}(w)$ for $0\leq i\leq l-1$. Since $\partial(x_{1},y)=p-1$, we have $p-1<\partial(z,x)\leq2+\partial(w,x_{l-1})=p$, which implies $l=p$ and $\partial(w,x)\geq p-1$. By $\partial(w,x)\leq1+\partial(w,x_{l-1})=p-1$, we get $\partial(w,x)=p-1$. Since $p_{(1,q-1),(q-1,1)}^{(1,p-1)}=0$, we obtain $q\neq p-1$. By Corollary \ref{bkb}, one has $(x,w)\in\Gamma_{2,p-1}$. Since $\Gamma_{1,p-1}^{2}=\{\Gamma_{2,p-1}\}$, we obtain $y'\in P_{(1,p-1),(1,p-1)}(x,w)$. By $p_{(1,p-2),(p-2,1)}^{(1,p-1)}=k_{1,p-2}$, one has $y',x\in\Gamma_{1,p-2}(y)$, contrary to $q>2$.

\textbf{Case 2}. $p_{(1,p-1),(1,p-1)}^{(2,p-2)}=k_{1,p-1}$.

Note that $p_{(1,q-1),(1,p-1)}^{(2,p-2)}=0$ and $l+2>p$. Lemma \ref{(1,q-1)} implies $\wz{\partial}(x_{p-2},y)=(1,p-1)$. Let $(x_{l-1},x_{l})\in\Gamma_{1,h-1}$. If $l=p-1$, from $p_{(1,p-1),(p-1,1)}^{(1,q-1)}=0$, then $h\neq p$, and so $h=q$ by Corollary \ref{bkb}, which imply that C($q$) exists and $p=q-1$ from Lemmas \ref{(1,q-1)}, \ref{mix-(2,q-2) 2} and Lemma \ref{mix-(2,q-2) 1} (ii), contrary to $(q,p)\notin\mathscr{B}$; If $l>p$, from $P_{(1,p-1),(1,p-1)}(x_{p-2},x_{p})=\Gamma_{1,p-1}(x_{p-2})$, then $\wz{\partial}(y,x_{p})=(1,p-1)$, contrary to the minimality of $l$. Then $l=p$.

Since $\Gamma_{1,q-1}^{2}\cap\Gamma_{1,p-1}\Gamma_{p-1,1}=\emptyset$ and $x_{p-2}\in P_{(p-1,1),(1,p-1)}(y,x_{p-1})$, one gets $h\notin\{q,p\}$. Since $(h,p)\in\mathscr{C}$ from the assumption, we have $\Gamma_{1,h-1}^{2}\cap\Gamma_{1,p-1}\Gamma_{p-1,1}\neq\emptyset$ or $\Gamma_{1,p-1}^{2}\cap\Gamma_{1,h-1}\Gamma_{h-1,1}\neq\emptyset$. If $\Gamma_{1,p-1}^{2}\cap\Gamma_{1,h-1}\Gamma_{h-1,1}\neq\emptyset$, from Lemma \ref{(1,q-1)}, then $P_{(1,h-1),(h-1,1)}(x_{p-3},x_{p-1})=\Gamma_{h-1,1}(x)$, which implies $(x_{p-3},x)\in\Gamma_{1,h-1}$, a contradiction. Thus, $\Gamma_{1,h-1}^{2}\cap\Gamma_{1,p-1}\Gamma_{p-1,1}\neq\emptyset$.

Note that $\partial(x_{p-2},x)=2$. If $\partial(x,x_{p-2})=p$, by $p_{(1,q-1),(1,p-1)}^{(2,l)}\neq0$, then there exists $x_{p-1}'\in P_{(1,p-1),(1,q-1)}(x_{p-2},x)$, which implies $x_{p-2}\in P_{(p-1,1),(1,p-1)}(x_{p-1}',y)$, contrary to $(q,p)\notin\mathscr{C}$. Since $p_{(1,p-1),(1,p-1)}^{(2,p-2)}\neq0$, from Lemma \ref{mix-(2,q-2) 1} (ii), $(1,p-1)$ is pure or D$(p)$ exists. By Lemma \ref{(1,q-1),(1,p-1)} (i) or (iii), one has $\partial(x,x_{p-2})=p-1$ and $p_{(1,p-2),(p-2,1)}^{(1,p-1)}=p_{(1,p-2),(p-2,1)}^{(1,h-1)}=k_{1,p-2}$. Then there exists a vertex $w$ such that $x,z,y,x_{p-1},x_{p-2}\in\Gamma_{p-2,1}(w)$, contrary to $p=\partial(z,x)\leq1+\partial(w,x)=p-1$.$\qed$

\begin{lemma}\label{mathscr B}
If there exists a shortest path from $z$ to $x$ consisting of arcs of type $(1,q-1)$ and $(1,p-1)$, then $(q,p)\in\mathscr{B}$.
\end{lemma}
\textbf{Proof.}~Suppose not. By Lemma \ref{leq l-2}, we may assume that $(z,x_{1},x_{2},x_{3},x_{4})$ is a path such that $(z,x_{1}),(x_{1},x_{2})\in\Gamma_{1,q-1}$, $(x_{2},x_{3}),(x_{3},x_{4})\in\Gamma_{1,p-1}$ and $\partial(x_{4},x)=l-4$. The minimality of $l$ implies $(x_{1},x_{3})\in\Gamma_{2,l}$. It follows that $(x,x_{1}),(x_{1},x_{4})\in\Gamma_{3,l-1}$. Observe $\Gamma_{3,l-1}\in\Gamma_{1,q-1}^{2}\Gamma_{1,p-1}\cap\Gamma_{1,q-1}\Gamma_{1,p-1}^{2}$, contrary to Lemma \ref{jiben2}.$\qed$

\begin{lemma}\label{bkb lemma}
There exists a shortest path from $z$ to $x$ not containing an arc of type $(1,h-1)$ with $\Gamma_{1,q-1},\Gamma_{1,p-1}\in\Gamma_{1,h-1}\Gamma_{h-1,1}$.
\end{lemma}
\textbf{Proof.}~Suppose for the contradiction that any shortest path from $z$ to $x$ containing an arc of type $(1,h-1)$ with $\Gamma_{1,q-1},\Gamma_{1,p-1}\in\Gamma_{1,h-1}\Gamma_{h-1,1}$ for some $h\in T$. Then there exists $x_{1}\in\Gamma_{1,h-1}(z)$ such that $\partial(x_{1},x)=l-1$ with $\Gamma_{1,p-1},\Gamma_{1,q-1}\in\Gamma_{1,h-1}\Gamma_{h-1,1}$. Then $h\notin\{q,p\}$. Lemma \ref{mix-(2,q-2) 1} (i) and (iii) imply $(q,p)\notin\mathscr{B}$. Choose vertices $x_{l-1}\in\Gamma_{p-1,1}(x)$ and $x_{l-1}'\in\Gamma_{q-1,1}(x)$. By $p_{(1,h-1),(h-1,1)}^{(1,q-1)}=p_{(1,h-1),(h-1,1)}^{(1,p-1)}=k_{1,h-1}$, we get $x_{l-1},x_{l-1}',x,y\in\Gamma_{h-1,1}(x_{1})$, and so $l=h$.  Observe $l-1\leq\partial(z,x_{l-1}),\partial(z,x_{l-1}')\leq1+\partial(x_{1},x_{l-1})=l$. By Lemma \ref{jiben2}, we may assume $(x_{l-1},z)\notin\Gamma_{3,l}$.

\textbf{Case 1}. $\partial(z,x_{l-1})=l-1$.

By the hypothesis, we may assume that there exists $x_{l-2}\in\Gamma_{r-1,1}(x_{l-1})$ such that $\partial(z,x_{l-2})=l-2$ and $\Gamma_{1,p-1},\Gamma_{1,q-1}\in\Gamma_{1,r-1}\Gamma_{r-1,1}$. It follows that $x,y,z\in\Gamma_{1,r-1}(x_{l-2})$. Since $l-2=\partial(z,x_{l-2})=r-1$, one has $r=l-1$. Note that $\partial(y,x_{l-1})\leq1+\partial(y,x_{l-2})=r=l-1$, contrary to the minimality of $l$.

\textbf{Case 2}. $\wz{\partial}(x_{l-1},z)=(1,l)$.

Since $\Gamma_{1,l}\in\Gamma_{1,l-1}\Gamma_{l-1,1}$, from Theorem \ref{Main1}, D($l+1$) exists. By $p>2$, $x_{l-1}\notin P_{(1,l),(1,l)}(x,z)$. Since $\wz{\partial}(x,z)=(2,l)$, from Lemma \ref{mix-(2,q-2) 2}, one has $\Gamma_{1,l}^{2}=\{\Gamma_{2,l-1}\}$. Pick vertices $x_{l-2}\in\Gamma_{l,1}(x_{l-1})$ and $x'\in P_{(1,q-1),(1,p-1)}(x_{l-1},y)$. By $p_{(1,l-1),(l-1,1)}^{(1,l)}=p_{(1,l-1),(l-1,1)}^{(1,q-1)}=k_{1,l-1}$, we get $x',x_{l-2}\in\Gamma_{l-1,1}(x_{1})$. Note that $l-1\leq\partial(x',x_{l-2})\leq1+\partial(x_{1},x_{l-2})=l$. Since $x\in P_{(1,p-1),(1,p-1)}(x_{l-1},y')$ and $(q,p)\notin\mathscr{C}$, we have $(x_{l-1},z)\notin\Gamma_{1,q-1}$, and so $q\neq h+1$. Since $p_{(1,l-1),(l-1,1)}^{(1,l)}=k_{1,l-1}$ and $p_{(1,l),(1,l)}^{(2,l-1)}=k_{1,l}$, one obtains $\partial(x',x_{l-2})=l$.

Suppose $\partial(x_{l-2},x')=1$. Since $x_{l-1}\in P_{(1,l),(1,l)}(x_{l-2},z)$, we get $\wz{\partial}(x',z)=(1,l)$. By $y\in P_{(1,p-1),(1,p-1)}(x',z)$, one has $x\in P_{(1,p-1),(1,p-1)}(x_{l-1},z)$, a contradiction.

Suppose $\partial(x_{l-2},x')=2$. Since $x_{l-1}\in P_{(1,l),(1,q-1)}(x_{l-2},x')$, there exists $y''\in P_{(1,q-1),(1,l)}(x,z)$. By $p_{(1,l-1),(l-1,1)}^{(1,q-1)}=k_{1,l-1}$, one has $\wz{\partial}(y'',x_{1})=(1,l-1)$. Note that $\partial(y'',x_{l-1})\leq1+\partial(x_{1},x_{l-1})=l$. Since $p_{(1,q-1),(q-1,1)}^{(1,p-1)}=p_{(1,p-1),(p-1,1)}^{(1,q-1)}=0$, from Corollary \ref{bkb} and the minimality of $l$, we get $\wz{\partial}(x_{l-1},y'')=(2,l)$. Since $z\in P_{(1,l),(l,1)}(x_{l-1},y'')$, we obtain $x_{l-1}\in P_{(1,l),(l,1)}(x_{l-2},x')$, contrary to $q\neq2$.

\textbf{Case 3}. $\wz{\partial}(x_{l-1},z)=(2,l)$.

Since $y'\in P_{(1,p-1),(1,q-1)}(x,z)$, there exists a vertex $y''\in P_{(1,p-1),(1,q-1)}(x_{l-1},z)$. By $P_{(1,p-1),(1,p-1)}(x_{l-1},y')=\Gamma_{1,p-1}(x_{l-1})$, we have $\wz{\partial}(y'',y')=(1,p-1)$, which implies $\Gamma_{1,p-1}\in\Gamma_{1,q-1}\Gamma_{q-1,1}$, a contradiction.$\qed$

By Lemma \ref{bkb lemma}, there exists a path $(z=x_{0},x_{1},\ldots,x_{l}=x)$ such that $\Gamma_{1,p-1}$ or $\Gamma_{1,q-1}\notin\Gamma_{1,h_{i}-1}\Gamma_{h_{i}-1,1}$ for $0\leq i\leq l-1$, where $h_{i}=\partial(x_{i+1},x_{i})+1$. Now we divide the proof of Lemma \ref{C} into two subsections according to separate assumptions, respectively.

\subsection{$\Gamma_{1,q-1}^{2}\cap\Gamma_{1,h_{i}-1}\Gamma_{h_{i}-1,1}=\emptyset$ or $\Gamma_{1,p-1}^{2}\cap\Gamma_{1,h_{i}-1}\Gamma_{h_{i}-1,1}=\emptyset$ for all $i$}

Since $h_{i}\in\{q,p\}$ or $\{(q,h_{i}),(p,h_{i})\}\subset\mathscr{C}$ for $0\leq i\leq l-1$, from Lemma \ref{tongyi}, we have $h_{i}\in\{q,p\}$, or $\Gamma_{1,h_{i}-1}^{2}\cap\Gamma_{1,p-1}\Gamma_{p-1,1}\neq\emptyset$ and $\Gamma_{1,h_{i}-1}^{2}\cap\Gamma_{1,q-1}\Gamma_{q-1,1}\neq\emptyset$ for $0\leq i\leq l-1$. If $h_{i}\in\{q,p\}$ for all $i$, from Lemma \ref{mathscr B}, then $(q,p)\in\mathscr{B}$. We only need to consider the case that $h_{i}\notin\{q,p\}$ for some $i$. Without loss of generality, we may assume $i=0$. Lemma \ref{mix-(2,q-2) 1} (i) and (iii) imply $(q,p)\notin\mathscr{B}$.

\begin{step}\label{5.1.1}
$(2,3)\notin\{(h_{i},p),(h_{i},q)\mid0\leq i\leq l-1\}$ or $(1,2)$ is pure.
\end{step}
\vspace{-1ex}

We prove it by contradiction. Without loss of generality, we may assume that $(2,3)=(h_{0},p)$ and $(1,2)$ is mixed. Since $(p,h_{0})\notin\mathscr{B}$, from Theorem \ref{Main1}, we have $p_{(1,1),(1,1)}^{(1,2)}=k_{1,1}$, and so $(y,x_{1})\in\Gamma_{1,1}$. Note that $\Gamma_{1,1}^{2}\cap\Gamma_{1,q-1}\Gamma_{q-1,1}\neq\emptyset$. By the minimality of $l$, we get $(y',x_{1})\notin\Gamma_{1,q-1}$, which implies $p_{(1,q-1),(q-1,1)}^{(1,1)}=0$. In view of Lemma \ref{(1,q-1),(1,p-1)}, we have $\Gamma_{1,1}\Gamma_{1,q-1}=\{\Gamma_{2,q-1}\}$ or $\{\Gamma_{2,q}\}$.

Suppose $\Gamma_{1,1}\Gamma_{1,q-1}=\{\Gamma_{2,q-1}\}$. Note that $(x,x_{1})\in\Gamma_{2,q-1}$ and $l=q$. By Lemma \ref{(1,q-1),(1,p-1)}, D$(q)$ exists and $p_{(1,q-2),(q-2,1)}^{(1,1)}=k_{1,q-2}$. It follows that there exists a vertex $u$ such that $x_{1},x,y,z\in\Gamma_{q-2,1}(u)$. Then $q=\partial(z,x)\leq1+\partial(u,x)=q-1$, a contradiction. Thus, $\Gamma_{1,1}\Gamma_{1,q-1}=\{\Gamma_{2,q}\}$ and $l=q+1$.

Suppose that $(1,q-1)$ is mixed. If C$(q)$ exists, from Lemma \ref{(1,q-1)}, then there exists a path $(y=y_{0},y_{1},\ldots,x=y_{q-1})$ such that $y_{1}\in P_{(q-1,1),(1,q-2)}(x,y)$ and $(y_{j},y_{j+1})\in\Gamma_{1,q-2}$ for $1\leq j\leq q-2$, which implies that $(z,x_{1},y_{0},y_{1},\ldots,y_{q-1})$ is a shortest path, contrary to $(q-1,q)\in\mathscr{C}\setminus\mathscr{B}$ from the assumption. By Theorem \ref{Main1}, D$(q)$ exists. Pick a vertex $u$ such that $y',z\in\Gamma_{q-2,1}(u)$. Since $q=\partial(z,x)-1\leq\partial(u,x)\leq2+\partial(u,y')=q$, one has $\partial(u,x)=q$. By $q>2$, we get $q+1\neq3$, which implies $(x,u)\in\Gamma_{2,q}$ from Corollary \ref{bkb}. Since $y'\in P_{(1,2),(1,q-2)}(x,u)$, there exists $y''\in P_{(1,q-2),(1,2)}(x,x_{1})$. Since $p_{(1,1),(1,1)}^{(1,2)}=k_{1,1}$ and $p_{(1,q-2),(q-2,1)}^{(1,q-1)}=k_{1,q-2}$, we obtain $\wz{\partial}(y,y'')=(1,1)=(1,q-2)$, contrary to $p\neq q$. Thus, $(1,q-1)$ is pure. Lemma \ref{(1,q-1)} implies that there exists a circuit $(y=y_{0},y_{1},\ldots,y_{q-1}=x)$ consisting of arcs of type $(1,q-1)$.

If $\Gamma_{1,2}^{2}=\{\Gamma_{2,1}\}$, then there exists $y_{q}\in P_{(1,2),(1,2)}(z,y)$, which implies that $(z,y_{q},y_{0},y_{1},\ldots,y_{q-1})$ is a shortest path consisting of arcs of types $(1,p-1)$ and $(1,q-1)$, contrary to Lemma \ref{mathscr B}. Since D($3$) exists, from Lemma \ref{mix-(2,q-2) 2}, we have $\Gamma_{1,2}^{2}=\{\Gamma_{2,2}\}$. Pick a circuit $(z,z_{1},z_{2},y)$ consisting of arcs of type $(1,2)$. Since $p_{(1,1),(1,1)}^{(1,2)}=k_{1,1}$, we have $z_{1},z_{2}\in\Gamma_{1,1}(x_{1})$.

Note that $\partial(z,y_{q-2})=l-1=q$. Since $q+1=l\leq\partial(y_{1},z_{2})\leq2+\partial(y_{q-2},z)+\partial(y_{1},y_{q-2})=q-1+\partial(y_{q-2},z)$, we have $\partial(y_{q-2},z)=2$ or $3$. If $\partial(y_{q-2},z)=2$, from $\Gamma_{1,1}\Gamma_{1,q-1}=\{\Gamma_{2,q}\}$, then there exists $y''\in P_{(1,q-1),(1,1)}(y_{q-2},z)$, which implies that $\wz{\partial}(y'',y)=(1,1)=(1,q-1)$ since $p_{(1,1),(1,1)}^{(1,2)}=k_{1,1}$ and $P_{(1,q-1),(1,q-1)}(y_{q-2},y)=\Gamma_{1,q-1,1}(y_{q-2})$, a contradiction. Thus, $\wz{\partial}(y_{q-2},z)=(3,q)$.

Since $(y_{1},y_{2},\ldots,y_{q-2},x,y,x_{1},z_{2})$ is a shortest path from the minimality of $l$, we have $\partial(x,z_{2})=3$. Note that $q+1=l=\partial(z,x)\leq\partial(z_{2},x)+2$. Then $(x,z_{2})\in\Gamma_{3,q-1}\cup\Gamma_{3,q}$. Observe that $\Gamma_{0,0}\in\Gamma_{1,1}^{2}\cap\Gamma_{1,q-1}\Gamma_{q-1,1}$ and $(1,q-1)$ is pure. Since $\Gamma_{1,q-1}\Gamma_{1,1}=\{\Gamma_{2,q}\}$, from Lemmas \ref{(1,q-1)} and \ref{Gamma 2,q-1}, one has $\Gamma_{1,q-1}^2\Gamma_{1,1}=\{\Gamma_{3,q-1}\}$. By Lemma \ref{jiben2}, we have $\wz{\partial}(x,z_{2})=(3,q)$. Since $x_{1}\in P_{(2,q),(1,1)}(x,z_{2})$, there exist vertices $w_{2}\in P_{(2,q),(1,1)}(y_{q-2},z)$ and $w_{1}\in P_{(1,q-1),(1,1)}(w_{1},y_{1})$. The fact $P_{(1,q-1),(1,q-1)}(y_{q-2},y)=\Gamma_{1,q-1}(y_{q-2})$ implies $(w_{1},y)\in\Gamma_{1,q-1}$. Note that $y\in P_{(1,q-1),(1,p-1)}(w_{1},z)$ and $(z,w_{2},w_{1})$ is a path, contrary to $l\geq3$.

\begin{step}\label{5.1.2}
$p_{(1,q-1),(q-1,1)}^{(1,h_{i}-1)}=p_{(1,p-1),(p-1,1)}^{(1,h_{i}-1)}=0$, $\Gamma_{1,h_{i}-1}\Gamma_{1,p-1}=\Gamma_{1,h_{0}-1}\Gamma_{1,p-1}$ and $\Gamma_{1,h_{i}-1}\Gamma_{1,q-1}=\Gamma_{1,h_{0}-1}\Gamma_{1,q-1}$ for $0\leq i\leq l-1$ with $h_{i}\neq\{q,p\}$.
\end{step}
\vspace{-1ex}

Suppose that $h_{i}\notin\{q,p\}$ for some $i\in\{0,1,\ldots,l-1\}$. Note that $\Gamma_{1,h_{i}-1}^{2}\cap\Gamma_{1,p-1}\Gamma_{p-1,1}\neq\emptyset$ and $\Gamma_{1,h_{i}-1}^{2}\cap\Gamma_{1,q-1}\Gamma_{q-1,1}\neq\emptyset$. There exists $x_{1}'\in\Gamma_{1,h_{i}-1}(z)$ such that $\partial(x_{1}',x)=l-1$. By the minimality of $l$, we get $(y,x_{1}')\notin\Gamma_{1,p-1}$ and $(y',x_{1}')\notin\Gamma_{1,q-1}$, which imply $p_{(1,q-1),(q-1,1)}^{(1,h_{i}-1)}=p_{(1,p-1),(p-1,1)}^{(1,h_{i}-1)}=0$. The desired result follows from Lemma \ref{(1,q-1),(1,p-1)uniform} and Step \ref{5.1.1}.

\vspace{3ex}

In the following, we reach a contradiction based on the above discussion.

Suppose $|\{i\mid h_{i}\notin\{q,p\}~\textrm{and}~0\leq i\leq l-1\}|>2$. Without loss of generality, we may assume $h_{1},h_{2}\notin\{q,p\}$. Pick a vertex $z'\in P_{(1,h_{0}-1),(1,p-1)}(y,x_{1})$. By Step \ref{5.1.2}, there exist vertices $x_{1}''\in P_{(1,h_{0}-1),(1,p-1)}(z',x_{2})$ and $x_{2}''\in P_{(1,h_{0}-1),(1,p-1)}(x_{1}'',x_{3})$. Since $\Gamma_{1,h_{0}-1}^{2}\cap\Gamma_{1,p-1}\Gamma_{p-1,1}\neq\emptyset$, there exists $z''\in P_{(1,h_{0}-1),(1,h_{0}-1)}(y,x_{1}'')$ such that $(z'',x_{3})\in\Gamma_{1,p-1}$, contrary to the minimality of $l$. Hence, there exists at most one integer $i\in\{1,2,\ldots,l-1\}$ such that $h_{i}\notin\{q,p\}$.

Suppose $\Gamma_{1,p-1}\Gamma_{1,h_{0}-1}=\{\Gamma_{2,l}\}$. In view of Lemma \ref{mathscr B}, there exists no shortest path from $z$ to $x$ consisting of arcs of types $(1,p-1)$ and $(1,q-1)$. By $\wz{\partial}(y,x_{1})=(2,l)$, we may assume $h_{1}\notin\{q,p\}$ and $h_{i}\in\{q,p\}$ for $2\leq i\leq l-1$. Since $y\in P_{(1,q-1),(1,p-1)}(x,z)$, there exists a vertex $z'\in P_{(1,q-1),(1,p-1)}(y,x_{1})$. In view of Step \ref{5.1.2}, one has $\wz{\partial}(z',x_{2})=(2,l)$, which implies that $(x_{2},x_{3},\ldots,x_{l}=x,y,z')$ is a shortest path consisting of arcs of types $(1,p-1)$ and $(1,q-1)$. Hence, there also exists a shortest path from $z$ to $x$ consisting of arcs of types $(1,q-1)$ and $(1,p-1)$, a contradiction. Then $\Gamma_{1,p-1}\Gamma_{1,h_{0}-1}\neq\{\Gamma_{2,l}\}$. Simialrly, $\Gamma_{1,q-1}\Gamma_{1,h_{0}-1}\neq\{\Gamma_{2,l}\}$.

Suppose $l\leq4$. Observe $\Gamma_{1,h_{0}-1}^{2}\cap\Gamma_{1,p-1}\Gamma_{p-1,1}\neq\emptyset$ and $\Gamma_{1,h_{0}-1}^{2}\cap\Gamma_{1,q-1}\Gamma_{q-1,1}\neq\emptyset$. By Steps \ref{5.1.1}, \ref{5.1.2} and Lemma \ref{(1,q-1),(1,p-1)}, we have $\partial(y,x_{1})=\partial(y',x_{1})=2$, $p-1\leq\partial(x_{1},y')<l$ and $q-1\leq\partial(x_{1},y)<l$. Since $p,q>2$, one gets $l=4$ and $\{q,p\}=\{3,4\}$. Since $(q,p)\notin\mathscr{B}\cup\mathscr{C}$, from Theorem \ref{Main1} and Lemma \ref{jiben3}, $(1,3)$ is pure. Lemma \ref{(1,q-1),(1,p-1)} (i) implies $4\leq\partial(x_{1},y')$ or $4\leq\partial(x_{1},y)$, a contradiction. Thus, $l\geq5$.

Suppose that the path $(x_{1},x_{2},\ldots,x_{l}=x)$ contains an arc of type $(1,p-1)$ and an arc of type $(1,q-1)$. Since $l>4$ and there exists at most one integer $i\in\{1,2,\ldots,l-1\}$ such that $h_{i}\notin\{q,p\}$, we may assume $h_{l-2}=q$ and $h_{l-3}=h_{l-1}=p$. In view of the minimality of $l$, one gets $\wz{\partial}(x_{l-3},x_{l-1})=(2,l)$, which implies $\wz{\partial}(x_{l-1},z)=\wz{\partial}(x,x_{1})=(3,l-1)$. Since $z\in P_{(2,l),(1,h_{0}-1)}(x,x_{1})$, there exist vertices $y''\in P_{(2,l),(1,h_{0}-1)}(x_{l-1},z)$ and $x'\in P_{(1,p-1),(1,q-1)}(x_{l-1},y'')$. By $P_{(1,p-1),(1,p-1)}(x_{l-1},y')=\Gamma_{1,p-1}(x_{l-1})$, one has $\wz{\partial}(x',y')=(1,p-1)$. In view of the minimality of $l$ again, we obtain $\wz{\partial}(x',z)=(2,l)$. Since $y''\in P_{(1,q-1),(1,h_{0}-1)}(x',z)$, we get $\Gamma_{2,l}\in\Gamma_{1,q-1}\Gamma_{1,h_{0}-1}$. Note that $\Gamma_{1,h_{0}-1}^{2}\cap\Gamma_{1,q-1}\Gamma_{q-1,1}\neq\emptyset$. By Steps \ref{5.1.1}, \ref{5.1.2} and Lemma \ref{(1,q-1),(1,p-1)}, one has $\Gamma_{1,q-1}\Gamma_{1,h_{0}-1}=\{\Gamma_{2,l}\}$, a contradiction. In view of Lemma \ref{leq l-2}, we may assume $h_{1}\notin\{q,p\}$ and $h_{i}=p$ for $2\leq i\leq l-1$. By Step \ref{5.1.2}, there exists $x_{2}'\in P_{(1,h_{0}-1),(1,p-1)}(x_{1},x_{3})$. Without loss of generality, we may assume $x_{2}=x_{2}'$.

Suppose $\Gamma_{1,p-1}\Gamma_{1,h_{0}-1}=\{\Gamma_{2,p-1}\}$. By Steps \ref{5.1.1}, \ref{5.1.2} and Lemma \ref{(1,q-1),(1,p-1)}, D$(p)$ exists and $p_{(1,p-2),(p-2,1)}^{(1,h_{0}-1)}=k_{1,p-2}$. Choose a vertex $w$ such that $\{x=x_{l},x_{l-1},\ldots,x_{0}=z\}\subseteq\Gamma_{p-2,1}(w)$. Since $\Gamma_{1,p-1}\Gamma_{1,h_{0}-1}=\{\Gamma_{2,p-1}\}$, from Lemma \ref{mix-(2,q-2) 2}, we get $p_{(1,p-1),(1,p-1)}^{(2,p-2)}=k_{1,p-1}$. By $\partial(z,x)\leq 1+\partial(w,x)=p-1$, one has $l=p-1$, contrary to $\Gamma_{1,p-1}\Gamma_{1,h_{0}-1}\neq\{\Gamma_{2,l}\}$. By Steps \ref{5.1.1}, \ref{5.1.2} and Lemma \ref{(1,q-1),(1,p-1)} again, one has $\Gamma_{1,p-1}\Gamma_{1,h_{0}-1}=\{\Gamma_{2,p}\}$ and $l>p$.

Suppose $p_{(1,p-1),(1,p-1)}^{(2,p-2)}=0$. Note that $(1,p-1)$ is mixed. Since $l\geq5$, from Theorem \ref{Main1} and Lemma \ref{jiben4}, D$(p)$ exists. Observe $\Gamma_{1,h_{0}-1}^{2}\cap\Gamma_{1,p-1}\Gamma_{p-1,1}\neq\emptyset$ and $\Gamma_{1,p-1}\Gamma_{1,h_{0}-1}=\{\Gamma_{2,p}\}$, contrary to Lemma \ref{(1,q-1),(1,p-1) D(p)}. Thus, $p_{(1,p-1),(1,p-1)}^{(2,p-2)}=k_{1,p-1}$.

In view of Lemma \ref{(1,q-1)}, we get $\wz{\partial}(x_{l-p+1},x)=(p-1,1)$. If $l>p+1$, from $P_{(p-1,1),(1,p-1)}(y',x_{l-p+1})=\Gamma_{p-1,1}(x_{l-p+1})$, then $(x_{l-p},y')\in\Gamma_{1,p-1}$, contrary to the minimality of $l$. Thus, $l=p+1$

Since $\Gamma_{1,p-1}\Gamma_{1,h_{0}-1}=\{\Gamma_{2,p}\}$, we have $\wz{\partial}(x_{1},x_{3})=(2,p)$. Then $l-1=\partial(x_{3},x_{1})\leq1+\partial(x_{3},x_{l-1})+\partial(x_{l-1},z)$ and $(x_{l-1},z)\in\Gamma_{2,p}\cup\Gamma_{3,p}$. If $\wz{\partial}(x_{l-1},z)=(2,p)$, then there exists $x'\in P_{(1,p-1),(1,h_{0}-1)}(x_{l-1},z)$, which implies $x'\in P_{(1,p-1),(1,p-1)}(x_{l-1},y')=\Gamma_{1,p-1}(x_{l-1})$, contrary to Corollary \ref{bkb}. Thus, $(x_{l-1},z)\in\Gamma_{3,p}$.

Observe $\Gamma_{1,h_{0}-1}^{2}\cap\Gamma_{1,p-1}\Gamma_{p-1,1}\neq\emptyset$ and $\Gamma_{1,p-1}\Gamma_{1,h_{0}-1}=\{\Gamma_{2,p}\}$. Lemma \ref{(1,q-1)} implies $\Gamma_{1,p-1}^{2}=\{\Gamma_{2,p-2}\}$. By Lemma \ref{Gamma 2,q-1}, one has $\Gamma_{1,p-1}^{2}\Gamma_{1,h_{0}-1}=\{\Gamma_{3,p-1}\}$. Note that $p-1\leq\partial(x_{3},x_{1})-1\leq\partial(x_{3},z)\leq l-1=p$. By Lemma \ref{jiben3}, we get $(z,x_{3})\in\Gamma_{3,p}$. Since $x\in P_{(1,p-1),(2,l)}(x_{l-1},z)$, there exist vertices $x_{1}'\in P_{(1,p-1),(2,l)}(z,x_{3})$ and $x_{2}'\in P_{(1,q-1),(1,p-1)}(x_{1}',x_{3})$, contrary to Lemma \ref{mathscr B}.

By the above discussion, we finish the proof of Lemma \ref{C} for this case.

\subsection{$\Gamma_{1,q-1}^{2}\cap\Gamma_{1,h_{i}-1}\Gamma_{h_{i}-1,1}\neq\emptyset$ and $\Gamma_{1,p-1}^{2}\cap\Gamma_{1,h_{i}-1}\Gamma_{h_{i}-1,1}\neq\emptyset$ for some $i$}

Without loss of generality, we may assume that $\Gamma_{1,q-1}^{2}\cap\Gamma_{1,h_{0}-1}\Gamma_{h_{0}-1,1}\neq\emptyset$ and $\Gamma_{1,p-1}^{2}\cap\Gamma_{1,h_{0}-1}\Gamma_{h_{0}-1,1}\neq\emptyset$. For convenience's sake, write $h$ instead of $h_{0}$. Since $q,p>2$, we have $h\notin\{q,p\}$. By Lemma \ref{mix-(2,q-2) 1} (i) and (iii), one gets $(q,p)\notin\mathscr{B}$.

Suppose $\Gamma_{1,p-1},\Gamma_{1,q-1}\notin\Gamma_{1,h-1}\Gamma_{h-1,1}$. Pick a vertex $x'\in P_{(p-1,1),(1,h-1)}(y,x_{1})$. Since $q,p>2$, from Lemma \ref{(1,q-1),(1,p-1)uniform}, we have $\Gamma_{1,q-1}\Gamma_{1,h-1}=\Gamma_{1,p-1}\Gamma_{1,h-1}$, which implies that there exists $z'\in P_{(1,q-1),(1,h-1)}(y,x_{1})$. By the minimality of $l$, one gets $\partial(x_{1},y)\leq l\leq\partial(z',x')\leq1+\partial(x_{1},x')=h$. If $\Gamma_{1,q-1}\Gamma_{1,h-1}=\{\Gamma_{2,l}\}$, then there exists $y''\in P_{(1,q-1),(1,h-1)}(x,z)$, which implies $\Gamma_{\wz{i}}\in\Gamma_{1,q-1}^{2}\Gamma_{1,h-1}\cap\Gamma_{1,q-1}\Gamma_{1,h-1}^{2}$ with $\wz{\partial}(x,x_{1})=\wz{i}$, contrary to Lemma \ref{jiben2}. Since $\partial(x_{1},y)<l\leq h$, from Lemma \ref{(1,q-1),(1,p-1)}, we have $\Gamma_{1,q-1}\Gamma_{1,h-1}=\Gamma_{1,p-1}\Gamma_{1,h-1}=\{\Gamma_{2,h-1}\}$ and $l=h$, which imply $p_{(1,h-2),(h-2,1)}^{(1,h-1)}=p_{(1,h-2),(h-2,1)}^{(1,q-1)}=p_{(1,h-2),(h-2,1)}^{(1,p-1)}=k_{1,h-2}$. Pick a vertex $w$ such that $x,y,z\in\Gamma_{h-2,1}(w)$. Note that $h=\partial(z,x)\leq1+\partial(w,x)=h-1$, a contradiction.

Since $\Gamma_{1,q-1}$ or $\Gamma_{1,p-1}\notin\Gamma_{1,h-1}\Gamma_{h-1,1}$, we may assume $\Gamma_{1,p-1}\in\Gamma_{1,h-1}\Gamma_{h-1,1}$ and $\Gamma_{1,q-1}\notin\Gamma_{1,h-1}\Gamma_{h-1,1}$. Since $p_{(1,h-1),(h-1,1)}^{(1,p-1)}=k_{1,h-1}$, we have $\wz{\partial}(y,x_{1})=(1,h-1)$. By Lemma \ref{(1,q-1),(1,p-1)}, we get $\Gamma_{1,q-1}\Gamma_{1,h-1}=\{\Gamma_{2,h-1}\}$ and $l=h$, or $\Gamma_{1,q-1}\Gamma_{1,h-1}=\{\Gamma_{2,h}\}$ and $l=h+1$. If $\Gamma_{1,q-1}\Gamma_{1,h-1}=\{\Gamma_{2,h-1}\}$ and $l=h$, then $p_{(1,h-2),(h-2,1)}^{(1,h-1)}=p_{(1,h-2),(h-2,1)}^{(1,q-1)}=k_{1,h-2}$ and there exists a vertex $w'$ such that $x,y,z,x_{1}\in\Gamma_{h-2,1}(w')$, which imply $h=\partial(z,x)\leq1+\partial(w',x)=h-1$, a contradiction.

Note that $\Gamma_{1,q-1}\Gamma_{1,h-1}=\{\Gamma_{2,h}\}$, $(x,x_{1})\in\Gamma_{2,h}$ and $l=h+1$. Since $\Gamma_{1,q-1}^{2}\cap\Gamma_{1,h-1}\Gamma_{h-1,1}\neq\emptyset$, there exists a vertex $x_{h}'\in P_{(q-1,1),(1,h-1)}(x,x_{1})$.  Next, we prove our result step by step.

\begin{stepp}\label{3,h}
$\wz{\partial}(x_{h}',z)=(3,h)$ and $\Gamma_{3,h}\in\Gamma_{1,q-1}^{2}\Gamma_{1,p-1}$.
\end{stepp}
\vspace{-1ex}

Since $h=\partial(z,x)-1\leq\partial(z,x_{h}')\leq1+\partial(x_{1},x_{h}')$, we have $\partial(z,x_{h}')=h$. It suffices to show that $\partial(x_{h}',z)=3$. Suppose $\partial(x_{h}',z)\leq2$. If $\partial(x_{h}',z)=2$, from $\Gamma_{1,q-1}\Gamma_{1,h-1}=\{\Gamma_{2,h}\}$, then there exists $x'\in P_{(1,q-1),(1,h-1)}(x_{h}',z)$, which implies $x_{h}'\in P_{(1,h-1),(1,h-1)}(x',x_{1})=\Gamma_{h-1,1}(x_{1})$, contrary to $q>2$. Then $\partial(x_{h}',z)=1$.

Since $x_{1}\in P_{(1,h-1),(h-1,1)}(x_{h}',z)$, from Theorem \ref{Main1}, D($h+1$) exists. Since $(q,p)\notin\mathscr{C}$, one has $(x_{h}',z)\notin\Gamma_{1,p-1}$, and so $p\neq h+1$. If $p_{(1,p-1),(1,h)}^{(2,h)}\neq0$ or $p_{(1,h),(1,h)}^{(2,h)}\neq0$, then there exists $y''\in\Gamma_{h,1}(x_{1})$ such that $(x,y'')\in\Gamma_{1,h}\cup\Gamma_{1,p-1}$, which implies $x,y''\in\Gamma_{1,h-1}(y)$ since $p_{(1,h-1),(h-1,1)}^{(1,h)}=p_{(1,h-1),(h-1,1)}^{(1,p-1)}=k_{1,h-1}$, a contradiction. Hence, $p_{(1,p-1),(1,h)}^{(2,h)}=p_{(1,h),(1,h)}^{(2,h)}=0$. By Lemma \ref{mix-(2,q-2) 2}, $\Gamma_{1,h}^{2}=\{\Gamma_{2,h-1}\}$. Pick a vertex $x_{1}'\in\Gamma_{1,h}(z)$. Since $p_{(1,h-1),(h-1,1)}^{(1,h)}=k_{1,h-1}$, one has $\wz{\partial}(x_{1}',x_{1})=(1,h-1)$. Note that $h-1\leq\partial(x_{1}',y)\leq1+\partial(x_{1},x_{1}')=h$. Since $p\neq h+1$, $p_{(1,h),(1,h)}^{(2,h-1)}=k_{1,h}$ and $p_{(1,p-1),(1,h)}^{(2,h)}=0$, we get $\partial(y,x_{1}')=1$. Corollary \ref{bkb} implies $(y,x_{1}')\in\Gamma_{1,h}$. By $P_{(1,h),(1,h)}(x_{h}',x_{1}')=\Gamma_{h,1}(x_{1}')$, one has $(x_{h}',y)\in\Gamma_{1,h}$ and $p_{(1,q-1),(1,q-1)}^{(1,h)}=k_{1,q-1}$, which imply $x\in P_{(1,q-1),(1,q-1)}(x_{h}',z)$, a contradiction.

\begin{stepp}\label{(1,h-1) is pure}
$(1,h-1)$ is pure.
\end{stepp}
\vspace{-1ex}

Suppose not.  If D($h$) exists, from $p_{(1,h-2),(h-2,1)}^{(1,h-1)}=k_{1,h-2}$, then there exists a vertex $w$ such that $x_{h}',z,x_{1}\in\Gamma_{h-2,1}(w)$, which implies $h+1=\partial(z,x)\leq2+\partial(w,x_{h}')=h$, a contradiction. By Theorem \ref{Main1}, C($h$) exists.

Choose a vertex $x_{2}'\in P_{(1,h-2),(h-1,1)}(x_{h}',x_{1})$. Since $(1,h-2)$ is pure, there exists a circuit $(x_{2}',x_{3}',\ldots,x_{h}')$ consisting of one type of arcs. The minimality of $l$ implies $\partial(y',x_{h-1}')=h$. If $(x_{h-1}',y')\in\Gamma_{1,h}$, then $(1,h)$ is mixed since $\partial(y',x_{h-1}')=3+\partial(x_{2}',x_{h-1}')=h$ and $\wz{\partial}(z,x_{1})=(1,h-1)$, which implies that $(1,h-1)$ is pure from Theorem \ref{Main1}, a contradiction. Thus, $(x_{h-1}',y')\in\Gamma_{2,h}\cup\Gamma_{3,h}$.

Suppose $(x_{h-1}',y')\in\Gamma_{2,h}$. Since $\Gamma_{1,q-1}\Gamma_{1,h-1}=\{\Gamma_{2,h}\}$, there exists a vertex $x_{h}''\in P_{(1,h-1),(1,q-1)}(x_{h-1}',y')$. Pick a vertex $x'\in P_{(1,p-1),(1,q-1)}(x_{h}',y')$. By $p_{(1,h-1),(1,h-1)}^{(1,h-2)}=p_{(1,h-1),(h-1,1)}^{(1,p-1)}=k_{1,h-1}$, one has $\wz{\partial}(x_{h}'',x_{h}')=\wz{\partial}(x_{h}'',x')=(1,h-1)$. In view of $y'\in P_{(1,q-1),(q-1,1)}(x',x_{h}'')$, we get $p_{(1,q-1),(q-1,1)}^{(1,h-1)}=k_{1,q-1}$, which implies $y'\in P_{(q-1,1),(1,q-1)}(z,x_{1})$, contrary to the minimality of $l$.

Suppose $(x_{h-1}',y')\in\Gamma_{3,h}$. By Step \ref{3,h}, there exists a path $(x_{h}',y_{1},y_{2},z)$ such that $(x_{h}',y_{1})\in\Gamma_{1,h-2}$, $(y_{1},y_{2})\in\Gamma_{1,q-1}$ and $(y_{2},z)\in\Gamma_{1,p-1}$. Since $p_{(1,h-1),(1,h-1)}^{(1,h-2)}=k_{1,h-1}$, we have $x_{1}\in P_{(1,h-1),(1,h-1)}(x_{h}',y_{1})$. Then $(y_{1},y_{2},z,x_{1})$ is a circuit consisting of arcs of types $(1,q-1),(1,p-1)$ and $(1,h-1)$, contrary to $2\notin\{q,p,h\}$.

\begin{stepp}\label{3,h+1}
$\Gamma_{3,h+1}\in\Gamma_{1,p-1}^{2}\Gamma_{1,q-1}$.
\end{stepp}
\vspace{-1ex}

Pick a path $(y_{1},y_{2},y_{3},y_{4},y_{5})$ such that $(y_{1},y_{2})\in\Gamma_{1,h-1}$, $(y_{2},y_{3}),(y_{3},y_{4})\in\Gamma_{1,p-1}$ and $(y_{4},y_{5})\in\Gamma_{1,q-1}$. Since $p_{(1,h-1),(h-1,1)}^{(1,p-1)}=k_{1,h-1}$, we have $y_{3},y_{4}\in\Gamma_{1,h-1}(y_{1})$. The fact $\Gamma_{1,q-1}\Gamma_{1,h-1}=\{\Gamma_{2,h}\}$ implies $\wz{\partial}(y_{1},y_{5})=(2,h)$. By the minimality of $l$, one gets $h\leq\partial(y_{5},y_{2})\leq h+1$. Since $\wz{\partial}(y_{1},y_{2})=(1,h-1)$, from Lemma \ref{(2,q-2)}, we have $\wz{\partial}(y_{2},y_{5})\neq(1,h+1)$. Lemma \ref{jiben2} and Step \ref{3,h} imply $(y_{2},y_{5})\in\Gamma_{1,h}\cup\Gamma_{2,h}\cup\Gamma_{2,h+1}\cup\Gamma_{3,h+1}$. It suffices to show that $(y_{2},y_{5})\in\Gamma_{3,h+1}$.

\textbf{Case 1}. $(y_{2},y_{5})\in\Gamma_{1,h}$.

Since $(y_{1},y_{5})\in\Gamma_{2,h}$, we have $p_{(1,h-1),(h-1,1)}^{(1,h)}=0$, which implies that D($h+1$) does not exist. By Lemma \ref{mix-(2,q-2) 1} (i), C($h+1$) does not exist. In view of Theorem \ref{Main1}, $(1,h)$ is pure. Since $p_{(1,h-1),(h-1,1)}^{(1,p-1)}=k_{1,h-1}$ and $(q,p)\notin\mathscr{C}$, we have $(y_{2},y_{5})\notin\Gamma_{1,p-1}\cup\Gamma_{1,q-1}$, and so  $h+1\notin\{q,p\}$.

Pick a vertex $y_{1}'\in\Gamma_{h,1}(y_{2})$. By Lemma \ref{(1,q-1)}, we have $\Gamma_{1,h}^{2}=\{\Gamma_{2,h-1}\}$ and $\partial(y_{5},y_{1}')=h-1=l-2$. Since $p_{(1,p-1),(p-1,1)}^{(1,q-1)}=p_{(1,q-1),(q-1,1)}^{(1,p-1)}=0$, from Corollary \ref{bkb} and the minimality of $l$, we get $(y_{3},y_{5})\in\Gamma_{2,h+1}$. If $\Gamma_{1,h}^{2}\cap\Gamma_{1,q-1}\Gamma_{q-1,1}\neq\emptyset$, then $y_{4}\in P_{(q-1,1),(1,q-1)}(y_{1}',y_{5})$, which implies that $(y_{2},y_{3},y_{4},y_{1}')$ is a circuit, contrary to $2\notin\{h+1,q,p\}$. Since $(q,h+1)\in\mathscr{C}$ from the assumption, we obtain $\Gamma_{1,q-1}^{2}\cap\Gamma_{1,h}\Gamma_{h,1}\neq\emptyset$. Since $p_{(1,h-1),(h-1,1)}^{(1,h)}=0$, one has $(y_{2},y_{4})\notin\Gamma_{1,h}$, and so $p_{(1,h),(h,1)}^{(1,q-1)}=0$. Lemma \ref{(1,q-1),(1,p-1)} (i) implies $\Gamma_{1,q-1}\Gamma_{1,h}=\{\Gamma_{2,h+1}\}$. It follows that there exists $y_{4}'\in P_{(1,q-1),(1,h)}(y_{3},y_{5})$. Similarly, $\wz{\partial}(y_{2},y_{4}')=(2,h+1)$. Since $y_{5}\in P_{(1,h),(h,1)}(y_{2},y_{4}')$, one obtains $p_{(1,h),(h,1)}^{(2,h+1)}=k_{1,h}$. Then $y_{2}\in P_{(h,1),(1,h)}(y_{3},y_{5})$, contrary to $h+1\neq p$.

\textbf{Case 2}. $(y_{2},y_{5})\in\Gamma_{2,h}\cup\Gamma_{2,h+1}$.

By $\Gamma_{1,q-1}\Gamma_{1,h-1}=\{\Gamma_{2,h}\}$ and $p_{(1,p-1),(1,q-1)}^{(2,h+1)}\neq0$, there exists $y_{4}''\in\Gamma_{q-1,1}(y_{5})$ such that $(y_{2},y_{4}'')\in\Gamma_{1,h-1}\cup\Gamma_{1,p-1}$. Since $P_{(1,h-1),(1,h-1)}(y_{1},y_{4}'')=\Gamma_{1,h-1}(y_{1})$ or $P_{(1,p-1),(1,p-1)}(y_{2},y_{4})=\Gamma_{1,p-1}(y_{2})$, we get $(y_{4},y_{4}'')\in\Gamma_{1,h-1}\cup\Gamma_{p-1,1}$, which implies $p_{(1,q-1),(q-1,1)}^{(1,h-1)}=k_{1,q-1}$ or $p_{(1,q-1),(q-1,1)}^{(1,p-1)}=k_{1,q-1}$. By $y_{4}\in P_{(1,h-1),(1,q-1)}(y_{1},y_{5})$ or Lemma \ref{jiben3}, one has $(y_{1},y_{5})\in\Gamma_{1,q-1}$ or $(q,p)\in\mathscr{C}$, a contradiction.

\begin{stepp}\label{4,h}
$\Gamma_{4,h}\in\Gamma_{1,p-1}^{2}\Gamma_{1,q-1}^{2}$.
\end{stepp}
\vspace{-1ex}

By Step \ref{3,h+1}, there exists a path $(y_{1},y_{2},y_{3},y_{4})$ such that $(y_{1},y_{2}),(y_{2},y_{3})\in\Gamma_{1,p-1}$, $(y_{3},y_{4})\in\Gamma_{1,q-1}$ and $(y_{1},y_{4})\in\Gamma_{3,h+1}$. Since $\Gamma_{1,q-1}^{2}\cap\Gamma_{1,h-1}\Gamma_{h-1,1}\neq\emptyset$, there exist vertices $y_{5},y_{6}$ such that $(y_{4},y_{5})\in\Gamma_{1,q-1}$ and $(y_{5},y_{6}),(y_{3},y_{6})\in\Gamma_{1,h-1}$. By $p_{(1,h-1),(h-1,1)}^{(1,p-1)}=k_{1,h-1}$, we have $y_{1},y_{2}\in\Gamma_{h-1,1}(y_{6})$. Since $h\leq\partial(y_{5},y_{1})\leq1+\partial(y_{6},y_{1})=h$, one gets $\partial(y_{5},y_{1})=h$. It suffices to show that $\partial(y_{1},y_{5})=4$.

\textbf{Case 1}. $(y_{1},y_{5})\in\Gamma_{1,h}$.

Since $y_{6}\in P_{(1,h-1),(h-1,1)}(y_{1},y_{5})$, D($h+1$) exists. If $p=h+1$, then $y_{5}\in P_{(1,h),(1,h)}(y_{1},y_{3})=\Gamma_{1,h}(y_{1})$ and $(y_{3},y_{4},y_{5})$ is a circuit, contrary to $2\notin\{q,p\}$. Hence, $p\neq h+1$. If $p_{(1,h),(1,h)}^{(2,h)}\neq0$, from $\Gamma_{1,q-1}\Gamma_{1,h-1}=\{\Gamma_{2,h}\}$, then $(y_{4},y_{6})\in\Gamma_{2,h}$ and there exists $y_{5}'\in P_{(1,h),(1,h)}(y_{4},y_{6})$, which implies $y_{5}',y_{4}\in\Gamma_{1,h-1}(y_{5})$ since $p_{(1,h-1),(h-1,1)}^{(1,h)}=k_{1,h-1}$, a contradiction. By Lemma \ref{mix-(2,q-2) 2}, one has $p_{(1,h),(1,h)}^{(2,h-1)}=k_{1,h}$. Then there exists a vertex $y_{0}\in P_{(h,1),(2,h-1)}(y_{1},y_{5})$. Since $p_{(1,h-1),(h-1,1)}^{(1,h)}=k_{1,h-1}$, we get $\wz{\partial}(y_{0},y_{6})=(1,h-1)$. Note that $h-1\leq\partial(y_{2},y_{0})\leq1+\partial(y_{6},y_{0})=h$.

Suppose $\partial(y_{0},y_{2})=1$. Since $p\neq h+1$, from Corollary \ref{bkb}, we get $(y_{0},y_{2})\in\Gamma_{1,h}$ and $p_{(h,1),(1,h)}^{(1,p-1)}=k_{1,h}$. Then $y_{2},y_{3}\in\Gamma_{h,1}(y_{5})$. Since $y_{4}\in P_{(1,q-1),(1,q-1)}(y_{3},y_{5})$, we get $p_{(1,q-1),(1,q-1)}^{(1,h)}=k_{1,q-1}$, and so $y_{4}\in P_{(1,q-1),(1,q-1)}(y_{1},y_{5})$, contrary to $(q,p)\notin\mathscr{C}$.

Suppose $\partial(y_{0},y_{2})=2$. Since $p\neq h+1$ and $p_{(1,h),(1,h)}^{(2,h-1)}=k_{1,h}$, we have $\partial(y_{2},y_{0})=h$. By $\Gamma_{1,q-1}\Gamma_{1,h-1}=\{\Gamma_{2,h}\}$, there exists $y_{1}'\in P_{(1,q-1),(1,h-1)}(y_{0},y_{2})$. In view of $p_{(1,h-1),(h-1,1)}^{(1,p-1)}=p_{(1,h-1),(h-1,1)}^{(1,h)}=k_{1,h-1}$, we get $y_{0},y_{1}\in\Gamma_{1,h-1}(y_{1}')$, a contradiction.

\textbf{Case 2}. $(y_{1},y_{5})\in\Gamma_{2,h}$.

Since $\Gamma_{1,q-1}\Gamma_{1,h-1}=\{\Gamma_{2,h}\}$, there exists $y_{3}'\in P_{(1,q-1),(1,h-1)}(y_{1},y_{5})$. The fact $P_{(1,h-1),(1,h-1)}(y_{3}',y_{6})=\Gamma_{h-1,1}(y_{6})$ implies $(y_{3}',y_{1})\in\Gamma_{1,h-1}$, a contradiction.

\textbf{Case 3}. $(y_{1},y_{5})\in\Gamma_{3,h}$.

By Step \ref{3,h}, there exists a path $(y_{1},y_{2}',y_{4}',y_{5})$ such that $\wz{\partial}(y_{1},y_{2}')=(1,p-1)$ and $\wz{\partial}(y_{2}',y_{4}')=\wz{\partial}(y_{4}',y_{5})=(1,q-1)$. Since $P_{(1,p-1),(1,p-1)}(y_{1},y_{3})=\Gamma_{1,p-1}(y_{1})$ and $P_{(1,q-1),(1,q-1)}(y_{3},y_{5})=\Gamma_{q-1,1}(y_{5})$, we have $\wz{\partial}(y_{2}',y_{3})=(1,p-1)$ and $\wz{\partial}(y_{3},y_{4}')=(1,q-1)$, contrary to $p_{(1,q-1),(q-1,1)}^{(1,p-1)}=0$.

\begin{stepp}\label{1,h or 4,h+1}
$q=h+1$, $(1,h)$ is pure and $\Gamma_{1,p-1}^{3}\cap\Gamma_{1,h}\Gamma_{h,1}\neq\emptyset$ if $\Gamma_{4,h+1}\notin\Gamma_{1,p-1}^{3}\Gamma_{1,q-1}$.
\end{stepp}
\vspace{-1ex}

By Step \ref{3,h+1}, there exists a path $(y_{3},y_{4},y_{5},y_{6})$ such that $\wz{\partial}(y_{3},y_{4})=\wz{\partial}(y_{4},y_{5})=(1,p-1)$, $\wz{\partial}(y_{5},y_{6})=(1,q-1)$ and $\wz{\partial}(y_{3},y_{6})=(3,h+1)$. Pick vertices $y_{1},y_{2}$ such that $\wz{\partial}(y_{1},y_{2})=(1,h-1)$ and $\wz{\partial}(y_{2},y_{3})=(1,p-1)$. Since $p_{(1,h-1),(h-1,1)}^{(1,p-1)}=k_{1,h-1}$, we have $y_{3},y_{4},y_{5}\in\Gamma_{1,h-1}(y_{1})$. By $\Gamma_{1,q-1}\Gamma_{1,h-1}=\{\Gamma_{2,h}\}$, we have $\wz{\partial}(y_{1},y_{6})=(2,h)$. In view of $\wz{\partial}(y_{3},y_{6})=(3,h+1)$, one gets $\partial(y_{6},y_{2})=h$ or $h+1$. Since $\wz{\partial}(y_{1},y_{2})=(1,h-1)$, from Lemma \ref{(2,q-2)}, we have $(y_{2},y_{6})\notin\Gamma_{1,h+1}$. By Lemma \ref{jiben2} and Step \ref{4,h}, we obtain $(y_{2},y_{6})\in\Gamma_{1,h}\cup\Gamma_{2,h}\cup\Gamma_{3,h}\cup\Gamma_{2,h+1}\cup\Gamma_{3,h+1}$.

\textbf{Case 1}. $(y_{2},y_{6})\in\Gamma_{2,h}\cup\Gamma_{2,h+1}$.

Since $\Gamma_{1,q-1}\Gamma_{1,h-1}=\{\Gamma_{2,h}\}$ and $\Gamma_{2,h+1}\in\Gamma_{1,q-1}\Gamma_{1,p-1}$, there exists a vertex $y_{4}'\in\Gamma_{q-1,1}(y_{6})$ such that $(y_{2},y_{4}')\in\Gamma_{1,h-1}\cup\Gamma_{1,p-1}$. By $p_{(1,h-1),(h-1,1)}^{(1,p-1)}=k_{1,h-1}$ or $P_{(1,p-1),(1,p-1)}(y_{2},y_{4})=\Gamma_{1,p-1}(y_{2})$, we have $y_{3},y_{4},y_{5}\in\Gamma_{h-1,1}(y_{4}')$ or $(y_{4}',y_{4})\in\Gamma_{1,p-1}$. It follows that $y_{6}\in P_{(1,q-1),(q-1,1)}(y_{5},y_{4}')$ and $p_{(1,q-1),(q-1,1)}^{(1,h-1)}=k_{1,q-1}$, or $y_{4}\in P_{(1,p-1),(1,p-1)}(y_{4}',y_{5})$ and $y_{6}\in P_{(1,q-1),(q-1,1)}(y_{4}',y_{5})$. Then $y_{6}\in P_{(1,q-1),(q-1,1)}(y_{1},y_{5})$ or $(q,p)\in\mathscr{C}$, a contradiction.

\textbf{Case 2}. $(y_{2},y_{6})\in\Gamma_{3,h}\cup\Gamma_{3,h+1}$.

By Steps \ref{3,h} and \ref{3,h+1}, there exists a path $(y_{2},y_{3}',y_{4}',y_{6})$ with $(y_{2},y_{3}')\in\Gamma_{1,p-1}$ and $(y_{4}',y_{6})\in\Gamma_{1,q-1}$ and $(y_{3}',y_{4}')\in\Gamma_{1,q-1}\cup\Gamma_{1,p-1}$. The fact $P_{(1,p-1),(1,p-1)}(y_{2},y_{4})=\Gamma_{1,p-1}(y_{2})$ implies $\wz{\partial}(y_{3}',y_{4})=(1,p-1)$. Since $P_{(1,q-1),(1,q-1)}(y_{3}',y_{6})=\Gamma_{q-1,1}(y_{6})$ or $P_{(1,p-1),(1,p-1)}(y_{3}',y_{5})=\Gamma_{1,p-1}(y_{3}')$, we have $\wz{\partial}(y_{3}',y_{5})=(1,q-1)$ or $\wz{\partial}(y_{4}',y_{5})=(1,p-1)$. By $p_{(1,q-1),(q-1,1)}^{(1,p-1)}=0$, one gets $\wz{\partial}(y_{3}',y_{5})=(1,q-1)$ and $\Gamma_{1,p-1}\in\Gamma_{1,q-1}^{2}$. Lemmas \ref{(1,q-1)}, \ref{mix-(2,q-2) 2} and Lemma \ref{mix-(2,q-2) 1} (ii) imply $(q,p)\in\mathscr{B}$, a contradiction.

\textbf{Case 3}. $(y_{2},y_{6})\in\Gamma_{1,h}$.

Since $(y_{1},y_{6})\in\Gamma_{2,h}$, we have $p_{(1,h-1),(h-1,1)}^{(1,h)}=0$, which implies that D($h+1$) does not exist. In view of Lemma \ref{mix-(2,q-2) 1} (i), C($h+1$) does not exist. By Theorem \ref{Main1}, $(1,h)$ is pure. Since $\Gamma_{1,q-1}^{2}\cap\Gamma_{1,h-1}\Gamma_{h-1,1}\neq\emptyset$, there exists $y_{7}\in P_{(1,h-1),(q-1,1)}(y_{1},y_{6})$.

Suppose $q\neq h+1$. Since $p_{(1,h-1),(h-1,1)}^{(1,h)}=p_{(1,h-1),(h-1,1)}^{(1,q-1)}=0$, we have $(y_{2},y_{7})\notin\Gamma_{1,h}\cup\Gamma_{1,q-1}$, which implies $\partial(y_{2},y_{7})=2$ from Corollary \ref{bkb}. The fact that $(1,h)$ is pure implies $p_{(1,h),(1,h)}^{(2,h-1)}=k_{1,h}$. Since $\partial(y_{7},y_{2})\leq 1+\partial(y_{7},y_{1})=h$, one gets $(y_{2},y_{7})\in\Gamma_{2,h}$. In view of $y_{1}\in P_{(h-1,1),(1,h-1)}(y_{2},y_{7})$, we obtain $p_{(h-1,1),(1,h-1)}^{(2,h)}=k_{1,h-1}$, which implies $y_{5}\in P_{(1,h-1),(h-1,1)}(y_{1},y_{6})$, contrary to $h\neq q$. Thus, the desired result holds.

\begin{stepp}\label{h-1}
$\Gamma_{i,h-1}\in\Gamma_{2,l}\Gamma_{1,q-1}^2$ for some $i\in\{1,2,3,4\}$.
\end{stepp}
\vspace{-1ex}

Pick a path $(y_{1},y_{2},y_{3},y_{4})$ such that $(y_{1},y_{2})\in\Gamma_{1,h-1}$, $(y_{2},y_{3})\in\Gamma_{1,p-1}$ and $(y_{3},y_{4})\in\Gamma_{1,q-1}$ with $(y_{2},y_{4})\in\Gamma_{2,l}$. Since $p_{(1,h-1),(h-1,1)}^{(1,p-1)}=k_{1,h-1}$, we have $\wz{\partial}(y_{1},y_{3})=(1,h-1)$. By $\Gamma_{1,q-1}^{2}\cap\Gamma_{1,h-1}\Gamma_{h-1,1}\neq\emptyset$, there exists a vertex $y_{5}\in P_{(1,h-1),(q-1,1)}(y_{1},y_{4})$. Choose a vertex $y_{6}\in\Gamma_{1,q-1}(x_{5})$. It suffices to show that $\partial(y_{6},y_{2})=h-1$. Suppose not. In view of $\Gamma_{1,q-1}\Gamma_{1,h-1}=\{\Gamma_{2,h}\}$, one gets $\wz{\partial}(y_{1},y_{6})=(2,h)$. Since $l=h+1$, from the minimality of $l$, we obtain $h\leq\partial(y_{6},y_{2})\leq1+\partial(y_{6},y_{1})=h+1$. Since $(y_{1},y_{2})\in\Gamma_{1,h-1}$, from Lemma \ref{(2,q-2)}, one has $(y_{2},y_{6})\notin\Gamma_{1,h+1}$. By Lemma \ref{jiben2} and Step \ref{4,h}, we get $(y_{2},y_{6})\in\Gamma_{1,h}\cup\Gamma_{2,h}\cup\Gamma_{3,h}\cup\Gamma_{2,h+1}\cup\Gamma_{3,h+1}\cup\Gamma_{4,h+1}$.

\textbf{Case 1}. $(y_{2},y_{6})\in\Gamma_{1,h}$.

Since $\wz{\partial}(y_{1},y_{6})=(2,h)$, from Lemma \ref{mix-(2,q-2) 1} (i) and Theorem \ref{Main1}, $(1,h)$ is pure. Since $\Gamma_{1,q-1}\Gamma_{1,h-1}=\{\Gamma_{2,h}\}$, one gets $\wz{\partial}(y_{1},y_{4})=(2,h)$, which implies that there exists $y_{3}'\in P_{(1,h-1),(1,h)}(y_{1},y_{4})$. If $q=h+1$, from $P_{(1,h),(1,h)}(y_{4},y_{6})=\Gamma_{h,1}(y_{6})$, then $y_{2}\in P_{(1,h),(1,h)}(y_{4},y_{6})$, and so $(y_{2},y_{3},y_{4})$ is a circuit, contrary to $2\notin\{q,p\}$. Hence, $q\neq h+1$. Since $h-1\leq\partial(y_{5},y_{3}')\leq1+\partial(y_{5},y_{1})=h$ and $(1,h)$ is pure, we obtain $\wz{\partial}(y_{3}',y_{5})=(2,h)$, which implies that there exists $y_{4}'\in P_{(1,h-1),(1,q-1)}(y_{3}',y_{5})$. By $P_{(1,h-1),(1,h-1)}{(y_{1},y_{4}')}=\Gamma_{1,h-1}(y_{1})$, we get $\wz{\partial}(y_{5},y_{4}')=(1,h-1)$, contrary to $q>2$.

\textbf{Case 2}. $(y_{2},y_{6})\in\Gamma_{2,h}\cup\Gamma_{2,h+1}$.

Since $\Gamma_{1,q-1}\Gamma_{1,h-1}=\{\Gamma_{2,h}\}$ and $\Gamma_{2,h+1}\in\Gamma_{1,q-1}\Gamma_{1,p-1}$, there exists a vertex $y_{4}'\in\Gamma_{q-1,1}(y_{6})$ such that $(y_{2},y_{4}')\in\Gamma_{1,h-1}\cup\Gamma_{1,p-1}$. By $P_{(1,q-1),(1,q-1)}(y_{4},y_{6})=\Gamma_{q-1,1}(y_{6})$, we have $\wz{\partial}(y_{4},y_{4}')=(1,q-1)$. Note that $\partial(y_{4},y_{2})\leq1+\partial(y_{4}',y_{2})=h$, or $y_{4}\in P_{(1,q-1),(1,q-1)}(y_{3},y_{4}')$ and $y_{2}\in P_{(p-1,1),(1,p-1)}(y_{3},y_{4}')$, contrary to the minimality of $l$ or $(q,p)\notin\mathscr{C}$.

\textbf{Case 3}. $(y_{2},y_{6})\in\Gamma_{3,h}\cup\Gamma_{3,h+1}$.

By Steps \ref{3,h} and \ref{3,h+1}, there exists a path $(y_{2},y_{3}',y_{4}',y_{6})$ such that $(y_{2},y_{3}')\in\Gamma_{1,p-1}$, $(y_{4}',y_{6})\in\Gamma_{1,q-1}$ and $(y_{3}',y_{4}')\in\Gamma_{1,q-1}\cup\Gamma_{1,p-1}$. It follows that $y_{4}'\in\Gamma_{q-1,1}(y_{6})=P_{(1,q-1),(1,q-1)}(y_{4},y_{6})$. Since $P_{(1,q-1),(1,q-1)}(y_{3},y_{4}')=\Gamma_{q-1,1}(y_{4}')$ and $p_{(1,p-1),(p-1,1)}^{(1,q-1)}=0$, we get $(y_{3}',y_{4}')\in\Gamma_{1,p-1}$. By $P_{(1,p-1),(1,p-1)}(y_{2},y_{4}')=\Gamma_{1,p-1}(y_{2})$, we get $(y_{3},y_{4}')\in\Gamma_{1,p-1}$. Then $\Gamma_{1,p-1}\in\Gamma_{1,q-1}^{2}$. In view of Lemmas \ref{(1,q-1)}, \ref{mix-(2,q-2) 2} and Lemma \ref{mix-(2,q-2) 1} (ii), $p=q-1$ and C($q$) exists, contrary to $(q,p)\notin\mathscr{B}$.

\textbf{Case 4}. $(y_{2},y_{6})\in\Gamma_{4,h+1}$.

Suppose $\Gamma_{4,h+1}\notin\Gamma_{1,p-1}^{3}\Gamma_{1,q-1}$. By Step \ref{1,h or 4,h+1}, $q=h+1$, $\Gamma_{1,p-1}^{3}\cap\Gamma_{1,h}\Gamma_{h,1}\neq\emptyset$ and $(1,h)$ is pure. Pick a path $(y_{0}',y_{1}',y_{2})$ consisting of arcs of type $(1,p-1)$ such that $\wz{\partial}(y_{0}',y_{4})=(1,h)$. By Lemma \ref{(1,q-1)}, one gets $\wz{\partial}(y_{0}',y_{6})=(3,h-2)$. Note that $h+1=\partial(y_{6},y_{2})\leq2+\partial(y_{6},y_{0}')=h$, a contradiction. Hence, $\Gamma_{4,h+1}\in\Gamma_{1,p-1}^{3}\Gamma_{1,q-1}$.

Let $(y_{2},y_{3}',y_{4}',y_{5}',y_{6})$ be a path such that $(y_{2},y_{3}'),(y_{3}',y_{4}'),(y_{4}',y_{5}')\in\Gamma_{1,p-1}$ and $(y_{5}',y_{6})\in\Gamma_{1,q-1}$. By $P_{(1,p-1),(1,p-1)}(y_{2},y_{4}')=\Gamma_{1,p-1}(y_{2})$ and $P_{(1,q-1),(1,q-1)}(y_{4},y_{6})=\Gamma_{q-1,1}(y_{6})$, we have $(y_{3},y_{4}')\in\Gamma_{1,p-1}$ and $(y_{4},y_{5}')\in\Gamma_{1,q-1}$. Hence, $\Gamma_{1,p-1}^{2}\cap\Gamma_{1,q-1}^{2}\neq\emptyset$. Since $q,p>2$, from Lemma \ref{(1,q-1)}, $(1,q-1)$ is mixed and $(1,p-1)$ is mixed.

If C($q$) exists, then there exists $y_{5}''\in P_{(1,q-1),(1,q-1)}(y_{4},y_{6})=\Gamma_{1,q-1}(y_{4})$ such that $(y_{3},y_{5}'')\in\Gamma_{1,q-2}$, contrary to $\partial(y_{2},y_{6})=4$. Hence, C$(q)$ does not exist. Similarly, C$(p)$ does not exist. By Theorem \ref{Main1}, D($q$) and D($p$) exist. By Lemma \ref{mix-(2,q-2) 2}, $q=p-1$ and $p_{(1,q-1),(q-1,1)}^{(1,p-1)}=k_{1,q-1}$, or $p=q-1$ and $p_{(1,p-1),(p-1,1)}^{(1,q-1)}=k_{1,p-1}$. Lemma \ref{jiben3} implies $(q,p)\in\mathscr{C}$, a contradiction.

\vspace{3ex}

In the following, we reach a contradiction based on the above discussion.

By Step \ref{h-1}, there exist distinct vertices $x_{1}',x_{2}'$ such that $\wz{\partial}(z,x_{1}')=\wz{\partial}(x_{1}',x_{2}')=(1,q-1)$ and $\partial(x_{2}',x)=h-1$. Since $(q,p)\notin\mathscr{B}$, from Lemma \ref{mathscr B}, there exists a vertex $x_{h}''\in\Gamma_{r-1,1}(x)$ such that $\partial(x_{2}',x_{h}'')=h-2$ with $r\notin\{q,p\}$.

If $\Gamma_{1,r-1}^{2}\cap\Gamma_{1,q-1}\Gamma_{q-1,1}=\emptyset$ or $\Gamma_{1,r-1}^{2}\cap\Gamma_{1,p-1}\Gamma_{p-1,1}=\emptyset$, then $r>2$ and $\Gamma_{1,q-1}^{2}\cap\Gamma_{1,r-1}\Gamma_{r-1,1}\neq\emptyset$ and $\Gamma_{1,p-1}^{2}\cap\Gamma_{1,r-1}\Gamma_{r-1,1}\neq\emptyset$ from Lemma \ref{tongyi}, which imply that there exist vertices $x'\in P_{(1,p-1),(1,r-1)}(x_{h}'',y')$ and $x_{1}''\in P_{(1,q-1),(1,q-1)}(z,x_{2}')$ such that $\wz{\partial}(x',x_{1}'')=(1,r-1)$, contrary to the minimality of $l$. Hence, $\Gamma_{1,r-1}^{2}\cap\Gamma_{1,q-1}\Gamma_{q-1,1}\neq\emptyset$ and $\Gamma_{1,r-1}^{2}\cap\Gamma_{1,p-1}\Gamma_{p-1,1}\neq\emptyset$. By the minimality of $l$ again, we have $\Gamma_{1,r-1}\notin\Gamma_{1,p-1}\Gamma_{p-1,1}\cup\Gamma_{1,q-1}\Gamma_{q-1,1}$. Note that $(x_{h}'',z)\in\Gamma_{1,h}\cup\Gamma_{2,h}\cup\Gamma_{3,h}$.

\textbf{Case 1}. $(x_{h}'',z)\in\Gamma_{1,h}$.

Suppose $q\neq h+1$. Since $\partial(x_{1}',x_{h}'')=h-1$ and $q\neq h$, from Theorem \ref{Main1} and Lemma \ref{mix-(2,q-2) 1} (iii),  D($h+1$) exists. By Corollary \ref{bkb}, one gets $\partial(x_{h}'',x_{1}')=2$ and $p_{(1,h),(1,h)}^{(2,h-1)}=0$. Lemma \ref{mix-(2,q-2) 2} implies $\Gamma_{1,h}^{2}=\{\Gamma_{2,h}\}$. Since $\wz{\partial}(x,x_{1})=(2,h)$, there exists a vertex $z''\in P_{(1,h),(1,h)}(x,x_{1})$. By $p_{(1,h-1),(h-1,1)}^{(1,h)}=k_{1,h-1}$, one has $z'',x\in\Gamma_{1,h-1}(y)$, contrary to $q>2$.

Suppose $q=h+1$. Observe $\Gamma_{1,h}\notin\Gamma_{1,h-1}\Gamma_{h-1,1}$ and $\Gamma_{1,h}\Gamma_{1,h-1}=\{\Gamma_{2,h}\}$. By Theorem \ref{Main1} and Lemma \ref{mix-(2,q-2) 1} (i), $(1,h)$ is pure. Since $\Gamma_{1,r-1}^{2}\cap\Gamma_{1,h}\Gamma_{h,1}\neq\emptyset$ and $\Gamma_{1,r-1}\notin\Gamma_{1,h}\Gamma_{h,1}$, from Lemma \ref{(1,q-1),(1,p-1)} (i), one has $(x_{h}'',y)\in\Gamma_{2,h+1}$. By $p_{(1,q-1),(1,p-1)}^{(2,h+1)}\neq0$, there exists $x''\in P_{(1,q-1),(1,p-1)}(x_{h}'',y)$, contrary to $(q,p)\notin\mathscr{C}$.

\textbf{Case 2}. $(x_{h}'',z)\in\Gamma_{2,h}$.

By $\Gamma_{1,q-1}\Gamma_{1,h-1}=\{\Gamma_{2,h}\}$, there exists a vertex $y''\in P_{(1,q-1),(1,h-1)}(x_{h}'',z)$. Since $p_{(1,h-1),(h-1,1)}^{(1,p-1)}=k_{1,h-1}$, we get $\wz{\partial}(y'',y)=(1,h-1)$ and $\wz{\partial}(x_{h}'',y)=(2,h)$. Since $(r,q)\notin\mathscr{B}$ from the assumption, $(r,q)\notin(2,3)$ or $(1,2)$ is pure by Theorem \ref{Main1}. Observe $\Gamma_{1,r-1}^{2}\cap\Gamma_{1,q-1}\Gamma_{q-1,1}\neq\emptyset$ and $\Gamma_{1,r-1}\notin\Gamma_{1,q-1}\Gamma_{q-1,1}$. Since $q\neq h$, from Lemma \ref{(1,q-1),(1,p-1)}, $\partial(y,x_{h}'')=q-1=h$ and D$(h+1)$ exists, contrary to $\Gamma_{1,q-1}\notin\Gamma_{1,h-1}\Gamma_{h-1,1}$.

\textbf{Case 3}. $(x_{h}'',z)\in\Gamma_{3,h}$.

By Step \ref{3,h}, there exists a path $(x_{h}'',x'',y'',z)$ such that $(x_{h}'',x'')\in\Gamma_{1,p-1}$ and $(x'',y''),(y'',z)\in\Gamma_{1,q-1}$. Since $P_{(1,q-1),(1,q-1)}(x'',z)=\Gamma_{q-1,1}(z)$, we have $\wz{\partial}(x'',y')=(1,q-1)$. Since $p_{(1,q-1),(q-1,1)}^{(1,p-1)}=p_{(1,p-1),(p-1,1)}^{(1,q-1)}=0$, from Corollary \ref{bkb} and the minimality of $l$, one has $\wz{\partial}(x_{h}'',y')=(2,h+1)$. Since $x\in P_{(1,r-1),(1,p-1)}(x_{h}'',y')$, there exists $y_{1}\in P_{(1,r-1),(1,p-1)}(x,z)$. Similarly, $\wz{\partial}(y_{1},x_{1}')=(2,h+1)$ and there exists $y_{2}\in P_{(1,r-1),(1,p-1)}(y_{1},x_{1}')$. By $\Gamma_{1,r-1}^{2}\cap\Gamma_{1,p-1}\Gamma_{p-1,1}\neq\emptyset$, there exists $y_{3}\in P_{(1,r-1),(1,r-1)}(x_{h}'',y_{1})$ with $\wz{\partial}(y_{3},x_{1}')=(1,p-1)$, contrary to the minimality of $l$.

We finish the proof of Lemma 5.1 for this case.

\section{Proof of Theorem \ref{Main4}}

Before proceeding with the details of the proof, we first give an outline of it. In Section 6.1, we prove the following result.

\begin{prop}\label{fenqingkuang}
Let $p$ be an integer such that $\Gamma_{1,s-1}^2\cap\Gamma_{1,p-1}\Gamma_{p-1,1}\neq\emptyset$ for all $s\in T\setminus I$, where $I=\{s\mid(s,p)\in\mathscr{B}~\textrm{or}~s=p\}$.\vspace{-0.3cm}
\begin{itemize}
\item [{\rm(i)}] One of the following holds:\vspace{-0.3cm}
\begin{itemize}
\item [{\rm C1)}] $p=q$, $I=\{q\}$ and $(1,q-1)$ is pure;\vspace{-0.3cm}

\item [{\rm C2)}] $p=q-1$, $I=\{q-1\}$ and {\rm C}$(q)$ does not exist;\vspace{-0.3cm}

\item [{\rm C3)}] $p\in\{q-1,q\}$, $I=\{q-1,q\}$ and {\rm C}$(q)$ exists.\vspace{-0.3cm}
\end{itemize}
\item [{\rm(ii)}] The following hold:\vspace{-0.3cm}
\begin{itemize}
\item [{\rm (a)}] If $(1,p-1)$ is pure, then $\Gamma_{1,p-1}^{i}\cap F_{T\setminus I}=\emptyset$ for $1\leq i\leq p-1$.\vspace{-0.3cm}

\item [{\rm (b)}] If $(1,p-1)$ is mixed, then $\Gamma_{1,p-1}^{i}\cap F_{T\setminus I}=\emptyset$ for $1\leq i\leq 2p-3$.
\end{itemize}
\end{itemize}
\end{prop}

According to separate assumptions based on Proposition \ref{fenqingkuang} (i), by using Proposition \ref{fenqingkuang} (ii), we show that $\Delta_{T\setminus I}$ is a thick weakly distance-regular subdigraph in Section 6.2, and determine the corresponding quotient digraph in Section 6.3.

\subsection{Proof of Proposition \ref{fenqingkuang}}

In order to prove Proposition \ref{fenqingkuang}, we need three lemmas.

\begin{lemma}\label{Gamma 1,q-1,Gamma 1,p-1 subseteq Gamma 1,p-1Gamma p-1,1}
Let $\Gamma_{1,h_{i}-1}^{2}\cap\Gamma_{1,s-1}\Gamma_{s-1,1}\neq\emptyset$ and $\Gamma_{1,h_{i}-1}\notin\Gamma_{1,s-1}\Gamma_{s-1,1}$ with $h_{i}\neq s$ for $i=1,2$. Suppose that $(2,3)\notin\{(h_{1},s),(h_{2},s)\}$ or $(1,2)$ is pure. Then the following hold:\vspace{-0.3cm}
\begin{itemize}
\item [{\rm(i)}] If $h_{1}\neq h_{2}$, then $\Gamma_{1,h_{1}-1}\Gamma_{1,h_{2}-1}\subseteq\Gamma_{1,s-1}\Gamma_{s-1,1}$.\vspace{-0.3cm}

\item [{\rm(ii)}] If $\Gamma_{1,h_{i}-1}^{2}\varsubsetneq\Gamma_{1,s-1}\Gamma_{s-1,1}$ for some $i\in\{1,2\}$, then $h_{1}=h_{2}$.
\end{itemize}
\end{lemma}
\textbf{Proof.}~Pick a path $(x,y,z,w)$ such that $(x,y),(y,z)\in\Gamma_{1,h_{1}-1}$ and $(z,w)\in\Gamma_{1,s-1}$. By Lemma \ref{(1,q-1),(1,p-1)uniform}, there exists a vertex $z'\in P_{(1,h_{2}-1),(1,s-1)}(y,w)$.

(i) Since $\Gamma_{1,h_{1}-1}^{2}\cap\Gamma_{1,s-1}\Gamma_{s-1,1}\neq\emptyset$, we may assume $(x,w)\in\Gamma_{1,s-1}$. By Theorem \ref{Main2}, if $(h_{1},h_{2})\in\mathscr{B}$, from Lemma \ref{mix-(2,q-2) 1} (i), then $|\Gamma_{1,h_{1}-1}\Gamma_{1,h_{2}-1}|=1$; if $(h_{1},h_{2})\in\mathscr{C}\setminus\mathscr{B}$, from Lemma \ref{(1,q-1),(1,p-1)}, then $|\Gamma_{1,h_{1}-1}\Gamma_{1,h_{2}-1}|=1$. Since $w\in P_{(1,s-1),(s-1,1)}(x,z')$ and $y\in P_{(1,h_{1}-1),(1,h_{2}-1)}(x,z')$, the desired result follows.

(ii) Since $\Gamma_{1,h_{i}-1}^{2}\varsubsetneq\Gamma_{1,s-1}\Gamma_{s-1,1}$  for some $i$, we may assume $(x,w)\notin\Gamma_{1,s-1}$. If $h_{1}\neq h_{2}$, from (i), then $w\in P_{(1,s-1),(s-1,1)}(x,z')=\Gamma_{1,s-1}(z')$, a contradiction.$\qed$

\begin{lemma}\label{both pure}
Let $s<h$. Suppose $\Gamma_{1,h-1}^{2}\cap\Gamma_{1,s-1}\Gamma_{s-1,1}\neq\emptyset$. Then $s=h-1$ and $(1,h-2)$ is pure. Moreover, if $\Gamma_{1,h-1}\notin\Gamma_{1,s-1}\Gamma_{s-1,1}$, then $h+1\notin T$.
\end{lemma}
\textbf{Proof.}~If $\Gamma_{1,h-1}\in\Gamma_{1,s-1}\Gamma_{s-1,1}$, then $h-1\leq s<h$, and so $s=h-1$, which imply that $(1,h-1)$ is mixed and $(1,h-2)$ is pure from Theorem \ref{Main1}. Now we consider the case that $\Gamma_{1,h-1}\notin\Gamma_{1,s-1}\Gamma_{s-1,1}$.

Let $x,y,z,w$ be vertices with $(x,y),(y,z)\in\Gamma_{1,h-1}$ and $(z,w),(x,w)\in\Gamma_{1,s-1}$. Lemma \ref{(1,q-1),(1,p-1)} implies $\partial(y,w)=2$. Since $h-2\leq\partial(w,y)\leq 1+\partial(w,x)=s$, from Lemma \ref{(2,q-2)}, we get $s=h-1$.

Suppose that $(1,h-2)$ is mixed. In view of Theorem \ref{Main1}, $(1,h-1)$ is pure. If $p_{(1,h-3),(h-3,1)}^{(1,h-2)}=k_{1,h-3}$ or $\Gamma_{1,h-2}^{2}=\{\Gamma_{1,h-3}\}$, then there exists a vertex $w'\in P_{(1,h-3),(h-3,1)}(x,w)\cup P_{(1,h-3),(h-2,1)}(x,w)$ such that $(z,w')\in\Gamma_{1,h-3}$, which implies $\partial(z,x)=1+\partial(w',x)=h-2$, contrary to the fact that $(1,h-1)$ is pure. By Theorem \ref{Main1} and Lemma \ref{mix-(2,q-2) 1} (ii), we have $p_{(1,h-2),(1,h-2)}^{(2,h-2)}=k_{1,h-2}$. Since $(1,h-1)$ is pure, from Lemma \ref{(1,q-1)}, one gets $(x,z)\in\Gamma_{2,h-2}$, which implies $w\in P_{(1,h-2),(1,h-2)}(x,z)$ and $h=3$, contrary to the fact that $(1,q-2)$ is mixed. Thus, $(1,h-2)$ is pure.

Suppose $h+1\in T$. In view of Lemma \ref{(1,q-1),(1,p-1)} (i), we obtain $\Gamma_{1,h-1}\Gamma_{1,h-2}=\{\Gamma_{2,h-1}\}$ and $\wz{\partial}(y,w)=(2,h-1)$. It follows that $(1,h)$ is mixed. If C$(h+1)$ exists, from Lemma \ref{mix-(2,q-2) 1} (i), then there exists a vertex $z'\in P_{(1,h),(1,h-1)}(y,w)$, which implies $\wz{\partial}(z',z)=(1,h)$ since $p_{(1,h),(1,h)}^{(1,h-1)}=k_{1,h}$, contrary to Corollary \ref{bkb}. By Theorem \ref{Main1}, D$(h+1)$ exists. Pick a vertex $x'\in\Gamma_{1,h}(x)$. The fact $p_{(1,h-1),(h-1,1)}^{(1,h)}=k_{1,h-1}$ implies that $(x',y)\in\Gamma_{1,h-1}$. Since $(1,h-1)$ is pure, from Lemma \ref{(1,q-1)}, one has $\Gamma_{1,h-1}^{2}=\{\Gamma_{2,h-2}\}$. By $\Gamma_{1,h-1}^{2}\cap\Gamma_{1,h-2}\Gamma_{h-2,1}\neq\emptyset$, we get $x'\in\Gamma_{h-2,1}(w)$, contrary to $h=\partial(x',x)\leq1+\partial(w,x)=h-1$. Thus, $h+1\notin T$.$\qed$

\begin{lemma}\label{(1,q-1),(1,1)}
Let $\Gamma_{1,h-1}^{2}\cap\Gamma_{1,s-1}\Gamma_{s-1,1}\neq\emptyset$ and $\Gamma_{1,h-1}^{2}\varsubsetneq\Gamma_{1,s-1}\Gamma_{s-1,1}$ with $h\neq s$. Suppose that $(h,s)\neq (2,3)$ or $(1,2)$ is pure. Then one of the following holds:\vspace{-0.3cm}
\begin{itemize}
\item [{\rm(i)}] {\rm C}$(h)$ exists and $\Gamma_{1,h-1}^{2}\setminus\Gamma_{1,s-1}\Gamma_{s-1,1}=\{\Gamma_{2,h-1}\}$;\vspace{-0.3cm}

\item [{\rm(ii)}] $h=2$ and $\Gamma_{1,1}^{2}\setminus\Gamma_{1,s-1}\Gamma_{s-1,1}=\{\Gamma_{2,2}\}$.
\end{itemize}
\end{lemma}
\textbf{Proof.}~In view of Lemmas \ref{(1,q-1)}, \ref{mix-(2,q-2) 2} and Lemma \ref{mix-(2,q-2) 1} (ii), C$(h)$ exists and $\Gamma_{1,h-1}^{2}=\{\Gamma_{1,h-2},\Gamma_{2,h-1}\}$, or $h=2$ and $\Gamma_{1,1}^{2}\subseteq\{\Gamma_{0,0},\Gamma_{1,2},\Gamma_{2,1},\Gamma_{2,2}\}$. Suppose that (i) does not hold. By Lemma \ref{(1,q-1),(1,p-1) C(q)}, one has $h=2$ and $\Gamma_{1,1}^{2}\subseteq\{\Gamma_{0,0},\Gamma_{1,2},\Gamma_{2,1},\Gamma_{2,2}\}$. It suffices to show that $\Gamma_{1,2}\in\Gamma_{1,s-1}\Gamma_{s-1,1}$ when $\Gamma_{1,2}\in\Gamma_{1,1}^{2}$.

Suppose $\Gamma_{1,2}\notin\Gamma_{1,s-1}\Gamma_{s-1,1}$ and $\Gamma_{1,2}\in\Gamma_{1,1}^{2}$. Note that $(1,2)$ is mixed and $s>3$. By Theorem \ref{Main2} and Lemma \ref{both pure}, we have $\Gamma_{1,2}^{2}\cap\Gamma_{1,s-1}\Gamma_{s-1,1}\neq\emptyset$. Lemma \ref{jiben3} implies $p_{(1,s-1),(s-1,1)}^{(1,1)}=0$. By Lemma \ref{(1,q-1),(1,p-1)uniform}, we have $\Gamma_{1,2}\Gamma_{1,s-1}=\Gamma_{1,1}\Gamma_{1,s-1}$. Let $x,y,y',z$ be vertices such that $y\in P_{(1,2),(1,s-1)}(x,z)$ and $y'\in P_{(1,1),(1,s-1)}(x,z)$. Since $p_{(1,1),(1,1)}^{(1,2)}=k_{1,1}$, one gets $(y,y')\in\Gamma_{1,1}$, contrary to $p_{(1,s-1),(s-1,1)}^{(1,1)}=0$.$\qed$

Now we are ready to give a proof of Proposition \ref{fenqingkuang}.

\noindent\textbf{Proof of Proposition \ref{fenqingkuang}.}~(i) Suppose that C1 does not hold. By Theorem \ref{Main1}, $(1,q-1)$ is pure or one of the configurations C$(q)$ and D($q$) exists. In view of Theorems \ref{Main1}, \ref{Main2} and Lemma \ref{both pure}, if $(1,q-1)$ is pure or D($q$) exists, then $\Gamma_{1,q-1}^{2}\cap\Gamma_{1,q-2}\Gamma_{q-2,1}\neq\emptyset$ and $\Gamma_{1,t-1}^{2}\cap\Gamma_{1,q-1}\Gamma_{q-1,1}\neq\emptyset$ for all $t\in T\setminus\{q-1,q\}$, which imply that C2 holds by Lemma \ref{tongyi}; if C$(q)$ exists, then C3 holds.

(ii) Let $m=1$ if $(1,p-1)$ is pure, and $m=2$ otherwise. If $(1,p-1)$ is mixed, from (i) and Lemma \ref{both pure}, then C3 holds and $p=q$. Assume for the contrary, namely, $\Gamma_{1,p-1}^{i}\cap F_{T\setminus I}\neq\emptyset$ for some $i\in\{1,2,\ldots,mp-2m+1\}$. Let $l$ be the minimum integer such that there exists a path from $x_{0}$ to $x_{l}$ of length $l$ in the subdigraph $\Delta_{T\setminus I}$ with $(x_{0},x_{l})\in\Gamma_{\wz{j}}$ for some $\Gamma_{\wz{j}}\in\Gamma_{1,p-1}^{i}$. Without loss of generality, we may assume that $(x_{0}=y_{0},y_{1},\ldots,y_{i}=x_{l})$ is a path consisting of arcs of type $(1,p-1)$ and $(x_{0},x_{1},\ldots,x_{l})$ is a path with $h_{j}\in T\setminus I$ for all $0\leq j\leq l-1$, where $\partial(x_{j+1},x_{j})+1=h_{j}$.

By the minimality of $l$, we have $(x_{1},y_{1})\notin\Gamma_{1,p-1}$. Then $p_{(1,p-1),(p-1,1)}^{(1,h_{0}-1)}=0$. Similarly, $p_{(1,p-1),(p-1,1)}^{(1,h_{j}-1)}=0$ for all $0\leq j\leq l-1$. Suppose that the path $(x_{0},x_{1},\ldots,x_{l})$ contains at least two types of arcs. Without loss of generality, we may assume $h_{0}\neq h_{1}$. Note that $(1,p-1)$ is pure or $p-1\notin\{h_{1},h_{2}\}$. Since $\Gamma_{1,h_{j}-1}^{2}\cap\Gamma_{1,p-1}\Gamma_{p-1,1}\neq\emptyset$ for $j=0,1$, we get $\Gamma_{1,h_{0}-1}\Gamma_{1,h_{1}-1}\subseteq\Gamma_{1,p-1}\Gamma_{p-1,1}$ from Lemma \ref{Gamma 1,q-1,Gamma 1,p-1 subseteq Gamma 1,p-1Gamma p-1,1} (i). It follows that $y_{1}\in P_{(1,p-1),(p-1,1)}(x_{0},x_{2})$, contrary to the minimality of $l$. Then the path $(x_{0},x_{1},\ldots,x_{l})$ consists of arcs of one type.

If $l\geq3$, then there exists $x_{2}'\in P_{(1,h_{0}-1),(1,h_{0}-1)}(x_{1},x_{3})$ such that $(x_{2}',y_{1})\in\Gamma_{1,p-1}$ since $\Gamma_{1,h_{0}-1}^{2}\cap\Gamma_{1,p-1}\Gamma_{p-1,1}\neq\emptyset$, contrary to the minimality of $l$. Since $(1,p-1)$ is pure or C$(p)$ exists, from Lemmas \ref{(1,q-1)}, \ref{Gamma{1,q-1}=2} and Lemma \ref{mix-(2,q-2) 1} (ii), one gets $l=2$.

If $i=1$, then $\Gamma_{1,h_{0}-1}^2\varsubsetneq\Gamma_{1,p-1}\Gamma_{p-1,1}$, and so $\partial(x_{0},y_{1})=2$ from Lemma \ref{(1,q-1),(1,1)}, a contradiction. Hence, $i>1$. If $(1,p-1)$ is pure, from Lemma \ref{(1,q-1)}, then $(x_{0},x_{2})\in\Gamma_{2,p-2}$ and $p_{(1,h_{0}-1),(1,h_{0}-1)}^{(2,p-2)}\neq0$, a contradiction. Then C3 holds and $p=q$. By Lemma \ref{mix-(2,q-2) 1} (i),(ii) and Lemma \ref{Gamma{1,q-1}=2}, we have $(x_{0},x_{2})\in\Gamma_{1,q-1}\cup\Gamma_{1,q-2}\cup\Gamma_{2,q-2}\cup\Gamma_{2,q-3}$. Since $h_{0}\notin\{q,q-1\}$, from Lemmas \ref{(1,q-1)}, \ref{mix-(2,q-2) 2} and Lemma \ref{mix-(2,q-2) 1} (ii), we have $\partial(x_{0},x_{2})\neq1$. Since $(1,q-2)$ is pure, from Lemma \ref{(1,q-1)}, one gets $p_{(1,h_{0}-1),(1,h_{0}-1)}^{(2,q-3)}=0$, which implies $(x_{0},x_{2})\in\Gamma_{2,q-2}$. By Lemma \ref{mix-(2,q-2) 1} (i), there exists $y_{1}'\in P_{(1,q-1),(1,q-2)}(x_{0},x_{2})$. It follows that $\Gamma_{1,h_{0}-1}^2\varsubsetneq\Gamma_{1,q-1}\Gamma_{q-1,1}$. Since $h_{0}\neq q-1$, from Lemma \ref{(1,q-1),(1,1)}, we have $\wz{\partial}(x_{0},x_{2})=(2,q-2)=(2,2)$ and $h_{0}=2$. Note that $(x_{0},y_{1}',x_{2},x_{1})$ is a circuit consisting of arcs of types $(1,1),(1,2)$ and $(1,3)$, contrary to Lemma \ref{mix-(2,q-2) 1} (iii).$\qed$

\subsection{Subdigraphs}

In this subsection, we discuss the subdigraph $\Delta_{T/I}$ under separate assumptions based on Proposition \ref{fenqingkuang} (i). Before proceeding it, we need three auxiliary lemmas.

\begin{lemma}\label{shortest path}
Let $J$ be a nonempty subset of $T$ and $s$ an integer in $J$ with $\Gamma_{1,i-1}^{2}\cap\Gamma_{1,s-1}\Gamma_{s-1,1}\neq\emptyset$ for all $i\in J\setminus\{s\}$. Suppose that $P$ is a shortest path of length $l$ between distinct vertices in the subdigraph $\Delta_{J}$. If $P$ contains an arc of type $(1,s-1)$ and an arc of type $(1,h-1)$ for some $h\in J\setminus\{s\}$, then $p_{(1,s-1),(s-1,1)}^{(1,h-1)}=0$, and $P$ contains $i$ arcs of type $(1,h-1)$ and $l-i$ arcs of type $(1,s-1)$, where $i=1$ or $2$.
\end{lemma}
\textbf{Proof.}~Without loss of generality, we may assume $P=(x_{0},x_{1},\ldots,x_{l})$ with $(x_{0},x_{1})\in\Gamma_{1,s-1}$ and $(x_{1},x_{2})\in\Gamma_{1,h-1}$. By the minimality of $l$, we have $(x_{0},x_{2})\notin\Gamma_{1,s-1}$. It follows that the first statement is valid.

Suppose $(x_{i},x_{i+1})\in\Gamma_{1,r-1}$ with $r\notin\{s,h\}$ for some $2\leq i\leq l-1$. Without loss of generality, we may assume $i=2$. Similarly, $p_{(1,s-1),(s-1,1)}^{(1,r-1)}=p_{(1,s-1),(s-1,1)}^{(1,h-1)}=0$. Observe $\Gamma_{1,h-1}^{2}\cap\Gamma_{1,s-1}\Gamma_{s-1,1}\neq\emptyset$ and $\Gamma_{1,r-1}^{2}\cap\Gamma_{1,s-1}\Gamma_{s-1,1}\neq\emptyset$. By Lemma \ref{both pure}, $(2,3)\notin\{(h,s),(r,s)\}$ or $(1,2)$ is pure. Lemma \ref{Gamma 1,q-1,Gamma 1,p-1 subseteq Gamma 1,p-1Gamma p-1,1} (i) implies $x_{0}\in P_{(s-1,1),(1,s-1)}(x_{1},x_{3})=\Gamma_{s-1,1}(x_{1})$, a contradiction. Thus, $(x_{i},x_{i+1})\in\Gamma_{1,h-1}\cup\Gamma_{1,s-1}$ for $2\leq i\leq l-1$.

Suppose that the path $P$ contains at least three arcs of type $(1,h-1)$. Without loss of generality, we may assume $(x_{2},x_{3}),(x_{3},x_{4})\in\Gamma_{1,h-1}$. Since  $\Gamma_{1,h-1}^{2}\cap\Gamma_{1,s-1}\Gamma_{s-1,1}\neq\emptyset$, there exists $x_{3}'\in P_{(1,h-1),(1,h-1)}(x_{2},x_{4})=\Gamma_{1,h-1}(x_{2})$ such that $(x_{0},x_{3}')\in\Gamma_{1,s-1}$, a contradiction. This proves the second statement of this lemma.$\qed$

\begin{lemma}\label{path in digraph}
Let $(1,s-1)$ be pure such that $\Gamma_{1,h-1}^{2}\cap\Gamma_{1,s-1}\Gamma_{s-1,1}\neq\emptyset$ for any $h\in T\setminus\{s,s+1\}$. If $\partial_{\Delta_{T\setminus\{s,s+1\}}}(x,y)>s$, then $(1,s-2)$ is pure and $\Gamma_{\wz{\partial}_{\Gamma}(x,y)}\in\Gamma_{1,s-2}\Gamma_{s-2,1}\Gamma_{1,h-1}^{2}$ for some $h\in T\setminus\{s,s+1\}$ with $\Gamma_{1,h-1}^{2}\cap\Gamma_{1,p-2}\Gamma_{p-2,1}\neq\emptyset$.
\end{lemma}
\textbf{Proof.}~In the subdigraph $\Delta_{T\setminus\{s,s+1\}}$, choose a shortest path $(x=x_{0},x_{1},\ldots,x_{l}=y)$. Let $(x_{i},x_{i+1})\in\Gamma_{1,h_{i}-1}$ for $0\leq i\leq l-1$. Lemma \ref{both pure} implies $h_{i}<s$ for $0\leq i\leq l-1$.

We claim that $(1,h_{i}-1)$ is mixed or $|\{j\mid h_{j}=h_{i}~\textrm{and}~0\leq j\leq l-1\}|<l-s+h_{i}$ for $0\leq i\leq l-1$. Suppose not. Without loss of generality, we may assume that $(1,h_{0}-1)$ is pure and $h_{0}=h_{1}=\cdots=h_{l-s+h_{0}-1}$. If $h_{0}=2$, then $l-s+2\leq2$ since $p_{(1,1),(1,1)}^{\wz{h}}=k_{1,1}$ for all $\Gamma_{\wz{h}}\in\Gamma_{1,1}^{2}$, contrary to $l>s$; if $h_{0}>2$, from Lemma \ref{(1,q-1)}, then $l=\partial_{\Delta_{T\setminus\{s,s+1\}}}(x,y)\leq s-h_{0}+\partial_{\Delta_{T\setminus\{s,s+1\}}}(x_{0},x_{l-s+h_{0}})\leq s-h_{0}+h_{0}=s$, a contradiction. Thus, our claim is valid.

Suppose $h_{0}=h_{1}=\cdots=h_{l-1}$. Since $l>s\geq3$, from Lemma \ref{jiben4}, C$(h_{0})$ does not exist. By the claim and Theorem \ref{Main1}, D$(h_{0})$ exists. Since $p_{(1,h_{0}-2),(h_{0}-2,1)}^{(1,h_{0}-1)}=k_{1,h_{0}-2}$, there exists a vertex $z\in\Gamma_{1,h_{0}-2}(x_{i})$ for all $0\leq i\leq l$. Lemma \ref{(1,q-1)} implies $\partial_{\Delta_{T\setminus\{s,s+1\}}}(x,y)\leq1+\partial_{\Delta_{T\setminus\{s,s+1\}}}(z,y)=h_{0}-1<s-1$, a contradiction.

Without loss of generality, we may assume that $h_{0}=\max\{h_{i}\mid0\leq i\leq l-1\}$ and $h_{1}=\max\{h_{i}\mid1\leq i\leq l-1~\textrm{and}~h_{i}\neq h_{0}\}$. Pick a vertex $x_{1}'\in P_{(1,h_{1}-1),(1,h_{0}-1)}(x_{0},x_{2})$. By Theorem \ref{Main2}, we consider two cases.

\textbf{Case 1.} $(h_{0},h_{1})\in\mathscr{B}$.

Note that $h_{1}=h_{0}-1$, C$(h_{0})$ exists and $(1,h_{0}-2)$ is pure. In view of Lemma \ref{jiben4}, the path $(x_{0},x_{1},\ldots,x_{l})$ contains at most two arcs of type $(1,h_{0}-1)$. By the claim, one has $|\{i\mid h_{i}=h_{0}-1~\textrm{and}~0\leq i\leq l-1\}|<l-s+h_{0}-1\leq l-2$. It follows that the path $(x_{0},x_{1},\ldots,x_{l})$ contains at least three types of arcs. Without loss of generality, we may assume $h_{2}=\max\{h_{t}\mid0\leq t\leq l-1, h_{t}\neq h_{0}, h_{t}\neq h_{1}\}$. Choose a vertex $x_{2}'\in P_{(1,h_{2}-1),(1,h_{0}-2)}(x_{1},x_{3})$.

By Theorem \ref{Main2} and Lemma \ref{both pure}, we have $\Gamma_{1,h-1}^{2}\cap\Gamma_{1,h_{0}-2}\Gamma_{h_{0}-2,1}\neq\emptyset$ for all $h\in\{i\in T\mid i\leq h_{2}\}$. In view of Lemma \ref{shortest path}, we have $p_{(1,h_{0}-2),(h_{0}-2,1)}^{(1,h_{2}-1)}=0$, $|\{i\mid h_{i}=h_{2}~\textrm{and}~2\leq i\leq l-1\}|\leq2$ and $h_{t}\in\{h_{0},h_{0}-1,h_{2}\}$ for all $t\in\{3,4,\ldots,l-1\}$.

By Lemma \ref{(1,q-1),(1,p-1)} (i), one has $(x_{1},x_{3})\in\Gamma_{2,h_{0}-1}$. If $\Gamma_{2,h_{0}-1}\in\Gamma_{1,h_{0}-1}^{2}$, then there exists $x_{2}''\in P_{(1,h_{0}-1),(1,h_{0}-1)}(x_{1},x_{3})=\Gamma_{1,h_{0}-1}(x_{1})$ such that $(x_{0},x_{2}'')\in\Gamma_{1,h_{0}-2}$, a contradiction. In view of Lemma \ref{mix-(2,q-2) 1} (ii), we get $\Gamma_{1,h_{0}-1}^{2}=\{\Gamma_{1,h_{0}-2}\}$, which implies $h_{i}\neq h_{0}$ for $3\leq i\leq l-1$.

Since the path $(x_{0},x_{1},\ldots,x_{l})$ contains at most $l-3$ arcs of type $(1,h_{0}-2)$ and at most $2$ arcs of type $(1,h_{2}-1)$, we may assume $h_{3}=h_{2}$ and $h_{i}=h_{0}-1$ for $4\leq i\leq l-1$.

Choose a vertex $x_{3}'\in P_{(1,h_{2}-1),(1,h_{0}-2)}(x_{2}',x_{4})$. Since $(1,h_{0}-2)$ is pure, from Lemma \ref{(1,q-1)}, one has $s+1\leq\partial_{\Delta_{T\setminus\{s,s+1\}}}(x,y)=3+\partial_{\Delta_{T\setminus\{s,s+1\}}}(x_{3}',x_{l})\leq h_{0}+2<s+2$. Hence, $l=h_{0}+2=s+1$ and $(x_{4},x_{l})\in\Gamma_{h_{0}-2,1}$. Since $(x_{0},x_{1},x_{2}',x_{3}',x_{4})$ is a path with $(x_{0},x_{1})\in\Gamma_{1,h_{0}-1}$ and $(x_{3}',x_{4})\in\Gamma_{1,h_{0}-2}$, there exists a path $(x_{0},y_{1},y_{2},y_{3},x_{4})$ such that $(y_{2},y_{3})\in\Gamma_{1,h_{0}-2}$ and $(y_{3},x_{4})\in\Gamma_{1,h_{0}-1}$. By $p_{(1,h_{0}-1),(1,h_{0}-1)}^{(1,h_{0}-2)}=k_{1,h_{0}-1}$, one has $(x_{l},y_{3}),(y_{2},x_{l})\in\Gamma_{1,h_{0}-1}$. It follows that $l\leq3$, contrary to $l>s\geq3$.

\textbf{Case 2.} $(h_{0},h_{1})\in\mathscr{C}\setminus\mathscr{B}$.

Since $\partial_{\Delta_{T\setminus\{s,s+1\}}}(x_{0},x_{2})=2$, we get $p_{(1,h_{0}-1),(h_{0}-1,1)}^{(1,h_{1}-1)}=p_{(1,h_{1}-1),(h_{1}-1,1)}^{(1,h_{0}-1)}=0$.

Suppose $\Gamma_{1,h_{0}-1}^{2}\cap\Gamma_{1,h_{1}-1}\Gamma_{h_{1}-1,1}\neq\emptyset$. By Lemma \ref{both pure}, $h_{1}=h_{0}-1$ and $(1,h_{0}-2)$ is pure. Theorem \ref{Main2} implies $\Gamma_{1,h_{i}-1}^{2}\cap\Gamma_{1,h_{0}-2}\Gamma_{h_{0}-2,1}\neq\emptyset$ for all $0\leq i\leq l-1$ with $h_{i}\neq h_{0}-1$. By Lemma \ref{shortest path}, the path $(x_{0},x_{1},\ldots,x_{l})$ contains at least $l-2$ arcs of type $(1,h_{0}-2)$. Then $|\{j\mid h_{j}=h_{0}-1~\textrm{and}~0\leq j\leq l-1\}|\geq l-2\geq l-s+h_{0}-1$, contrary to the claim.

Suppose $\Gamma_{1,h_{1}-1}^{2}\cap\Gamma_{1,h_{0}-1}\Gamma_{h_{0}-1,1}\neq\emptyset$. Since $h_{0}>h_{1}$, from Theorem \ref{Main2} and Lemma \ref{both pure}, one gets $\Gamma_{1,h_{i}-1}^{2}\cap\Gamma_{1,h_{0}-1}\Gamma_{h_{0}-1,1}\neq\emptyset$ for all $0\leq i\leq l-1$ with $h_{i}\neq h_{0}$. In view of Lemma \ref{shortest path}, the path $(x_{0},x_{1},\ldots,x_{l})$ consisting of $i$ arcs of type $(1,h_{1}-1)$ and $l-i$ arcs of type $(1,h_{0}-1)$, where $i=1$ or $2$. Without loss of generality, we may assume $h_{2}=h_{3}=\cdots=h_{l-2}=h_{0}$.

Since $l>s>h_{0}>h_{1}$, $l\geq5$. By Lemma \ref{jiben4}, C$(h_{0})$ does not exist. If D($h_{0}$) exists, then there exists a vertex $z$ such that $x_{1}',x_{2},x_{3},\ldots,x_{l-1}\in\Gamma_{h_{0}-2,1}(z)$, which implies $l=\partial_{\Delta_{T\setminus\{s,s+1\}}}(x,y)=\partial_{\Delta_{T\setminus\{s,s+1\}}}(x_{1}',x_{l-1})+2\leq3+\partial_{\Delta_{T\setminus\{s,s+1\}}}(z,x_{l-1})=h_{0}+1\leq s$ from Lemma \ref{(1,q-1)}, a contradiction. By Theorem \ref{Main1}, $(1,h_{0}-1)$ is pure.

In view of Lemma \ref{(1,q-1)}, we have $\partial_{\Delta_{T\setminus\{s,s+1\}}}(x_{1}',x_{l-1})\leq h_{0}$. Since $s<l=2+\partial_{\Delta_{T\setminus\{s,s+1\}}}(x_{1}',x_{l-1})\leq h_{0}+2\leq s+1$, we get $l=h_{0}+2=s+1$ and $h_{l-1}=h_{1}$. By Lemma \ref{(1,q-1)}, one has $\Gamma_{\wz{\partial}_{\Gamma}(x,y)}\in\Gamma_{1,s-2}\Gamma_{s-2,1}\Gamma_{1,h_{1}-1}^{2}$ with $\Gamma_{1,h_{1}-1}^{2}\cap\Gamma_{1,s-2}\Gamma_{s-2,1}\neq\emptyset$.

This completes the proof of this lemma.$\qed$

\begin{lemma}\label{thick}
Let $J\subseteq T$. If $\partial_{\Gamma}(x,y)=\partial_{\Gamma}(x',y')$ for all $x,y,x',y'\in V\Delta_{J}$ with $\wz{\partial}_{\Delta_{J}}(x,y)=\wz{\partial}_{\Delta_{J}}(x',y')$, then $\Delta_{J}$ is a thick weakly distance-regular digraph.
\end{lemma}
\textbf{Proof.}~If $\wz{\partial}_{\Gamma}(x,y)=\wz{\partial}_{\Gamma}(x',y')$ for $x,y,x',y'\in V\Delta_{J}$, from the weakly distance-regularity of $\Gamma$, then $\wz{\partial}_{\Delta_{J}}(x,y)=\wz{\partial}_{\Delta_{J}}(x',y')$. By the assumption, we have $\Gamma_{\wz{i}}(x)=\{y\mid\wz{\partial}_{\Delta_{J}}(x,y)=\wz{i}'\}$, where $\wz{i}=\wz{\partial}_{\Gamma}(x,y)$ for some $y$ with $\wz{\partial}_{\Delta_{J}}(x,y)=\wz{i}'$. The desired result follows.$\qed$

In the following, we divide the proof into two subsubsections according to separate assumptions based on Proposition \ref{fenqingkuang} (i).

\subsubsection{The cases C1 and C3}

\begin{prop}\label{digraph C1 4}
 If {\rm C1} or {\rm C3} holds, then $\Delta_{T\setminus I}$ is a thick weakly distance-regular digraph.
\end{prop}
\textbf{Proof.}~Let $p$ be the maximum integer such that $(1,p-1)$ is pure. If C1 holds, then $p=q$; if C3 holds, then $p=q-1$. In order to prove this proposition, we only need to show that $\Delta_{T\setminus I}$ is a thick weakly distance-regular digraph. In view of Lemma \ref{thick}, it suffices to show that $\partial_{\Gamma}(x_{0},y_{0})=\partial_{\Gamma}(x_{1},y_{1})$ when $\partial_{\Delta_{T\setminus I}}(x_{0},y_{0})=\partial_{\Delta_{T\setminus I}}(x_{0},y_{1})$ for $x_{0},y_{0},x_{1},y_{1}\in V\Delta_{T\setminus I}$. Suppose, to the contrary that $\partial_{\Gamma}(x_{0},y_{0})<\partial_{\Gamma}(x_{1},y_{1})\leq\partial_{\Delta_{T\setminus I}}(x_{0},y_{0})=\partial_{\Delta_{T\setminus I}}(x_{1},y_{1})$. It follows that, in $\Gamma$, any shortest path from $x_{0}$ to $y_{0}$ containing an arc of type $(1,p-1)$ or $(1,p)$. By Proposition \ref{fenqingkuang} (ii) and Lemmas \ref{(1,q-1)}--\ref{Gamma{1,q-1}=2}, we have $p\leq\partial_{\Gamma}(x_{0},y_{0})<\partial_{\Delta_{T\setminus I}}(x_{0},y_{0})$.

In view of Lemmas \ref{(1,q-1)} and \ref{path in digraph}, $\partial_{\Delta_{T\setminus I}}(x_{0},y_{0})=\partial_{\Delta_{T\setminus I}}(x_{1},y_{1})=p+1$, $\partial_{\Gamma}(x_{0},y_{0})=p$ and $(1,p-2)$ is pure. For each $i\in\{0,1\}$, there exist vertices $u_{i},v_{i},w_{i}$ such that $(u_{i},x_{i}),(u_{i},v_{i})\in\Gamma_{1,p-2}$ and $(v_{i},w_{i}),(w_{i},y_{i})\in\Gamma_{1,h_{i}-1}$ with $\Gamma_{1,h_{i}-1}^2\cap\Gamma_{1,p-2}\Gamma_{p-2,1}\neq\emptyset$ and $h_{i}<p-1$. Since $(u_{i},y_{i})\notin\Gamma_{1,p-2}$, we get $\Gamma_{1,h_{i}-1}^{2}\subsetneq\Gamma_{1,p-2}\Gamma_{p-2,1}\neq\emptyset$, and so $\Gamma_{1,h_{i}-1}\notin\Gamma_{1,p-2}\Gamma_{p-2,1}$ by Lemma \ref{jiben3}. In view of Lemma \ref{both pure} (ii), one has $h_{0}=h_{1}$. By Lemma \ref{(1,q-1),(1,1)}, one gets $\wz{\partial}_{\Gamma}(v_{0},y_{0})=\wz{\partial}_{\Gamma}(v_{1},y_{1})$.

Since $\partial_{\Gamma}(x_{0},y_{0})=p$, from Proposition \ref{fenqingkuang} (ii) and  Lemmas \ref{(1,q-1)}--\ref{Gamma{1,q-1}=2} again, any shortest path from $x_{0}$ to $y_{0}$ in $\Gamma$ consisting of arcs of type $(1,p-1)$, which implies $\Gamma_{\wz{\partial}_{\Gamma}(x_{0},y_{0})}\in\Gamma_{1,p-1}\Gamma_{p-1,1}$. Choose vertices $z_{0}\in P_{(1,p-1),(p-1,1)}(x_{0},y_{0})$, $z_{1}\in\Gamma_{1,p-1}(x_{1})$, $x'\in\Gamma_{p-2,1}(u_{1})$ and $x\in\Gamma_{p-2,1}(u_{0})$. Note that $\Gamma_{1,p-2}^{2}\cap\Gamma_{1,p-1}\Gamma_{p-1,1}\neq\emptyset$. Since $(1,p-2)$ is pure, from Lemma \ref{(1,q-1)}, we get $p_{(1,p-1),(p-1,1)}^{(2,p-3)}=k_{1,p-1}$, $x,v_{0}\in\Gamma_{p-1,1}(z_{0})$ and $x',v_{1}\in\Gamma_{p-1,1}(z_{1})$. Since $z_{0}\in P_{(1,p-1),(p-1,1)}(v_{0},y_{0})$ and $\wz{\partial}_{\Gamma}(v_{0},y_{0})=\wz{\partial}_{\Gamma}(v_{1},y_{1})$,
we get $z_{1}\in P_{(1,p-1),(p-1,1)}(v_{1},y_{1})$, contrary to $\wz{\partial}_{\Gamma}(x_{0},y_{0})<\wz{\partial}_{\Gamma}(x_{1},y_{1})$. $\qed$

\subsubsection{The case C2}

In this subsubsection, we show that $\Delta_{T\setminus I}$ is a thick weakly distance-regular digraph. By Lemma \ref{both pure}, $(1,q-2)$ is pure. In view of Theorem \ref{Main2},  $\Gamma_{1,h-1}^{2}\cap\Gamma_{1,q-1}\Gamma_{q-1,1}\neq\emptyset$ for all $h\in T\setminus\{q-1,q\}$.

Next, we give some auxiliary lemmas.

\begin{lemma}\label{partial_Gamma(x,y)}
Let $(x,y)\in\Gamma_{a,b}$. Suppose $\partial_{\Delta_{T\setminus I}}(x,y)>a$. Then one of the following holds:\vspace{-0.3cm}
\begin{itemize}
\item [{\rm (i)}] $a=q-1$ and $\Gamma_{a,b}\in\Gamma_{1,q-2}\Gamma_{q-2,1}$;\vspace{-0.3cm}

\item [{\rm (ii)}] $a=q$ and $\Gamma_{a,b}\in\Gamma_{1,q-2}\Gamma_{q-2,1}\Gamma_{1,h-1}$;\vspace{-0.3cm}

\item [{\rm (iii)}] $a=q+1$ and $\Gamma_{a,b}\in\Gamma_{1,q-2}\Gamma_{q-2,1}\Gamma_{1,h-1}^{2}$;\vspace{-0.3cm}
\end{itemize}
Here, $\Gamma_{1,h-1}\notin\Gamma_{1,q-2}\Gamma_{q-2,1}$ for some $h\in T\setminus I$.
\end{lemma}
\textbf{Proof.}~Let $(x=x_{0},x_{1},\ldots,x_{a}=y)$ be a shortest path in $\Gamma$. Since $\partial_{{T\setminus I}}(x,y)>a$, the path $(x_{0},x_{1},\ldots,x_{a})$ contains an arc of type $(1,q-2)$. By Lemma \ref{(1,q-1)} and Proposition \ref{fenqingkuang} (a), the path $(x_{0},x_{1},\ldots,x_{a})$ contains exactly $q-1$ arcs of type $(1,q-2)$. Without loss of generality, we may assume $(x_{i},x_{i+1})\in\Gamma_{1,q-2}$ for $0\leq i\leq q-2$. Lemma \ref{(1,q-1)} implies $(x_{1},x_{q-1})\in\Gamma_{q-2,1}$. If $a=q-1$, then (i) holds.

Now suppose $a>q-1$. Write $\partial_{\Gamma}(x_{q},x_{q-1})=h-1$. Since C2 holds, we have $\Gamma_{1,i-1}^{2}\cap\Gamma_{1,q-2}\Gamma_{q-2,1}\neq\emptyset$ for $i\in T\setminus I$. By Lemma \ref{shortest path}, we have $\Gamma_{1,h-1}\notin\Gamma_{1,q-2}\Gamma_{q-2,1}$, $a\in\{q,q+1\}$ and $(x_{a-1},x_{a})\in\Gamma_{1,h-1}$. Thus, (ii) or (iii) holds.$\qed$

\begin{lemma}\label{Gamma_1,q_n-1^2=Gamma_2,q_n-2}
If $\partial_{\Gamma}(x,y)>q-1$ for some $x,y\in V\Delta_{T\setminus I}$, then $\Gamma_{1,q-1}^{2}=\{\Gamma_{2,q-2}\}$.
\end{lemma}
\textbf{Proof.}~Suppose not. By Lemma \ref{(1,q-1)}, D($q$) exists. In view of Lemma \ref{mix-(2,q-2) 2}, we have $\Gamma_{1,q-1}^{2}=\{\Gamma_{2,q-1}\}$.

Suppose that $p_{(1,q-2),(q-2,1)}^{(1,h-1)}=0$ for some $h\in T\setminus\{q-1,q\}$. Since $\Gamma_{1,h-1}^{2}\cap\Gamma_{1,q-2}\Gamma_{q-2,1}\neq\emptyset$, from Lemma \ref{(1,q-1),(1,p-1)} (i), one has $\Gamma_{1,h-1}\Gamma_{1,q-2}=\{\Gamma_{2,q-1}\}$. Then there exist vertices $y_{0},y_{1},y_{1}',y_{2}$ such that $y_{1}\in P_{(1,h-1),(1,q-2)}(y_{0},y_{2})$ and $y_{1}'\in P_{(1,q-1),(1,q-1)}(y_{0},y_{2})$. Since $p_{(1,q-2),(q-2,1)}^{(1,q-1)}=k_{1,q-2}$, we get $(y_{1},y_{1}')\in\Gamma_{1,q-2}$, contrary to Corollary \ref{bkb}. Thus, $\Gamma_{1,q-2}\Gamma_{1,h-1}=\{\Gamma_{1,q-2}\}$ for all $h\in T\setminus I$.

Pick a path $(x=x_{0},x_{1},\ldots,x_{l}=y)$ in the digraph $\Delta_{T\setminus I}$ and a vertex $z\in\Gamma_{q-2,1}(x)$. It follows that $x_{i}\in\Gamma_{1,q-2}(z)$ for all $0\leq i\leq l$. Then $q-1<\partial_{\Gamma}(x,y)\leq\partial_{\Gamma}(x,z)+1=q-1$, a contradiction.$\qed$

\begin{lemma}\label{partial Delta}
Let $(x,y)\in\Gamma_{a,b}$. Suppose $\partial_{\Delta_{T\setminus I}}(x,y)>q-1$. If $\Gamma_{1,q-1}^{2}=\{\Gamma_{2,q-2}\}$, then one of the following holds: \vspace{-0.3cm}
\begin{itemize}
\item [{\rm (i)}] $\partial_{\Delta_{T\setminus I}}(x,y)=q$ and $\Gamma_{a,b}\in\Gamma_{1,q-1}\Gamma_{q-1,1}\cup\Gamma_{q-1,1}\Gamma_{1,h-1}\cup\Gamma_{q-2,2}\Gamma_{1,h-1}^{2}$;\vspace{-0.3cm}

\item [{\rm (ii)}] $\partial_{\Delta_{T\setminus I}}(x,y)=q+1$ and $\Gamma_{a,b}\in\Gamma_{1,q-1}\Gamma_{q-1,1}\Gamma_{1,h-1}\cup\Gamma_{q-1,1}\Gamma_{1,h-1}^{2}$;\vspace{-0.3cm}

\item [{\rm (iii)}] $\partial_{\Delta_{T\setminus I}}(x,y)=q+2$ and $\Gamma_{a,b}\in\Gamma_{1,q-1}\Gamma_{q-1,1}\Gamma_{1,h-1}^{2}$.\vspace{-0.3cm}
\end{itemize}
Here, $p_{(1,q-1),(q-1,1)}^{(1,h-1)}=0$ and $h\notin\{q-1,q\}$.
\end{lemma}
\textbf{Proof.}~Pick a shortest path $(x=x_{0},x_{1},\ldots,x_{l}=y)$ in the subdigraph $\Delta_{T\setminus I}$.

Suppose $(x_{i},x_{i+1})\notin\Gamma_{1,q-1}$ for all $0\leq i\leq l-1$. Then $\partial_{\Delta_{T\setminus\{q-1,q\}}}(x,y)=l>q-1$. By Lemma \ref{path in digraph}, $(1,q-3)$ is pure, $\Gamma_{a,b}\in\Gamma_{1,q-3}\Gamma_{q-3,1}\Gamma_{1,h-1}^{2}$ and $\Gamma_{1,h-1}^{2}\cap\Gamma_{1,q-3}\Gamma_{q-3,1}\neq\emptyset$ for some $h<q-2$. Lemma \ref{(1,q-1)} implies $\Gamma_{1,h-1}\notin\Gamma_{1,q-3}\Gamma_{q-3,1}$. By Lemma \ref{(1,q-1),(1,p-1)} (i), we have $\Gamma_{1,q-3}\Gamma_{1,h-1}=\{\Gamma_{2,q-2}\}$. Pick a path $(y_{0},y_{1},y_{2})$ such that $y_{1}\in P_{(1,q-3),(1,h-1)}(y_{0},y_{2})$. Since $(y_{0},y_{2})\in\Gamma_{2,q-2}$ and $p_{(1,q-1),(1,q-1)}^{(2,q-2)}=k_{1,q-1}$, there exists $y_{1}'\in P_{(1,q-1),(1,q-1)}(y_{0},y_{2})$. By $q-2>h\geq2$ and Lemma \ref{(1,q-1)}, there exists a path $(y_{2},y_{3},y_{4})$ consisting of arcs of type $(1,q-1)$ such that $(y_{0},y_{4})\in\Gamma_{4,q-4}$. Since $(y_{1},y_{3})\in\Gamma_{2,q-2}$, one gets $y_{2}
\in P_{(1,q-1),(1,q-1)}(y_{1},y_{3})=\Gamma_{q-1,1}(y_{3})$, a contradiction.

Since $\Gamma_{1,h-1}^{2}\cap\Gamma_{1,q-1}\Gamma_{q-1,1}\neq\emptyset$ for all $h\in T\setminus\{q-1,q\}$, from Lemma \ref{shortest path}, the path $(x_{0},x_{1},\ldots,x_{l})$ contains $i$ arcs of type $(1,h-1)$ and $l-i$ arcs of type $(1,q-1)$ for some $h\in T\setminus\{q-1,q\}$ and $i\in\{0,1,2\}$ with $p_{(1,q-1),(q-1,1)}^{(1,h-1)}=0$.

Without loss of generality, we may assume $(x_{j},x_{j+1})\in\Gamma_{1,q-1}$ for $0\leq j\leq l-i-1$. Since $\Gamma_{1,q-1}^{2}=\{\Gamma_{2,q-2}\}$, from Lemma \ref{(1,q-1)}, one has $l-i\leq q$. The fact $q\leq l$ implies $l=q$ and $i\in\{0,1,2\}$, $l=q+1$ and $i\in\{1,2\}$, or $l=q+2$ and $i=2$.

Suppose $l=q$. By Lemma \ref{(1,q-1)}, if $i=0$, then $(x_{1},x_{q})\in\Gamma_{q-1,1}$ and $\Gamma_{a,b}\in\Gamma_{1,q-1}\Gamma_{q-1,1}$; if $i=1$, then $(x_{0},x_{q-1})\in\Gamma_{q-1,1}$ and $(x_{q-1},x_{q})\in\Gamma_{1,h-1}$, which imply $\Gamma_{a,b}\in\Gamma_{q-1,1}\Gamma_{1,h-1}$; if $i=2$, then $(x_{0},x_{q-2})\in\Gamma_{q-2,2}$ and $(x_{q-2},x_{q-1}),(x_{q-1},x_{q})\in\Gamma_{1,h-1}$, which imply $\Gamma_{a,b}\in\Gamma_{q-2,2}\Gamma_{1,h-1}^{2}$. Thus, (i) holds.

Suppose $l=q+1$. By Lemma \ref{(1,q-1)}, if $i=1$, then $(x_{1},x_{q})\in\Gamma_{q-1,1}$ and $(x_{q},x_{q+1})\in\Gamma_{1,h-1}$, which imply $\Gamma_{a,b}\in\Gamma_{1,q-1}\Gamma_{q-1,1}\Gamma_{1,h-1}$; if $i=2$, then $(x_{0},x_{q-1})\in\Gamma_{q-1,1}$ and $(x_{q-1},x_{q}),(x_{q},x_{q+1})\in\Gamma_{1,h-1}$, which imply $\Gamma_{a,b}\in\Gamma_{q-1,1}\Gamma_{1,h-1}^{2}$. Thus, (ii) holds.

Suppose $l=q+2$. Since $i=2$, from Lemma \ref{(1,q-1)}, we get $(x_{1},x_{q})\in\Gamma_{q-1,1}$ and $(x_{q},x_{q+1}),(x_{q+1},x_{q+2})\in\Gamma_{1,h-1}$. Thus, $\Gamma_{a,b}\in\Gamma_{1,q-1}\Gamma_{q-1,1}\Gamma_{1,h-1}^{2}$ and (iii) holds.$\qed$

\begin{lemma}\label{Gamma 3,q}
If $\Gamma_{1,h-1}^{2}\varsubsetneq\Gamma_{1,q-2}\Gamma_{q-2,1}$ for some $h\in T\setminus I$, then $\Gamma_{3,q}\in\Gamma_{1,h-1}^{2}\Gamma_{1,q-2}$ and $\Gamma_{4,q-1}\in\Gamma_{1,h-1}^{2}\Gamma_{1,q-2}^{2}$.
\end{lemma}
\textbf{Proof.}~Since C2 holds, from Lemmas \ref{(1,q-1)} and \ref{mix-(2,q-2) 2}, we have $\Gamma_{1,q-1}^{2}\subseteq\Gamma_{1,q-2}\Gamma_{q-2,1}$. It follows that $q-1>h\geq2$.

Pick distinct vertices $x,y,z,w$ such that $(x,y),(y,z)\in\Gamma_{1,h-1}$, $(z,w)\in\Gamma_{1,q-2}$ and $(x,w)\notin\Gamma_{1,q-2}$. Let $l$ be the minimum integer such that $\Gamma_{i,l}\in\Gamma_{\wz{\partial}_{\Gamma}(x,z)}\Gamma_{1,q-2}$. Without loss of generality, we may assume $\partial_{\Gamma}(w,x)=l$. Choose $u\in\Gamma_{h-1,1}(x)$ with $u\neq y$. Since $\Gamma_{1,h-1}^{2}\cap\Gamma_{1,q-2}\Gamma_{q-2,1}\neq\emptyset$, there exists $x'\in P_{(1,h-1),(1,h-1)}(u,y)$ such that $(x',w)\in\Gamma_{1,q-2}$. By the assumption and Lemma \ref{jiben3}, we get $p_{(1,q-2),(q-2,1)}^{(1,h-1)}=0$. Since $(1,q-2)$ is pure, by Lemma \ref{(1,q-1),(1,p-1)} (i), one has $\Gamma_{1,h-1}\Gamma_{1,q-2}=\{\Gamma_{2,q-1}\}$, which implies $u,y\in\Gamma_{q-1,2}(w)$. Then $q-2\leq\partial_{\Gamma}(w,x)=l\leq1+\partial_{\Gamma}(u,w)=q$.

Assume $l\in\{q-2,q-1\}$. Choose a path $(w=y_{0},y_{1},\ldots,y_{l}=x)$. Suppose $(y_{0},y_{1})\in\Gamma_{1,h'-1}$ with $h'\neq q-1$. By the minimality of $l$, we have $(z,y_{1})\notin\Gamma_{1,q-2}$, and so $p_{(1,q-2),(q-2,1)}^{(1,h'-1)}=0$.  Since C2 holds, one gets $\Gamma_{1,h'-1}^{2}\cap\Gamma_{1,q-2}\Gamma_{q-2,1}\neq\emptyset$. Note that $\Gamma_{1,h-1}^{2}\cap\Gamma_{1,q-2}\Gamma_{q-2,1}\neq\emptyset$ and $\Gamma_{1,h-1}\notin\Gamma_{1,q-2}\Gamma_{q-2,1}$. Lemma \ref{Gamma 1,q-1,Gamma 1,p-1 subseteq Gamma 1,p-1Gamma p-1,1} (ii) implies $h=h'$. Pick vertices $w'\in P_{(1,h-1),(1,q-2)}(z,y_{1})$ and $y'\in P_{(1,h-1),(1,h-1)}(x,z)$ such that $(y',y_{1})\in\Gamma_{1,q-2}$. Since $\Gamma_{1,h-1}\Gamma_{1,q-2}=\{\Gamma_{2,q-1}\}$, one has $(x,y_{1})\in\Gamma_{2,q-1}$, contrary to $\partial_{\Gamma}(y_{1},x)=l-1<q-1$. Hence, $(y_{i},y_{i+1})\in\Gamma_{1,q-2}$ for $0\leq i\leq l-1$. Since $(x,w)\notin\Gamma_{1,q-2}$, from Lemma \ref{(1,q-1)}, we get $l=q-1$ and $(z,x)\in\Gamma_{1,q-2}$. The minimality of $l$ implies $l=0$ and $q=1$, a contradiction. Thus, $\partial_{\Gamma}(w,x)=l=q$.

Pick a vertex $v\in P_{(h-1,1),(1,q-2)}(u,w)$. Since $(1,q-2)$ is pure, from Lemma \ref{(1,q-1)}, there exists a path $(w=z_{0},z_{1},\ldots,z_{q-2}=v)$ consisting of arcs of type $(1,q-2)$. By $(x,w)\notin\Gamma_{1,q-2}$, we have $P_{(1,q-2),(q-2,1)}(v,x)=\emptyset$. In view of Lemma \ref{(1,q-1),(1,1)}, one gets $\wz{\partial}_{\Gamma}(v,x)=\wz{\partial}_{\Gamma}(x,z)$. Since $(x,y,z,w=z_{0},z_{1},\ldots,z_{q-3})$ is a path, from the minimality of $l$, we get $\partial_{\Gamma}(x,z_{q-3})=q$, which implies $(x,w)\in\Gamma_{3,q}$ and $\Gamma_{3,q}\in\Gamma_{1,h-1}^{2}\Gamma_{1,q-2}$. By $q>3$, $(x,z_{1})\in\Gamma_{4,q-1}$. Thus, $\Gamma_{4,q-1}\in\Gamma_{1,h-1}^{2}\Gamma_{1,q-2}^{2}$.$\qed$

\begin{lemma}\label{b leq 3}
If $\Gamma_{a,b}\in\Gamma_{q-1,1}^{i}\Gamma_{1,h-1}^{j}$ for some $h\notin\{q-1,q\}$ and $i,j\in\{1,2\}$, then $\Gamma_{a,b}\in\Gamma_{q-1,1}^{i}\Gamma_{h-1,1}^{j}$.
\end{lemma}
\textbf{Proof.}~Let $(x_{i},y_{j})\in\Gamma_{a,b}$. Pick paths $(x_{0},x_{1},\ldots,x_{i})$ consisting of arcs of type $(1,q-1)$ and $(x_{0}=y_{0},y_{1},\ldots,y_{j})$ consisting of arcs of type $(1,h-1)$. Since $\Gamma_{1,h-1}^{2}\cap\Gamma_{1,q-1}\Gamma_{q-1,1}\neq\emptyset$, there exists a vertex $x_{0}'\in P_{(1,h-1),(1,q-1)}(y_{1},x_{1})$. If $j=1$, then the desired result holds.

Suppose $j=2$. By $y_{1}\in P_{(h-1,1),(1,h-1)}(x_{0}',y_{2})$, there exists a vertex $y_{3}\in P_{(1,h-1),(h-1,1)}(x_{0}',y_{2})$. Since $\Gamma_{1,h-1}^{2}\cap\Gamma_{1,q-1}\Gamma_{q-1,1}\neq\emptyset$ again, there also exists a vertex $x_{0}''\in P_{(1,h-1),(1,q-1)}(y_{3},x_{1})$, which implies $\Gamma_{a,b}\in\Gamma_{q-1,1}^{i}\Gamma_{h-1,1}^{2}$.$\qed$

\begin{lemma}\label{b>3}
Let $\Gamma_{1,q-1}^{2}=\{\Gamma_{2,q-2}\}$. Suppose $p_{(1,q-2),(q-2,1)}^{(1,q-1)}=0$ or $p_{(1,q-2),(q-2,1)}^{(1,h-1)}=0$, and $p_{(1,q-1),(q-1,1)}^{(1,h-1)}=0$ for some $h\in T\setminus\{q-1,q\}$. If there exists $(x,y)\in\Gamma_{q-1,b}$ such that $\partial_{\Delta_{T\setminus I}}(x,y)\geq q$, then $b>3$.
\end{lemma}
\textbf{Proof.}~Suppose not. Since $q\leq\partial_{\Delta_{T\setminus I}}(x,y)$ and $\Gamma_{1,q-1}^{2}=\{\Gamma_{2,q-2}\}$, from Lemma \ref{(1,q-1)}, we have $b=2$ or $3$. By Lemma \ref{partial_Gamma(x,y)} (i), there exists $z\in P_{(1,q-2),(q-2,1)}(x,y)$.

Suppose $b=2$. Note that $\Gamma_{1,q-1}$ or $\Gamma_{1,h-1}\notin\Gamma_{1,q-2}\Gamma_{q-2,1}$. Since C2 holds, from Lemma \ref{(1,q-1),(1,p-1)} (i), we have $\Gamma_{1,q-1}\Gamma_{1,q-2}=\{\Gamma_{2,q-1}\}$ or $\Gamma_{1,h-1}\Gamma_{1,q-2}=\{\Gamma_{2,q-1}\}$. It follows that there exists $w\in\Gamma_{1,q-2}(y)$ such that $(w,x)\in\Gamma_{1,q-1}\cup\Gamma_{1,h-1}$. By $P_{(1,q-2),(q-2,1)}(x,y)=\Gamma_{1,q-2}(y)$, one gets $(x,w)\in\Gamma_{1,q-2}$, contrary to $q-1\notin\{q,h\}$.

Suppose $b=3$. Since $\Gamma_{1,h-1}^{2}\cap\Gamma_{1,q-1}\Gamma_{q-1,1}\neq\emptyset$ and $p_{(1,q-1),(q-1,1)}^{(1,h-1)}=0$, from the assumption and Lemma \ref{(1,q-1),(1,p-1)}, one gets $\Gamma_{1,q-1}\Gamma_{1,h-1}=\{\Gamma_{2,q}\}$. By Lemma \ref{Gamma 2,q-1}, one gets $\Gamma_{1,q-1}^{2}\Gamma_{1,h-1}=\{\Gamma_{3,q-1}\}$, which implies that there exist vertices $u,v$ such that $(y,u),(u,v)\in\Gamma_{1,q-1}$ and $(v,x)\in\Gamma_{1,h-1}$. Since $\Gamma_{1,h-1}^{2}\cap\Gamma_{1,q-1}\Gamma_{q-1,1}\neq\emptyset$, there exists $w\in P_{(1,h-1),(q-1,1)}(x,u)$. By Lemma \ref{(1,q-1)}, there exists a path $(w=x_{0},x_{1},\ldots,x_{q-2}=y)$ consisting of arcs of type $(1,q-1)$. It follows that $(x,x_{0},x_{1},\ldots,x_{q-2}=y)$ is a path of length $q-1$ in the subdigraph $\Delta_{T\setminus I}$, contrary to $\partial_{\Delta_{T\setminus I}}(x,y)\geq q$.$\qed$

\begin{prop}\label{c2 c3}
The digraph $\Delta_{T\setminus I}$ is a thick weakly distance-regular digraph.
\end{prop}
\textbf{Proof.}~By Lemma \ref{thick}, it suffices to show that $\partial_{\Gamma}(x,y)=\partial_{\Gamma}(x',y')$ for any vertices $x,x',y,y'\in V\Delta_{T\setminus I}$ with $\wz{\partial}_{\Delta_{T\setminus I}}(x,y)=\wz{\partial}_{\Delta_{T\setminus I}}(x',y')$. Suppose, to the contrary that $a<a'$, where $(x,y)\in\Gamma_{a,b}$ and $(x',y')\in\Gamma_{a',b'}$. By Lemma \ref{partial_Gamma(x,y)}, we have $q-1\leq a<a'\leq\partial_{\Delta_{T\setminus I}}(x,y)$. Lemma \ref{Gamma_1,q_n-1^2=Gamma_2,q_n-2} implies $\Gamma_{1,q-1}^{2}=\{\Gamma_{2,q-2}\}$.

We claim $\Gamma_{1,q-1}\Gamma_{q-1,1}\subseteq\Gamma_{1,q-2}\Gamma_{q-2,1}$. Choose vertices $x_{0},x_{1},x_{2}$ such that $x_{1}\in P_{(1,q-1),(q-1,1)}(x_{0},x_{2})$. Let $x_{1}'\in\Gamma_{1,q-1}(x_{1})$ and $w\in\Gamma_{1,q-2}(x_{0})$. Since $\Gamma_{1,q-1}^{2}=\{\Gamma_{2,q-2}\}$, we have $x_{0},x_{2}\in\Gamma_{q-2,2}(x_{1}')$. Since C2 holds, one gets $\Gamma_{2,q-2}\in\Gamma_{1,q-2}\Gamma_{q-2,1}$, which implies $x_{1}',x_{2}\in\Gamma_{q-2,1}(w)$. Thus, our claim is valid.

Since $\partial_{\Delta_{T\setminus I}}(x,y)>q-1$, from Lemma \ref{partial Delta}, we consider three cases.

\textbf{Case 1}. $\partial_{\Delta_{T\setminus I}}(x,y)=\partial_{\Delta_{T\setminus I}}(x',y')=q+2$.

Pick vertices $u,u',v,v',w,w'$ such that $(u,x),(u',x'),(u,v),(u',v')\in\Gamma_{1,q-1}$, $(v,w),(w,y)\in\Gamma_{1,h-1}$ and $(v',w'),(w',y')\in\Gamma_{1,h'-1}$ for some $h,h'\in T\setminus\{q-1,q\}$. By Lemma \ref{(1,q-1)}, we have $(u,y),(u',y')\notin\Gamma_{1,q-1}$, which implies $\Gamma_{1,h-1}^{2},\Gamma_{1,h'-1}^{2}\subsetneq\Gamma_{1,q-1}\Gamma_{q-1,1}$. Note that $\Gamma_{1,h-1}^{2}\cap\Gamma_{1,q-1}\Gamma_{q-1,1}\neq\emptyset $ and $\Gamma_{1,h'-1}^{2}\cap\Gamma_{1,q-1}\Gamma_{q-1,1}\neq\emptyset$. In view of Lemma \ref{jiben3} and Lemma \ref{Gamma 1,q-1,Gamma 1,p-1 subseteq Gamma 1,p-1Gamma p-1,1} (ii), one has $h=h'$.

Since $\Gamma_{1,h-1}^{2}\subsetneq\Gamma_{1,q-1}\Gamma_{q-1,1}$, from Lemma \ref{(1,q-1),(1,1)}, we have $\wz{\partial}_{\Gamma}(v,y)=\wz{\partial}_{\Gamma}(v',y')$. Pick vertices $z,z'$ such that $(x,z),(x',z')\in\Gamma_{1,q-2}$. By the claim, one obtains $(v,z),(v',z')\in\Gamma_{1,q-2}$. Lemma \ref{partial_Gamma(x,y)} implies $q-1\leq a<a'\leq q+1$. Since $(y',z')\notin\Gamma_{1,q-2}$, we have $P_{(1,q-2),(q-2,1)}(v',y')=\emptyset$. The fact $\wz{\partial}_{\Gamma}(v,y)=\wz{\partial}_{\Gamma}(v',y')$ implies $(y,z)\notin\Gamma_{1,q-2}$, and so $\Gamma_{a,b}\notin\Gamma_{1,q-2}\Gamma_{q-2,1}$. By Lemma \ref{partial_Gamma(x,y)}, we have $a=q$.

By Lemma \ref{partial_Gamma(x,y)} (ii), there exist two vertices $y_{0},y_{1}$ such that $(x,y_{0}),(y_{1},y_{0})\in\Gamma_{1,q-2}$ and $(y_{1},y)\in\Gamma_{1,h''-1}$ for some $h''\in T\setminus I$ with $p_{(1,q-2),(q-2,1)}^{(1,h''-1)}=0$. Since $(y',z')\notin\Gamma_{1,q-2}$, we have $\Gamma_{1,h-1}^2\subsetneq\Gamma_{1,q-2}\Gamma_{q-2,1}$, which implies $p_{(1,q-2),(q-2,1)}^{(1,h-1)}=0$ from Lemma \ref{jiben3}. Since $\Gamma_{1,s-1}^{2}\cap\Gamma_{1,q-2}\Gamma_{q-2,1}\neq\emptyset$ for all $s\in T\setminus I$, from Lemma \ref{Gamma 1,q-1,Gamma 1,p-1 subseteq Gamma 1,p-1Gamma p-1,1} (ii), one gets $h''=h$. The fact that $P_{(1,h-1),(1,h-1)}(v,y)=\Gamma_{h-1,1}(y)$ implies $(v,y_{1})\in\Gamma_{1,h-1}$. By the claim, we obtain $(v,y_{0})\in\Gamma_{1,q-2}$, contrary to $p_{(1,q-2),(q-2,1)}^{(1,h-1)}=0$.

\textbf{Case 2}. $\partial_{\Delta_{T\setminus I}}(x,y)=\partial_{\Delta_{T\setminus I}}(x',y')=q+1$.

By Lemma \ref{partial Delta} (ii), we have $\Gamma_{a,b}\in\Gamma_{q-1,1}\Gamma_{1,q-1}\Gamma_{1,h-1}\cup\Gamma_{q-1,1}\Gamma_{1,h-1}^{2}$ for some $h\in T\setminus\{q-1,q\}$ with $\Gamma_{1,h-1}\notin\Gamma_{1,q-1}\Gamma_{q-1,1}$, and $\Gamma_{a',b'}\in\Gamma_{q-1,1}\Gamma_{1,q-1}\Gamma_{1,h'-1}\cup\Gamma_{q-1,1}\Gamma_{1,h'-1}^{2}$ for some $h'\in T\setminus\{q-1,q\}$ with $\Gamma_{1,h'-1}\notin\Gamma_{1,q-1}\Gamma_{q-1,1}$. Note that $\Gamma_{1,h-1}^{2}\cap\Gamma_{1,q-1}\Gamma_{q-1,1}\neq\emptyset$ and $\Gamma_{1,h'-1}^{2}\cap\Gamma_{1,q-1}\Gamma_{q-1,1}\neq\emptyset$. Since $h,h'<q-1$, from Lemmas \ref{(1,q-1)} and \ref{(1,q-1),(1,p-1)uniform}, we get $\Gamma_{q-1,1}\Gamma_{1,q-1}\Gamma_{1,h-1}=\Gamma_{q-1,1}\Gamma_{1,q-1}\Gamma_{1,h'-1}$ and $\Gamma_{q-1,1}\Gamma_{1,h-1}^{2}=\Gamma_{q-1,1}\Gamma_{1,h'-1}^{2}$. Without loss of generality, we may assume $h=h'$.

Suppose $\Gamma_{a,b},\Gamma_{a',b'}\in\Gamma_{1,q-1}\Gamma_{q-1,1}\Gamma_{1,h-1}$. Pick vertices $u,u',v,v',z,z'$ such that $(x,u),(v,u),(x',u'),(v',u')\in\Gamma_{1,q-1}$, $(v,y),(v',y')\in\Gamma_{1,h-1}$ and $(x,z),(x',z')\in\Gamma_{1,q-2}$. By the claim, we get $(v,z),(v',z')\in\Gamma_{1,q-2}$. Then $a=q-1$ and $a'=q$. Lemma \ref{partial_Gamma(x,y)} (i) implies $z\in P_{(1,q-2),(q-2,1)}(x,y)=\Gamma_{1,q-2}(x)$. Hence, $p_{(1,q-2),(q-2,1)}^{(1,h-1)}=k_{1,q-2}$ and $(y',z')\in\Gamma_{1,q-2}$, contrary to $a'=q$. Then $\Gamma_{a,b}$ or $\Gamma_{a',b'}\in\Gamma_{q-1,1}\Gamma_{1,h-1}^{2}$.

By Lemma \ref{b leq 3}, we have $\partial_{\Delta_{T\setminus I}}(y,x)=\partial_{\Delta_{T\setminus I}}(y',x')\leq3$, and so $b,b'\leq3$. If $\Gamma_{1,q-1},\Gamma_{1,h-1}\in\Gamma_{1,q-2}\Gamma_{q-2,1}$, then $\Gamma_{a,b},\Gamma_{a',b'}\in\Gamma_{1,q-2}\Gamma_{q-2,1}$ and $a=a'=q-1$, a contradiction. Then $\Gamma_{1,q-1}$ or $\Gamma_{1,h-1}\notin\Gamma_{1,q-2}\Gamma_{q-2,1}$. Since $\Gamma_{1,h-1}\notin\Gamma_{1,q-1}\Gamma_{q-1,1}$ and $\Gamma_{1,q-1}^{2}=\{\Gamma_{2,q-2}\}$, from Lemma \ref{b>3}, one has $a=q$ and $a'=q+1$.

If $\Gamma_{a',b'}\in\Gamma_{1,q-1}\Gamma_{q-1,1}\Gamma_{1,h-1}$, from the claim, then $\Gamma_{a',b'}\in\Gamma_{1,q-2}\Gamma_{q-2,1}\Gamma_{1,h-1}$, contrary to $a'=q+1$. Then $\Gamma_{a',b'}\in\Gamma_{q-1,1}\Gamma_{1,h-1}^{2}$. If $\Gamma_{1,h-1}^{2}\subseteq\Gamma_{1,q-2}\Gamma_{q-2,1}$, then $\Gamma_{a',b'}\in\Gamma_{q-1,1}\Gamma_{q-2,1}\Gamma_{1,q-2}$, which implies $\Gamma_{a',b'}\in\Gamma_{q-2,1}\Gamma_{1,q-2}\cup\Gamma_{q-1,2}\Gamma_{1,q-2}$ from Lemma \ref{(1,q-1),(1,p-1)} (i), contrary to $a'=q+1$. Hence, $\Gamma_{1,h-1}^{2}\subsetneq\Gamma_{1,q-2}\Gamma_{q-2,1}$.

By Lemma \ref{partial_Gamma(x,y)} (ii), there exist vertices $z,w$ such that $(z,x),(z,w)\in\Gamma_{1,q-2}$ and $(w,y)\in\Gamma_{1,h''-1}$ with $\Gamma_{1,h''-1}\notin\Gamma_{1,q-2}\Gamma_{q-2,1}$ for some $h''\notin T\setminus I$. Since $\Gamma_{1,h-1}^{2}\subsetneq\Gamma_{1,q-2}\Gamma_{q-2,1}$, from Lemma \ref{jiben3}, we have $\Gamma_{1,h-1}\notin\Gamma_{1,q-2}\Gamma_{q-2,1}$. In view of Lemma \ref{Gamma 1,q-1,Gamma 1,p-1 subseteq Gamma 1,p-1Gamma p-1,1} (ii), one gets $h=h''$. Observe that $(a,b)=(q,2)$ or $(q,3)$.

Suppose $(a,b)=(q,2)$. Note that $\Gamma_{1,h-1}^{2}\cap\Gamma_{1,q-1}\Gamma_{q-1,1}\neq\emptyset$ and $\Gamma_{1,h-1}\notin\Gamma_{1,q-1}\Gamma_{q-1,1}$. Since $\Gamma_{1,q-1}$ or $\Gamma_{1,h-1}\notin\Gamma_{1,q-2}\Gamma_{q-2,1}$, from Lemma \ref{(1,q-1),(1,p-1)}, we have $\Gamma_{1,q-1}\Gamma_{1,h-1}=\{\Gamma_{2,q}\}$. Then there exist vertices $u\in P_{(1,h-1),(1,q-1)}(y,x)$ and $v\in P_{(q-1,1),(1,h-1)}(x,y)$. Since $\Gamma_{1,q-1}^{2}=\{\Gamma_{2,q-2}\}$, from Lemma \ref{(1,q-1)}, there exists a path $(x_{0}=x,x_{1},\ldots,x_{q-1}=v)$ consisting of arcs of type $(1,q-1)$. It follows that $(x=x_{0},x_{1},\ldots,x_{q-1},y)$ is a path of length $q$ in the subdigraph $\Delta_{T\setminus I}$, contrary to $\partial_{\Delta_{T\setminus I}}(x,y)=q+1$.

Suppose $(a,b)=(q,3)$. Since $\Gamma_{1,h-1}^{2}\subsetneq\Gamma_{1,q-2}\Gamma_{q-2,1}$, from Lemma \ref{Gamma 3,q}, there exist vertices $u,v$ such that $(y,u),(u,v)\in\Gamma_{1,h-1}$ and $(v,x)\in\Gamma_{1,q-2}$. The fact $\Gamma_{1,h-1}^{2}\cap\Gamma_{1,q-2}\Gamma_{q-2,1}\neq\emptyset$ implies that there exists $y_{0}\in P_{(1,h-1),(1,h-1)}(w,u)=\Gamma_{1,h-1}(w)$ such that $(y_{0},x)\in\Gamma_{1,q-2}$. Since $P_{(q-2,1),(1,q-2)}(x,w)=\Gamma_{q-2,1}(x)$ implies $(y_{0},w)\in\Gamma_{1,q-2}$, contrary to $q-1\neq h$.

\textbf{Case 3}. $\partial_{\Delta_{T\setminus I}}(x,y)=\partial_{\Delta_{T\setminus I}}(x',y')=q$.

Since $q-1\leq a<a'\leq \partial_{\Delta_{T\setminus I}}(x',y')$, we have $a=q-1$, $a'=q$. If $\Gamma_{q,b'}\in\Gamma_{1,q-1}\Gamma_{q-1,1}$, from the claim, then $\Gamma_{q,b'}\in\Gamma_{1,q-2}\Gamma_{q-2,1}$, a contradiction. Lemma \ref{(1,q-1)} and Lemma \ref{partial Delta} (i) imply $\Gamma_{q,b'}\in\Gamma_{q-1,1}\Gamma_{1,h-1}\cup\Gamma_{q-2,2}\Gamma_{1,h-1}^{2}$ for some $h\in T\setminus\{q-1,q\}$ with $p_{(1,q-1),(q-1,1)}^{(1,h-1)}=0$.

Suppose $\Gamma_{1,q-1},\Gamma_{1,h-1}\in\Gamma_{1,q-2}\Gamma_{q-2,1}$. Since $\Gamma_{1,q-1}^{2}=\{\Gamma_{2,q-2}\}$, we have $\Gamma_{q,b'}\in\Gamma_{1,q-2}\Gamma_{q-2,1}\Gamma_{1,h-1}\cup\Gamma_{1,q-2}\Gamma_{q-2,1}\Gamma_{1,h-1}^{2}$, which implies $\Gamma_{q,b'}\in\Gamma_{1,q-2}\Gamma_{q-2,1}$, contrary to $a'=q$. Hence, $p_{(1,q-2),(q-2,1)}^{(1,q-1)}=0$ or $p_{(1,q-2),(q-2,1)}^{(1,h-1)}=0$.

Since $\Gamma_{1,q-1}^{2}=\{\Gamma_{2,q-2}\}$, from Lemma \ref{b>3}, we get $b>3$. Since $\partial_{\Delta_{T\setminus I}}(y,x)=\partial_{\Delta_{T\setminus I}}(y',x')\geq b>3$, from Lemma \ref{b leq 3}, we have $\Gamma_{q,b'}\notin\Gamma_{q-1,1}\Gamma_{1,h-1}$. By $\Gamma_{1,q-1}^2=\{\Gamma_{2,q-2}\}$, one obtains $\Gamma_{q,b'}\in\Gamma_{q-2,2}\Gamma_{1,h-1}^{2}=\Gamma_{q-1,1}^{2}\Gamma_{1,h-1}^{2}$ and $b=\partial_{\Delta_{T\setminus I}}(y,x)=4$.

If $\Gamma_{1,h-1}^{2}\subseteq\Gamma_{1,q-2}\Gamma_{q-2,1}$, then $\Gamma_{q,b'}\in\Gamma_{1,q-2}\Gamma_{q-2,1}$ since $\Gamma_{1,q-1}^{2}=\{\Gamma_{2,q-2}\}\subseteq\Gamma_{1,q-2}\Gamma_{q-2,1}$, contrary to $a'=q$. Hence, $\Gamma_{1,h-1}^{2}\subsetneq\Gamma_{1,q-2}\Gamma_{q-2,1}$. By Lemma \ref{Gamma 3,q}, we have $\Gamma_{4,q-1}\in\Gamma_{1,h-1}^{2}\Gamma_{1,p-2}^{2}$.

By Lemma \ref{partial_Gamma(x,y)} (i), there exists $z\in P_{(1,q-2),(q-2,1)}(x,y)$. Since $(x,y)\in\Gamma_{q-1,4}$, there exists a path $(y,u,v,w,x)$ such that $(y,u),(u,v)\in\Gamma_{1,q-2}$ and $(v,w),(w,x)\in\Gamma_{1,h-1}$. The fact $P_{(1,q-2),(q-2,1)}(x,y)=\Gamma_{1,q-2}(y)$ implies $(x,u)\in\Gamma_{1,q-2}$. Since $(x,u,v,w)$ is a circuit and $(1,q-2)$ is pure, one gets $q=4$ and $h=2$. By Lemma \ref{(1,q-1)}, we obtain $\Gamma_{1,2}^{2}=\{\Gamma_{2,1}\}$ and $(v,x)\in\Gamma_{1,2}$, contrary to $\partial_{\Gamma}(y,x)=b=4$.$\qed$

\subsection{Quotient digraphs}

In this subsection, we determine quotient digraphs of $\Gamma$ over $F_{T/I}$ under two separate assumptions based on Proposition \ref{fenqingkuang} (i).

\begin{prop}\label{quotient 1}
If {\rm C1} or {\rm C2} holds, then $\Gamma/F_{T\setminus I}$ is isomorphic to one of the digraphs in Theorem {\rm\ref{Main4} (i)--(ii)}.
\end{prop}
\textbf{Proof.}~By Lemma \ref{both pure}, $(1,p-1)$ is pure. Lemma \ref{(1,q-1)} implies that there exists a circuit $(x_{0,0},x_{1,0},\ldots,x_{p-1,0})$ consisting of arcs of type $(1,p-1)$, where the first subscription of $x$ could be read modulo $p$. In view of Proposition \ref{fenqingkuang} (a), $(F_{T\setminus I}(x_{0,0}),F_{T\setminus I}(x_{1,0}),\ldots,F_{T\setminus I}(x_{p-1,0}))$ is a circuit in the quotient digraph $\Gamma/F_{T\setminus I}$. If $V\Gamma=\dot{\bigcup}_{i}F_{T\setminus I}(x_{i,0})$, then $\Gamma/F_{T\setminus I}$ is isomorphic to the digraph in Theorem \ref{Main4} (i).

Suppose that $V\Gamma\neq\dot{\bigcup}_{i}F_{T\setminus I}(x_{i,0})$. By the weakly distance-regularity of $\Gamma$, there exists a vertex $x_{i,1}\in\Gamma_{1,p-1}(x_{i-1,0})$ such that $x_{i,1}\notin\dot{\bigcup}_{j}F_{T\setminus I}(x_{j,0})$ for each $i$. It follows that $x_{i,1}\in P_{(1,p-1),(1,p-1)}(x_{i-1,0},x_{i+1,0})=\Gamma_{1,p-1}(x_{i-1,0})$ and $x_{i,1}\in P_{(1,p-1),(1,p-1)}(x_{i-1,0},x_{i+1,1})=\Gamma_{1,p-1}(x_{i-1,0})$. Since $x_{i,1}\notin F_{T\setminus I}(x_{i,0})$, from Lemma \ref{(1,q-1)}, one gets $(x_{i,0},x_{i,1})\in\Gamma_{p,p}$. If there exists $x_{i,2}\in\Gamma_{1,p-1}(x_{i-1,0})$ such that $x_{i,2}\notin\dot{\bigcup}_{j}(F_{T\setminus I}(x_{j,0})\cup F_{T\setminus I}(x_{j,1}))$, then $(x_{i,1},x_{i,2}),(x_{i,2},x_{i,0})\in\Gamma_{p,p}$, which implies $x_{i,0}\in\Gamma_{p,p}(x_{i,1})=P_{(p,p),(p,p)}(x_{i,1},x_{i,0})$, a contradiction. Thus, $V\Gamma=\dot{\bigcup}_{i,j}F_{T\setminus I}(x_{i,j})$, and the quotient digraph $\Gamma/F_{T\setminus I}$ is isomorphic to the digraph in Theorem \ref{Main4} (ii).$\qed$

\begin{prop}\label{quotient 2}
If {\rm C3} holds, then $\Gamma/F_{T\setminus I}$ is isomorphic to one of the digraphs in Theorem {\rm\ref{Main4} (iii)--(vi)}.
\end{prop}
\textbf{Proof.}~Let $(x_{0,0},x_{0,1},x_{1,0},x_{1,1},\ldots,x_{q-2,0},x_{q-2,1})$ be a circuit consisting of arcs of type $(1,q-1)$ such that $(x_{i,j},x_{i+1,j})\in\Gamma_{1,q-2}$ for $0\leq i\leq q-2$ and $0\leq j\leq 1$, where the first subscription of $x$ could be read modulo $q-1$. By Proposition \ref{fenqingkuang} (b), $(F_{T\setminus I}(x_{0,0}),F_{T\setminus I}(x_{0,1}),\ldots,F_{T\setminus I}(x_{q-2,0}),F_{T\setminus I}(x_{q-2,1}))$ is a circuit in the quotient digraph $\Gamma/F_{T\setminus I}$ such that $\partial(F_{T\setminus I}(x_{i,j}),F_{T\setminus I}(x_{i+1,j}))=1$ for any $i,j$.

If $V\Gamma=\dot{\bigcup}_{i,j}F_{T\setminus I}(x_{i,j})$, then $\Gamma/F_{T\setminus I}$ is isomorphic to the digraph in Theorem \ref{Main4} (iii). We only need to consider the case $V\Gamma\neq\dot{\bigcup}_{i,j}F_{T\setminus I}(x_{i,j})$.

\textbf{Case 1.} $\partial(F_{T\setminus I}(x_{0,0}),F_{T\setminus I}(x))\neq2$ for all $x\in\Gamma_{1,q-1}(x_{0,1})$.

By the weakly distance-regularity of $\Gamma$, there exists a vertex $x_{i,j}'$ such that $x_{i,j}'\notin\dot{\bigcup}_{a,b}F_{T\setminus I}(x_{a,b})$ and $\partial(F_{T\setminus I}(x_{i,j}'),F_{T\setminus I}(x_{i+j,j+1}))=1$ for all $i,j$, where the second subscription of $x$ could be read modulo $2$ for this case. Since $p_{(1,q-1),(1,q-1)}^{(1,q-2)}=k_{1,q-1}$, we may assume $(x_{i,j}',x_{i+j,j+1})\in\Gamma_{1,q-1}$. It follows that $x_{i,j}'\in \Gamma_{q-1,1}(x_{i+j,j+1})=P_{(1,q-1),(1,q-1)}(x_{i+j-1,j+1},x_{i+j,j+1})$. By $\Gamma_{1,q-1}(x_{i-1,j})=P_{(1,q-1),(1,q-1)}(x_{i-1,j},x_{i,j}')$, we get $(x_{i+j-1,j+1}',x_{i,j}')\in\Gamma_{1,q-1}$. Then $(F_{T\setminus I}(x_{i-1,j}),F_{T\setminus I}(x_{i,j}'))$ is an arc. Since $P_{(1,q-1),(1,q-1)}(x_{i-1,j},x_{i,j}')=\Gamma_{q-1,1}(x_{i,j}')$, from Lemma \ref{mix-(2,q-2) 1} (ii) and Proposition \ref{fenqingkuang} (b), there exists $y_{i,j}\in F_{T\setminus I}(x_{i,j}')$ such that $(x_{i-1,j},y_{i,j})\in\Gamma_{1,q-2}$ and $x_{i+j-1,j+1},x_{i+j-1,j+1}'\in P_{(1,q-1),(1,q-1)}(x_{i-1,j},y_{i,j})$. Without loss of generality, we may assume $y_{i,j}=x_{i,j}'$. The fact $P_{(1,q-2),(1,q-2)}(x_{i-1,j},x_{i+1,j})=\Gamma_{1,q-2}(x_{i-1,j})$ implies $(x_{i,j}',x_{i+1,j})\in\Gamma_{1,q-2}$. Since $x_{i-1,j}\in P_{(q-2,1),(1,q-2)}(x_{i,j},x_{i,j}')$, one gets $(x_{i-1,j}',x_{i,j}')\in\Gamma_{1,q-2}$. By $x_{i,j}'\notin F_{T\setminus I}(x_{i,j})$, from Lemmas \ref{mix-(2,q-2) 1}, \ref{Gamma{1,q-1}=2} and Proposition \ref{fenqingkuang} (ii), one has $(x_{i,j},x_{i,j}')\in\Gamma_{q-1,q-1}$.

Suppose that $V\Gamma\neq\dot{\bigcup}_{i,j}(F_{T\setminus I}(x_{i,j})\cup F_{T\setminus I}(x_{i,j}'))$. Similarly, there exists a vertex $x_{i,j}''$ such that $x_{i,j}''\notin \dot{\bigcup}_{a,b}(F_{T\setminus I}(x_{a,b})\cup F_{T\setminus I}(x_{a,b}'))$ and $x_{i,j}''\in\Gamma_{q-1,q-1}(x_{i,j})\cap\Gamma_{q-1,q-1}(x_{i,j}')$. It follows that $x_{i,j}\in\Gamma_{q-1,q-1}(x_{i,j}')=P_{(q-1,q-1),(q-1,q-1)}(x_{i,j},x_{i,j}')$, a contradiction. Then $\Gamma/F_{T\setminus I}$ is isomorphic to the digraph in Theorem \ref{Main4} (vi).

\textbf{Case 2.} $\partial(F_{T\setminus I}(x_{0,0}),F_{T\setminus I}(x))=2$ for some $x_{0,2}\in\Gamma_{1,q-1}(x_{0,1})$.

Note that $(x_{0,0},x_{0,2})\in\Gamma_{2,q-1}$. Lemma \ref{mix-(2,q-2) 1} (ii) implies $\Gamma_{1,q-1}^{2}=\{\Gamma_{1,q-2},\Gamma_{2,q-1}\}$. Since $p_{(1,q-1),(1,q-1)}^{(1,q-2)}=k_{1,q-1}$, we have $x_{0,2}\in P_{(1,q-1),(1,q-1)}(x_{0,1},x_{1,1})$. By the weakly distance-regularity of $\Gamma$, there exists $x_{i,2}\in\Gamma_{1,q-1}(x_{i,1})\cap\Gamma_{q-1,1}(x_{i+1,1})$ such that $(x_{i,0},x_{i,2})\in\Gamma_{2,q-1}$ and $\partial(F_{T\setminus I}(x_{i,0}),F_{T\setminus I}(x_{i,2}))=2$ for all $i$. The fact that $x_{i-1,0}\in P_{(q-1,2),(1,q-2)}(x_{i-1,2},x_{i,0})$ implies $x_{i,2}\notin P_{(2,q-1),(q-1,2)}(x_{i-1,2},x_{i,0})$ and $(x_{i-1,2},x_{i,2})\in\Gamma_{1,q-2}$. Since $P_{(1,q-2),(1,q-2)}(x_{i,2},x_{i+2,2})=\Gamma_{1,q-2}(x_{i,2})$, one obtains $(x_{i,2},x_{i+2,0})\in\Gamma_{2,q-1}$. The fact that $\partial(F_{T\setminus I}(x_{i,0}),F_{T\setminus I}(x_{i,2}))=2$ implies $\partial(F_{T\setminus I}(x_{i,2}),F_{T\setminus I}(x_{i+2,0}))=2$. By Lemma \ref{Gamma{1,q-1}=2} and Proposition \ref{fenqingkuang} (b), $F_{T\setminus I}(x_{i,j})$ are all distinct with $0\leq i\leq q-2$ and $0\leq j\leq 2$.

Note that there exists a vertex $x_{i,3}\in\Gamma_{1,q-1}(x_{i,2})\cap\Gamma_{q-1,1}(x_{i+1,2})$ such that $(x_{i-1,3},x_{i,3})\in\Gamma_{1,q-2}$, $(x_{i,1},x_{i,3})\in\Gamma_{2,q-1}$ and $\partial(F_{T\setminus I}(x_{i,1}),F_{T\setminus I}(x_{i,3}))=2$ for all $i$. Similarly, $(x_{i,3},x_{i+2,1})\in\Gamma_{2,q-1}$ and $\partial(F_{T\setminus I}(x_{i,3}),F_{T\setminus I}(x_{i+2,1}))=2$. Since $P_{(1,q-1),(1,q-1)}(x_{i,1},x_{i,3})=\Gamma_{1,q-1}(x_{i,1})$ and $p_{(1,q-1),(1,q-1)}^{(1,q-2)}=k_{1,q-1}$, one gets $(x_{i,0},x_{i-1,3}),(x_{i,3},x_{i+2,0})\in\Gamma_{1,q-1}$. By Lemma \ref{Gamma{1,q-1}=2} and Proposition \ref{fenqingkuang} (b) again, $F_{T\setminus I}(x_{i,j})$ are all distinct with $0\leq i\leq q-2$ and $0\leq j\leq 3$.

If $V\Gamma=\dot{\bigcup}_{i,j}F_{T\setminus I}(x_{i,j})$, from \cite[Proposition 9]{YYF16}, then $\Gamma/F_{T\setminus I}$ is isomorphic to the digraph in Theorem \ref{Main4} (v). Now suppose $V\Gamma\neq\dot{\bigcup}_{i,j}F_{T\setminus I}(x_{i,j})$.

Note that there exists a vertex $x_{0,1}'$ such that $\partial(F_{T\setminus I}(x_{0,0}),F_{T\setminus I}(x_{0,1}'))=1$ with $x_{0,1}'\notin\dot{\bigcup}_{a,b}F_{T\setminus I}(x_{a,b})$. Since $p_{(1,q-1),(1,q-1)}^{(1,q-2)}=k_{1,q-1}$, we may assume that $(x_{0,0},x_{0,1}')\in\Gamma_{1,q-1}$. The fact that $p_{(1,q-1),(1,q-1)}^{(2,q-1)}=p_{(1,q-1),(1,q-1)}^{(1,q-2)}=k_{1,q-1}$ implies that $(x_{0,1}',x_{1,0}),(x_{0,1}',x_{0,2})\in\Gamma_{1,q-1}$. If $(x_{0,1}',x_{3,0}),(x_{0,1}',x_{1,1})\in\Gamma_{i,j}$ with $(i,j)\in\{(1,q-2),(2,q-1)\}$, from $\Gamma_{j,i}(x_{0,3})=\Gamma_{j,i}(x_{1,1})=P_{(j,i),(i,j)}(x_{0,3},x_{1,1})$, then $\wz{\partial}_{\Gamma}(x_{0,1},x_{1,1})=\wz{\partial}_{\Gamma}(x_{0,1},x_{0,3})$, a contradiction. Without loss of generality, we assume $(x_{0,1}',x_{1,1})\in\Gamma_{1,q-2}$ and $(x_{0,1}',x_{0,3})\in\Gamma_{2,q-1}$.

Observe that there exists a vertex $x_{i,j}'$ such that  $(x_{i,j}',x_{i+1,j})\in\Gamma_{1,q-2}$ and $\partial(F_{T\setminus I}(x_{i,j}'),F_{T\setminus I}(x_{i+1,j}))=1$ with $x_{i,j}'\notin F_{T\setminus I}(x_{i,j})$. The fact that $p_{(1,q-1),(1,q-1)}^{(1,q-2)}=k_{1,q-1}$ implies that $(x_{i,j}',x_{i,j+1}),(x_{i,3}',x_{i+2,0}),(x_{i,j'}',x_{i+1,j'-1}),(x_{i,0}',x_{i-1,3})\in\Gamma_{1,q-1}$ for all $j\in\{0,1,2\}$ and $j'\in\{1,2,3\}$. In view of $p_{(1,q-1),(1,q-1)}^{(2,q-1)}=k_{1,q-1}$, one obtains $(x_{i,j},x_{i,j+1}'),(x_{i,3},x_{i+2,0}'),(x_{i,j'},x_{i+1,j'-1}'),(x_{i,0},x_{i-1,3}')\in\Gamma_{1,q-1}$ for all $j\in\{0,1,2\}$ and $j'\in\{1,2,3\}$. Then $(x_{i,j}',x_{i,j+1}'),(x_{i,3}',x_{i+2,0}'),(x_{i,j'}',x_{i+1,j'-1}'),(x_{i,0}',x_{i-1,3}')\in\Gamma_{1,q-1}$. By $P_{(1,q-2),(1,q-2)}(x_{i-1,j},x_{i+1,j})=\Gamma_{q-2,1}(x_{i+1,j})$, we have $(x_{i-1,j},x_{i,j}')\in\Gamma_{1,q-2}$ for $j\in\{0,1,2,3\}$. The fact $x_{i+1,j}\in P_{(1,q-2),(1,q-2)}(x_{i,j}',x_{i+2,j})$ implies $(x_{i,j}',x_{i+1,j}')\in\Gamma_{1,q-2}$. By Lemma \ref{Gamma{1,q-1}=2} and Proposition \ref{fenqingkuang} (b), $F_{T\setminus I}(x_{i,j})$ and $F_{T\setminus I}(x_{i,j}')$ are all distinct $0\leq i\leq q-2$ and $0\leq j\leq 3$.

Suppose that there exists a vertex $x_{i,j}''$  such that $x_{i,j}''\notin F_{T\setminus I}(x_{i,j})\cup F_{T\setminus I}(x_{i,j}')$ and $(x_{i,j}'',x_{i+1,j})\in\Gamma_{1,q-2}$ for some $i\in\{0,1,\ldots,q-2\}$ and $j\in\{0,1,2,3\}$. By Lemma \ref{(1,q-1)}, we get $P_{(1,q-2),(1,q-2)}(x_{i-1,j},x_{i+1,j})=\Gamma_{q-2,1}(x_{i+1,j})$, which implies $(x_{i-1,j},x_{i,j}'')\in\Gamma_{1,q-2}$. Since $F_{T\setminus I}(x_{i,j}), F_{T\setminus I}(x_{i,j}')$ and $F_{T\setminus I}(x_{i,j}'')$ are distinct, from Lemma \ref{Gamma{1,q-1}=2} and Proposition \ref{fenqingkuang} (ii), one obtains $(x_{i,j},x_{i,j}'),(x_{i,j}',x_{i,j}''),(x_{i,j}'',x_{i,j})\in\Gamma_{q-1,q-1}$, contrary to $x_{i,j}'\in\Gamma_{q-1,q-1}(x_{i,j})=P_{(q-1,q-1),(q-1,q-1)}(x_{i,j},x_{i,j}')$.

Suppose that $V\Gamma\neq\dot{\bigcup}_{i,j}(F_{T\setminus I}(x_{i,j})\cup F_{T\setminus I}(x_{i,j}'))$. Then there exists a vertex $y$ such that $y\notin\dot{\bigcup}_{i,j}(F_{T\setminus I}(x_{i,j})\cup F_{T\setminus I}(x_{i,j}'))$ and $(x_{0,0},y)\in\Gamma_{1,q-1}$. The fact that $P_{(1,q-1),(1,q-1)}(x_{0,0},x_{0,2})=\Gamma_{1,q-1}(x_{0,0})$ implies that $(y,x_{0,2})\in\Gamma_{1,q-1}$, and so $x_{0,3},x_{1,1}\in\Gamma_{2,q-1}(y)$. Then $x_{1,3}\in\Gamma_{2,q-1}(x_{1,1})=P_{(2,q-1),(q-1,2)}(x_{0,3},x_{1,1})$, a contradiction. Thus, $V\Gamma=\dot{\bigcup}_{i,j}(F_{T\setminus I}(x_{i,j})\cup F_{T\setminus I}(x_{i,j}'))$ and the quotient digraph $\Gamma/F_{T\setminus I}$ isomorphic to the digraph in Theorem \ref{Main4} (vi).$\qed$

Combining Proposition \ref{fenqingkuang} (i) and Propositions \ref{digraph C1 4}, \ref{c2 c3}, \ref{quotient 1}, \ref{quotient 2}, we complete the proof of Theorem \ref{Main4}.

\section*{Acknowledgement}
Y. Yang is supported by the Fundamental Research Funds for the Central Universities (2652017141), K. Wang is supported by NSFC (11671043).

\end{CJK*}

\begin{thebibliography}{00}
\bibitem{ZA99} Z. Arad, E. Fisman and M. Muzychuk, Generalized table algebras, Israel J. Math. 114 (1999) 29--60.

\bibitem{EB84} E. Bannai and T. Ito, Algebraic Combinatorics I: Association Schemes,
Benjamin/Cummings, California, 1984.

\bibitem{AEB98} A.E. Brouwer, A.M. Cohen and A. Neumaier, Distance-Regular Graphs, Springer-Verlag, New York, 1998.

\bibitem{DKT16} E.R. van Dam, J.H. Koolen and H. Tanaka, Distance-regular graphs, Electron. J. Combin. (2016) DS22.

\bibitem{HS04} H. Suzuki, Thin weakly distance-regular digraphs, J. Combin. Theory Ser. B 92 (2004) 69--83.

\bibitem{KSW03} K. Wang and H. Suzuki, Weakly distance-regular digraphs, Discrete Math. 264 (2003) 225--236.

\bibitem{KSW04} K. Wang, Commutative weakly distance-regular digraphs of girth 2, European J. Combin. 25 (2004) 363--375.

\bibitem{YYF16} Y. Yang, B. Lv and K. Wang, Weakly distance-regular digraphs of valency three, \Rmnum{1}, Electron. J.
Combin. 23(2) (2016), Paper 2.12.

\bibitem{YYF18} Y. Yang, B. Lv and K. Wang, Weakly distance-regular digraphs of valency three, \Rmnum{2}, J. Combin. Theory Ser. A 160 (2018) 288--315.

\bibitem{YYF16+} Y. Yang, B. Lv and K. Wang, Quasi-thin weakly distance-regular digraphs, J. Algebraic Combin. 51 (2020) 19--50.

\bibitem{MY15} M. Yoshikawa, On regular association schemes of order $pq$, Discrete Math. 338 (2015) 111--113.

\bibitem{MY16} M. Yoshikawa, On association schemes of finite exponent, European J. Combin. 51 (2016) 433--442.

\bibitem{MY18} M. Yoshikawa, On higher indicators of regular association schemes, Discrete Math. 341 (2018) 2028--2034.

\bibitem{MY20} M. Yoshikawa,  On closed subsets generated by a regular relation of an association scheme, Discrete Math. 343 (2020) 111706.

\bibitem{PHZ96} P.H. Zieschang, An Algebraic Approach to Assoication Schemes, in: Lecture Notes in Mathematics, Vol.1628, Springer, Berlin, Heidelberg, 1996.

\bibitem{PHZ05} P.H. Zieschang, Theory of Association Schemes, Springer Monograph in Mathematics, Springer, Berlin, 2005.
\end{thebibliography}
\end{document}